\documentclass[a4paper]{article}
\usepackage{RR}
\RRNo{6515}
\usepackage{hyperref}
\usepackage{graphicx}
\usepackage{algorithm}
\usepackage{algorithmic}
\usepackage{url}
\newtheorem{definition}{Definition}
\newcounter{step}
\newenvironment{code}{\begin{list}{{\em %
  Step~\arabic{step}:}}{\usecounter{step}\setlength{\labelsep}{1em}%
  \setlength{\leftmargin}{4em}\setlength{\itemsep}{0em}%
  \setlength{\rightmargin}{0em}\setlength{\labelwidth}{3em}}}{\end{list}}
\RRdate{April 2008}

\RRauthor{Arnaud~Liefooghe \and Laetitia~Jourdan \and El-Ghazali~Talbi}
\authorhead{Liefooghe et al.}
\RRtitle{Metaheuristiques et leur hybridation pour r\'esoudre le Ring Star Problem bi-objectif:\\une \'etude comparative}
\RRetitle{Metaheuristics and Their Hybridization to\\Solve the Bi-objective Ring Star Problem:\\a Comparative Study}
\titlehead{Metaheuristics and Hybridization for the Bi-objective Ring Star Problem}
\RRnote{Arnaud.Liefooghe@lifl.fr, Laetitia.Jourdan@lifl.fr, El-Ghazali.Talbi@lifl.fr}
\RRresume{
Ce document pr\'esente et exp\'erimente diff\'erentes approches pour la r\'esolution d'un nouveau probl\`eme de routage bi-objectif appel\'e le ring star problem.
Il consiste \`a localiser un cycle simple dans un sous-ensemble de noeuds d'un graphe en optimisant deux types de co\^ut. 
Le premier objectif est la minimisation d'un co\^ut d'anneau qui est li\'e \`a la longueur du cycle. 
Le deuxi\`eme est la minimisation d'un co\^ut d'assignement depuis les noeuds non visit\'es vers les noeuds visit\'es. 
En d\'epit de sa formulation bi-objectif \'evidente, ce probl\`eme a toujours \'et\'e \'etudi\'e de fa\c{c}on mono-objectif.
Pour attaquer le ring star problem bi-objectif, nous \'etudions d'abord diff\'erentes m\'ethodes de recherche.
Ensuite, nous proposons deux strat\'egies coop\'eratives qui combinent deux m\'etaheuristiques multi-objectif :
un algorithme \'evolutionnaire \'elitiste et une recherche locale \`a base de population. 
Nous appliquons ces approches hybrides à des instances de test r\'eput\'ees 
et d\'emontrons leur efficacit\'e par rapport aux algorithmes non-hybrides et \`a des m\'ethodes classiques de la litt\'erature.
}
\RRabstract{
This paper presents and experiments approaches to solve a new bi-objective routing problem called the ring star problem.
It consists of locating a simple cycle through a subset of nodes of a graph while optimizing two kinds of cost. 
The first objective is the minimization of a ring cost that is related to the length of the cycle.
The second one is the minimization of an assignment cost from non-visited nodes to visited ones.
In spite of its obvious bi-objective formulation, this problem has always been investigated in a single-objective way.
To tackle the bi-objective ring star problem, we first investigate different stand-alone search methods.
Then, we propose two cooperative strategies that combines two multiple objective metaheuristics: an elitist evolutionary algorithm and a population-based local search.
We apply this new hybrid approaches to well-known benchmark test instances 
and demonstrate their effectiveness in comparison to non-hybrid algorithms and to state-of-the-art methods.
}
\RRmotcle{Optimisation combinatoire multi-objectif, algorithmes \'evolution\-naires, exploration locale, approches hybrides, routage.}
\RRkeyword{Multi-objective combinatorial optimization, evolutionary algorithms, local exploration, hybrid approaches, routing.}
\RRprojet{Dolphin}  
\RRtheme{\THNum} 
\RCLille 

\begin{document}
\makeRR   

\section{Introduction}

First introduced by Labb\'e et al.~\cite{LL+:99}, the \emph{Ring Star Problem} (RSP), also known as \emph{Median Cycle Problem} (MCP), is a recent routing problem.
It aims to locate a simple cycle (the ring) containing a subset of nodes of a graph (including the depot), 
while nodes that do not belong to the ring must be connected to a visited one.
Two costs are considered here, what naturally leads to a bi-criteria formulation.
The first one is the ring cost associated to the ring itself 
and the second one is the assignment cost associated to the arcs directed from non-visited nodes to visited ones.
Even so, the RSP has always been investigated in a single-objective way; either where both criterion are combined, or where a criteria is treated as a constraint.
Thus, despite its numerous industrial interests, this is the first time that such a problem is regarded as a multi-objective one, 
probably due to its complexity.
Indeed, it is highly combinatorial as, once is decided which nodes have to be visited or not, a traveling salesman problem is still to be solved.

Recently, we studied the effectiveness of existing metaheuristics for solving the bi-objective RSP~\cite{LJ+:08}.
In this paper, we provide a deeper analysis of this work, and propose new stand-alone and hybrid approaches.
As an initial step, we investigate four metaheuristics to approximate the set of efficient solutions for the problem under consideration.
First, IBMOLS, a population-based local search recently proposed in~\cite{BB:07}, is designed with a variable neighborhood.
Next, we present steady-state variations of two well-known multi-objective search methods, namely IBEA~\cite{ZK:04} and NSGA-II~\cite{DA+:02}.
Last, we propose a new simple elitist evolutionary algorithm, SEEA, able to solve any kind of multi-objective combinatorial optimization problems.
We compare all these metaheuristics to each other on state-of-the-art benchmarks test instances. 
As a second step, we propose two new cooperative strategies between the local search method and SEEA,
based on a periodic and on an auto-adaptive scheme, that try to benefit of the advantages of each method it is compound of.
The performance of these hybrid metaheuristics is shown in comparison to non-hybrid approaches.

The paper is organized as follows.
Section~\ref{sec:rsp} is devoted to the bi-objective RSP.
Then, four metaheuristics for multi-objective optimization are introduced in Section~\ref{sec:meta}.
A general presentation of these methods and their application to the RSP is followed by a parameter analysis and a comparative case study.
In Section~\ref{sec:hybrid}, we propose new cooperation schemes to solve multi-objective optimization problems.
We experiment the resulting search methods on the bi-objective RSP and we compare the obtained computational results to the previous ones. 
At last, conclusions and perspectives are drawn in the last section.

\section{The Bi-objective Ring Star Problem}
\label{sec:rsp}

In this section, we first present some basic concepts, notation and definitions related to multi-objective optimization.
Next, we provide a formulation of the Ring Star Problem (RSP) as a bi-objective problem.
Finally, we survey the literature relating to the RSP and discuss its industrial concerns.

\subsection{Multi-objective Optimization}

A \emph{Multi-objective Optimization Problem} (MOP) aims to optimize a set of $n \geq 2$ objective functions $f_1 , f_2, \dots , f_n$ simultaneously.
Each objective function can be either minimized or maximized.
Let~$X$ denote the set of feasible solutions in the \emph{decision space}, and $Z$ the set of feasible points in the \emph{objective space}.
A decision vector $x = ( x_1 , x_2 , \ldots , x_k )$ is represented by a vector of $k$~decision variables.
In the case of a \emph{Multi-objective Combinatorial Optimization Problem} (MCOP), a decision vector~$x \in X$ has a finite set of possible values.
To each decision vector $x \in X$ is assigned exactly one objective vector~$z \in Z$ on the  basis of a vector function 
$f : X \rightarrow Z$ with $z = f(x) = (f_1(x) , f_2(x), \dots , f_n(x))$.
We will assume, throughout the paper, that objective values are normalized.
To achieve this, the minimum and the maximum value of each objective function are used in order to adaptively replace each objective function by its corresponding normalized function, so that its values lie in the interval~$[0,1]$.
Without loss of generality, we assume that $Z \subseteq \Re^n$ and that all $n$ objective functions have to be minimized.
Therefore, a MOP can be formulated as follows:
\begin{equation}
(MOP) = 
\left\{
\begin{array}{rl}
\mbox{`min'}      & f(x) = ( f_1(x) , f_2(x) , \ldots , f_n(x) ) \\
\mbox{subject to} & x \in X 
\end{array}
\right.\
\end{equation}
\begin{definition}
An objective vector $z \in Z$ \emph{weakly dominates} another objective vector $z' \in Z$ if and only if 
$\forall i \in \{1,2,\dots,n\}$, $z_i \leq z_i'$.
\end{definition}
\begin{definition}
An objective vector $z \in Z$ \emph{dominates}%
\footnote{We will also say that a decision vector $x \in X$ \emph{dominates} a decision vector $x' \in X$ if~$f(x)$ dominates $f(x')$.}
another objective vector $z' \in Z$ if and only if 
$\forall i \in \{1,2,\dots,n\}$, $z_i \leq z_i'$ and $\exists j \in \{1,2,\dots,n\} $ such as $z_j < z_j'$. 
\end{definition}
\begin{definition}
An objective vector $z \in Z$ is \emph{non-dominated} if and only if there does not exist another objective vector $z' \in Z$ such that $z'$ dominates $z$.
\end{definition}
A solution $x \in X$ is said to be \emph{efficient} (or \emph{Pareto optimal}, \emph{non-dominated}) if its mapping in the objective space results in a non-dominated point.
The set of all efficient solutions is the \emph{efficient} (or \emph{Pareto optimal}) \emph{set}, denoted by~$X_E$.
The set of all non-dominated vectors is the \emph{non-dominated front} (or the \emph{trade-off surface}), denoted by~$Z_N$.
A common approach in solving MOPs is to find or to approximate the minimal set of efficient solutions,
{\itshape i.e.} a solution $x \in X_E$ for each non-dominated vector~$z \in Z_N$ such as $f(x) = z$ 
(in case there exists multiple solutions mapping to the same non-dominated point).
But, generating the entire set of Pareto optimal solutions is usually infeasible due to the complexity of the underlying problem or to the large number of optima.
Therefore, the overall goal is often to identify a good approximation of it. 
Population-based metaheuristics are commonly used to this end as they naturally find multiple and well-spread non-dominated solutions in a single simulation run.
A reasonable basic introduction to multi-objective optimization can be found in~\cite{Deb:05}.
But, the interested reader could refer to~\cite{CVL:02,Deb:01} for more details about evolutionary multi-objective optimization 
and to~\cite{EG:04} for more details about multi-objective combinatorial optimization.

\subsection{Problem Definition}

\begin{figure}[htbp]
\centering
\includegraphics[width=3.2in]{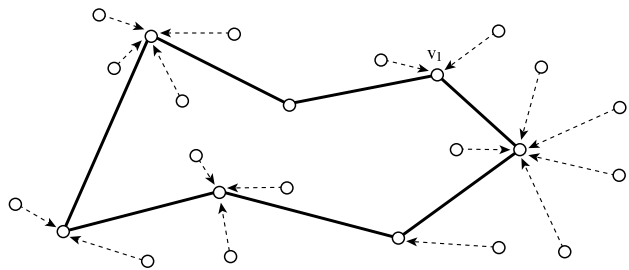}
\caption{An example of solution for the ring star problem.}
\label{fig:rsp}
\end{figure}

The \emph{Ring Star Problem} (RSP) can be described as follows.
Let $G = ( V , E , A )$ be a complete mixed graph where 
$V = \{ v_1 , v_2 , \dots , v_n \}$ is a set of vertexes, 
$E = \{ [v_i , v_j] | v_i , v_j \in V , i < j \}$ is a set of edges, and
$A = \{ (v_i , v_j) | v_i , v_j \in V \}$ is a set of arcs.
Vertex $v_1$ is the depot.
To each edge $[ v_i , v_j ]$ we assign a non-negative \emph{ring cost}~$c_{ij}$, 
and to each arc $( v_i , v_j )$ we assign a non-negative \emph{assignment cost}~$d_{ij}$.
The RSP consists of locating a simple cycle through a subset of nodes $V' \subset V$ (with $v_1 \in V'$) while 
($i$) minimizing the sum of the ring costs related to all edges that belong to the cycle, and
($ii$) minimizing the sum of the assignment costs of arcs directed from every non-visited node to a visited one so that the associated cost is minimum.
An example of solution is given in Figure~\ref{fig:rsp}, where solid lines represent edges that belong to the ring and dashed lines represent arcs of the assignments.

\noindent
The first objective is called the \emph{ring cost} and is defined as:
\begin{equation}
	\sum_{[ v_i ,v_j ] \in E}{c_{ij} x_{ij}}
\end{equation}
where $x_{ij}$ is a binary variable equal to 1 if and only if the edge $[ v_i , v_j ]$ belongs to the cycle.
The second objective, the \emph{assignment cost}, can be computed as follows:
\begin{equation}
	\sum_{v_i \in V \setminus V'}{\min_{v_j \in V'}{d_{ij}}}
\end{equation}
Let us remark that these two objectives are comparable only if we assume that the ring cost and the assignment cost are commensurate one to another, 
what is rarely the case in practice.
Furthermore, the fact of privileging a cost compared to the other is closely related to the decision-maker preferences.
However, 
the RSP is an NP-hard combinatorial problem since the particular case of visiting the whole set of nodes is equivalent to a traditional traveling salesman problem.

\subsection{Related Works}

The RSP belongs to the class of location-allocation problems aiming at locating structures in a graph (see~\cite{LLR:98} for a review).
It has initially been formulated by Labb\'e et al.~\cite{LL+:99} in two different ways.
In the first formulation (denoted MCP1 in~\cite{LL+:99} and more often called `ring star problem'), the sum of both objectives is to be optimized.
In the second formulation (MCP2~\cite{LL+:99}, usually called `median cycle problem'), the ring cost is to be minimized, 
while the assignment cost is bounded by a prefixed value.
Even if it was not explicitly noticed by the original authors, 
these two kinds of formulation are commonly employed to convert a multi-objective problem into a single-objective one by using two well-known scalar approaches; 
respectively an \emph{aggregation method} and an \emph{epsilon-constraint method}%
\footnote{The reader is referred to~\cite{Mie:99} for more details about scalar methods for solving MOPs.}.
The first formulation of the problem has been more widely studied in~\cite{LL+:04}. 
The authors used a branch-and-bound method and successfully solved TSPLIB and randomly generated instances involving up to $200$ nodes in less than two hours.
In~\cite{LL+:05}, the same authors investigated the second formulation of the problem using a very close method.
Finally, one or both versions of the problem have been heuristically solved by 
a variable neighborhood tabu search~\cite{MMR:03}, 
an evolutionary algorithm~\cite{RBL:04},
a multi-start greedy add heuristic~\cite{RBL:04}, 
and a variable neighborhood tabu search hybridized with a greedy randomized adaptive search procedure~\cite{DS+:06}.

A very similar problem, namely the \emph{capacitated m-ring star problem}, has been introduced in~\cite{BDS:07}.
It extends the first formulation by bounding the number of nodes that can be assigned to a ring using a prefixed parameter (the ring capacity).
This problem has been exactly solved by a branch-and-bound method and experiments have been performed on random and real-world instances.
Note that the same problem has also been tackled in~\cite{MN+:07} by a hybrid metaheuristic.
Another variant of the RSP is the \emph{Steiner Ring Star Problem}, where only optional nodes can be visited~\cite{LCS:98}.
In addition, Beasley and Nascimento~\cite{BN:96} introduced a \emph{Single Vehicle Routing-Allocation Problem}~(SVRAP)
in which a non-visited node can be either assigned to the cycle or isolated.
In this case, a corresponding isolated cost has to be minimized as well.
Vogt et al.~\cite{VPB:07} proposed a tabu search to solve a single-objective problem resulting of this formulation.

As reflected in the survey of Jozefowiez et al.~\cite{JST:08}, an increasing number of multi-objective routing problems appeared in the literature in recent years.
However, in spite of its numerous industrial applications (see subsection~\ref{sec:rsp:industrial}), 
the RSP has never been explicitly investigated in a multi-objective fashion.
Nevertheless, Current and Schilling~\cite{CS:94} defined two multiple criteria variants of a very similar problem: 
the \emph{median tour problem} and the \emph{maximal covering tour problem}.
In both, one criterion is the minimization of the total length of the tour, while another one is the maximization of the access to the tour for non-visited nodes.
To tackle these problems, the authors used a kind of lexicographic method. 
This approach consists of assigning a priority value to each objective function, and in solving the problem one objective after another. 
Also, Dorner et al.~\cite{DFG:07} recently formulated a tri-objective optimization problem of tour planning for mobile health care facilities in Senegal,
closely related to the SVRAP~\cite{BN:96} introduced earlier.
A mobile facility has to visit a subset of nodes.
Non-visited nodes are then assigned to their closest tour stop or are regarded as unable to reach a tour stop (within a predefined maximum distance).
The objectives are 
($i$)~to minimize the ratio between medical working time and total working time, 
($ii$)~to minimize the average distance to the nearest tour stops and 
($iii$)~to maximize a coverage criterion.
To do so, a Pareto ant colony optimization algorithm as well as two genetic algorithms (namely VEGA~\cite{Sch:85} and MOGA~\cite{FF:93}) 
were designed to solve real-world instances.

\subsection{Industrial Interests}
\label{sec:rsp:industrial}

The RSP has a wide range of industrial interests, especially in telecommunications.
Of course, practical applications may contain additional constraints.
For instance, as noticed by Labb\'e et al.~\cite{LL+:04,LL+:05}, 
a ring-based network is designed to interconnect a set of hubs in the case of digital data services~\cite{XCG:99}.
Some concentrators are installed on a subset of locations and are interconnected on a ring network (the Internet) 
while the remaining locations are assigned to this concentrators (the Intranet).
Very close problems arise in rapid transit systems planning~\cite{LL+:04} or while designing optical networks~\cite{BDS:07}.
Also, other kinds of applications appears in the postal collection or delivery routes design, 
where the distance between a customer and a collection point have to be reasonable.
For instance, post-box location while taking both the collection cost and the user inconvenience into account has been studied in~\cite{LL:86}.
Besides, other applications closely related to the RSP are the location of circular shaped transportation infrastructure (such as metro lines or motorways),
the location of recyclable garbage collection bins, and school bus routing.
Finally, the routing of essential health care services, already investigated in~\cite{AS:92,DFG:07} among other authors, 
consists of a mobile clinic servicing an area without being able to visit every population nodes.
Then, unvisited ones have to reach the nearest tour stop by their own to be medically treated.

\section{Metaheuristics to Solve the Bi-objective Ring Star Problem}
\label{sec:meta}

Four metaheuristics are proposed to tackle the bi-objective RSP: 
a variable neighborhood iterative population-based Local Search~(LS) and three Evolutionary Algorithms (EAs).
These algorithms are respectively steady-state variations of IBMOLS~\cite{BB:07}, IBEA~\cite{ZK:04} and NSGA-II~\cite{DA+:02}, 
as well as a new search method, SEEA, introduced in this paper for the first time. 
IBMOLS and IBEA are both recent indicator-based metaheuristics, whereas NSGA-II is one of the most often used Pareto-based resolution approach.
In this section, RSP-specific components are described after we presented the main characteristics of the LS and of the EAs.

\subsection{A Multi-objective Local Search}

Since they are easily adaptable to the multi-objective context, many of the search algorithms proposed to tackle MOPs are EAs.
However, LS algorithms are known to be effective metaheuristics for solving real-world applications~\cite{BK:05,GK:03}.
Several neighborhood search methods to solve MOPs have been proposed in the literature.
Most of them are based on a set of aggregations or scalarizations of the objective functions, what is the case, for instance, in~\cite{GMF:97,Han:97,IM:98,UT+:99}.
Pareto-based LS algorithms, like the one proposed by Paquete et al.~\cite{PCS:04}, are more rare.
An Indicator-Based Multi-Objective Local Search (IBMOLS for short) has recently been presented in~\cite{BB:07}.
IBMOLS is a generic population-based multi-objective LS dealing with a fixed population size.
This allows to obtain a set of efficient solutions in a single run without specifying any mechanism to control the number of solutions during the search process.
Moreover, IBMOLS presents an alternative to aggregation- and Pareto-based multi-objective metaheuristics.
Indeed, as proposed in~\cite{ZK:04}, it is assumed that the optimization goal is given in terms of a binary quality indicator $I$~\cite{ZT+:03} 
which can be regarded as an extension of the Pareto dominance relation.
A value~$I(A,B)$ quantifies the difference in quality between two approximated efficient sets~$A$ and $B$. 
So, if $R$ denotes a reference set (which can be the Pareto-optimal set~$X_E$ or any other set), the overall optimization goal can be formulated as:
\begin{equation}
  \mbox{arg min}_{A \in \Omega}\ I(A,R) \ ,
\end{equation}
where $\Omega$ denotes the space of all efficient set approximations. 
As noted in~\cite{ZK:04}, $R$ does not have to be known, it is just required in the formalization of the optimization goal.
Since $R$ is fixed, $I$ actually represents a unary function that assigns a real number reflecting the quality of each approximation set according to the optimization goal.
If $I$ is dominance preserving~\cite{ZK:04}, $I(A,R)$ is minimum for $A=R$.
One of the main advantages of indicator-based optimization is that no specific diversity preservation mechanism is generally required, 
according to the indicator being used.

The IBMOLS algorithm maintains a population~$P$ of size~$N$.
Then, it generates the neighborhood of a solution contained in~$P$ until a good solution is found 
({\itshape i.e.} one that is better than at least one solution of $P$ in terms of the indicator being used). 
By iterating this simple principle to every solution of $P$, we obtain a local search step. 
The whole local search procedure stops when the archive of potentially efficient solutions has not received any new item during a complete local search step. 
Moreover, as local search methods are usually performed in an iterative way, a population re-initialization scheme has to be designed after each local search.
As noticed in~\cite{BB:07}, several strategies can be used within an Iterative IBMOLS (I-IBMOLS).
Solutions can be re-initialized randomly, and crossover or random noise can be applied to solutions belonging to the efficient set approximation.
The interested reader could refer to~\cite{BB:07} for more details about IBMOLS and I-IBMOLS.

A beneficial feature of this LS is the low number of parameters that are required.
In addition to the population size, 
the binary quality indicator to be used and the population re-initialization strategy (between each local search) are the only two other problem-independent parameters.
%
Indeed, several quality indicators can be used within IBMOLS.
Some examples can be found in~\cite{HJ:98,KC:02,ZT+:03}.
One of them is the binary additive $\epsilon$-indicator ($I_{\epsilon+}$)~\cite{ZK:04,ZT+:03}, 
inspired by the concept of $\epsilon$-dominance introduced by Laumanns et al.~\cite{LT+:02}.
This indicator is particularly well-adapted to indicator-based search and seems to be efficient on different kinds of problems (see, for instance, \cite{BB:07,ZK:04}).
It is capable of obtaining both a well-converged and a well-diversified Pareto front approximation.
$I_{\epsilon+}$ computes the minimum value by which a solution~$x \in X$ can be translated in the objective space to weakly~dominate another solution~$x' \in X$.
For a minimization problem, it is defined as follows:
\begin{equation}
I_{\epsilon +} (x, x') = \max_{i \in \{1,\dots,n\}} (f_i(x) - f_i(x')) \ .
\label{eq:eps}
\end{equation}
Furthermore, to evaluate the quality of a solution $x \in X$ according to a whole population~$P$ and a binary quality indicator~$I$, 
and then to compute the fitness value of $x$, different approaches exist.
As proposed in~\cite{ZK:04}, 
we will here consider an additive technique that amplifies the influence of solutions mapping to dominating points over solutions mapping to dominated ones,
what can be outlined as follows:
\begin{equation}
I(x,P\setminus \{x\})=\sum_{x^\star\in P\setminus \{x\}}-e^{-I(x^\star,x)/\kappa} \ ,
\label{eq:fitness}
\end{equation}
where $\kappa > 0$ is a scaling factor.
However, the initial experiments were not satisfactory because the algorithm was not able to find the extreme points of the trade-off surface.
As pointed out in~\cite{HS+:07}, this is known to be one of the drawbacks of the $\epsilon$-dominance relation, apparently due to the high convexity of the front.
Indeed, the authors illustrate that one of the limitation of the $\epsilon$-dominance is that the extreme points of the Pareto front are usually lost.
To tackle this problem, we add a condition preventing the deletion of solutions corresponding to the extreme non-dominated vectors during the replacement step of IBMOLS.
%
Besides, the population re-initialization scheme used between each local search procedure is based on random noise, 
such as in the basic simulated annealing algorithm~\cite{BK:05}.
This noise consists of multiple mutations applied to $N$ different randomly chosen solutions contained in the archive of potentially efficient solutions.
If the size of the archive is less than $N$, the population is filled with random solutions.

\subsection{Multi-objective Evolutionary Algorithms}

The first two multi-objective EAs designed to tackle the bi-objective RSP 
are variations of two state-of-the-art search methods, namely IBEA~\cite{ZK:04} and NSGA-II~\cite{DA+:02}.
Some minor modifications have been carried out to improve the algorithms for the particular case of the addressed problem, 
for which a preliminary investigation revealed that the set of non-dominated points was, in general, very large.
Finally, a generic approach to solve MOPs, called SEEA, is proposed in this paper for the first time and is presented in details.

\subsubsection{IBEA}
Introduced by Zitzler and K\"unzli~\cite{ZK:04}, the \emph{Indicator-Based Evolutionary Algorithm} (IBEA) is, like IBMOLS, an indicator-based metaheuristic.
The fitness assignment scheme of this EA is based on a pairwise comparison of solutions contained in a population by using a binary quality indicator.
As within IBMOLS, no diversity preservation technique is required, according to the indicator being used.
The selection scheme for reproduction is a binary tournament between randomly chosen individuals.
The replacement strategy is an environmental one that consists of deleting, one-by-one, the worst individuals, 
and in updating the fitness values of the remaining solutions each time there is a deletion; this is continued until the required population size is reached.
Moreover, an archive stores solutions mapping to potentially non-dominated points in order to prevent their loss during the stochastic search process. 
However, in our case, and in contrast to the IBEA defined in~\cite{ZK:04}, this archive is updated at each generation since the beginning of the EA, 
so that the output size is not necessarily less than or equal to the population size.
Just like for the IBMOLS algorithm, the indicator used within IBEA is the additive $\epsilon$-indicator; 
and the same mechanism has been used to prevent the loss of the extreme points on the trade-off surface.
More information about this algorithm can be found in~\cite{ZK:04}.

\subsubsection{NSGA-II}
At each generation of NSGA-II (the \emph{Non-dominated Sorting Genetic Algorithm II} introduced by Deb et al. in~\cite{DA+:02}), 
the solutions contained in the population are ranked into several classes.
Individuals mapping to vectors from the first front all belong to the best efficient set; 
individuals mapping to vectors from the second front all belong to the second best efficient set; and so on.
Two values are computed for every solution of the population.
The first one corresponds to the \emph{rank} the corresponding solution belongs to, and represents the quality of the solution in terms of convergence.
The second one, the \emph{crowding distance}, consists of estimating the density of solutions surrounding a particular point of the objective space, 
and represents the quality of the solution in terms of diversity.
A solution is said to be better than another if it has a best rank value, or in case of equality, if it has the best crowding distance.
The selection strategy is a deterministic tournament between two random solutions.
At the replacement step, only the best individuals survive, with respect to the population size.
Likewise, an external population is added to the steady-state NSGA-II in order to store every potentially efficient solution found during the search.
The reader is referred to~\cite{DA+:02} for more details about NSGA-II.

\subsubsection{SEEA}
\label{sec:seea}
If evaluating a solution in the objective space is not too much time consuming (what is the case for our problem), 
computing fitness values and diversity information are generally the most computationally expensive steps of a multi-objective EA.
Based on this observation, we propose here a simple search method for which none of these phases is required.
The resulting EA, called \emph{Simple Elitist Evolutionary Algorithm} (SEEA for short), is detailed in Algorithm~\ref{algo:SEEA}.
An archive of potentially efficient solutions is updated at each generation, 
and the individuals contained in the main population are generated by applying variation operators to randomly chosen archive members.
The replacement step is a generational one, {\itshape i.e.} the parent population is replaced by the offspring one.
Note that the initial population can, for instance, be filled with random solutions.
Then, as proposed in~\cite{ZT:99} among other authors, 
the archive is not only used as an external storage, but it is integrated into the optimization process during the selection phase of the EA, what is called \emph{elitism}.
The use of elitism is an important issue within evolutionary multi-objective optimization~\cite{LZT:01} 
and SEEA is in somehow related to other elitist multi-objective EAs such as PAES~\cite{KC:00}, PESA~\cite{CKO:00} or SEAMO~\cite{Val:02}.
But, contrary to other approaches, 
no strategy to preserve diversity or to manage the size of the archive is involved here, as solutions are selected randomly and the archive is unbounded.
The biggest advantage of this EA is that the population (or the population size if solutions are randomly initialized) is the only problem-independent parameter.
If non-dominated solutions are relatively close to each other in the decision space and if the archive is not too small compared to the main population, 
we believe that SEEA may convergence to a good approximation of the efficient set in a very short run-time.
However, in some case, this method may prematurely converge or may appear inefficient if promising solutions are far from each other.
\begin{algorithm}
\footnotesize
\caption{Simple Elitist Evolutionary Algorithm (SEEA)}
\label{algo:SEEA}
\begin{tabbing}
\makebox[4em]{}  \=  \makebox[3em]{}  \=  \kill                          \\
{\em Input:}     \>  $P$              \> Initial population              \\
{\em Output:}    \>  $A$              \> Efficient set approximation
\end{tabbing}
\begin{code}
\item{\bfseries Initialization.} $A \leftarrow$ non-dominated individuals of $P$; $N \leftarrow |P|$; $P' \leftarrow \emptyset$.
\item{\bfseries Selection.} Repeat until $|P'|=N$: randomly select an individual from $A$ and add it to the offspring population~$P'$.
\item{\bfseries Variation.} Apply crossover and mutation operators to individuals of the offspring population~$P'$.
\item{\bfseries Replacement.} $P \leftarrow P'$; $P' \leftarrow \emptyset$.
\item{\bfseries Elitism.} $A \leftarrow$ non-dominated individuals of $A \cup P$.
\item{\bfseries Termination.} If a stopping criteria is satisfied return $A$, else go to Step $2$.
\end{code}
\end{algorithm}

\subsection{Application to the Ring Star Problem}

This section provides the problem-specific steps of the metaheuristics introduced earlier.
Components designed for the particular case of the bi-objective RSP, such as the encoding mechanism, 
the population initialization as well as the neighborhood, mutation and crossover operators, are described below.

\subsubsection{Solution Encoding}
The encoding mechanism used to represent a RSP solution, for both the LS and the EAs, is based on the random keys concept proposed by Bean~\cite{Bea:94}.
This implementation has already been successfully applied for a single-objective version of the RSP in~\cite{RBL:04}.
To each node $v_i$ belonging to the ring we assign exactly one \emph{random key} $x_i \in [0,1[$, where $x_1 = 0$.
A special value is assigned to unvisited nodes.
Thus, the ring route associated to a solution corresponds to the nodes ordered according to their random keys in an increasing way;
{\itshape i.e.} if $x_i < x_j$, then $v_j$ comes after $v_i$.
As an example, a possible representation for the cycle~($v_1$,$v_7$,$v_4$,$v_9$,$v_2$,$v_6$) is given in Figure~\ref{fig:randomkey}.
vertexes~$v_3$,$v_5$,$v_8$ and $v_{10}$ are assigned to a visited node in such a way that the associated assignment cost is minimum.
\begin{figure}[htb]
\footnotesize
\centering
\setlength{\tabcolsep}{2.5pt}
\begin{tabular}{p{0.8in}|p{0.2in}|p{0.2in}|p{0.2in}|p{0.2in}|p{0.2in}|p{0.2in}|p{0.2in}|p{0.2in}|p{0.2in}|p{0.2in}p{0in}}
	Vertex				&	\centering $v_1$	& \centering $v_2$	&	\centering $v_3$	&	\centering $v_4$	&	\centering $v_5$
								&	\centering $v_6$	&	\centering $v_7$	& \centering $v_8$	&	\centering $v_9$	&	\centering $v_{10}$	&
								\\
	\hline
	Random key		&	\centering 0			& \centering 0.7		& \centering - 			&	\centering 0.3		&	\centering -
								&	\centering 0.8		&	\centering 0.2		&	\centering -			&	\centering 0.5		&	\centering -	 &
								\\
\end{tabular}
\caption{A RSP solution represented by random keys.}
\label{fig:randomkey}
\end{figure}

\subsubsection{Population Initialization}
For each optimization method, the initial population has been generated randomly.
Each node has a probability $p=0.5$ that it will be visited or not, and to each visited vertex we associate a key value randomly generated between $0$ and $1$.

\subsubsection{Neighborhood and Mutation Operators}
As the RSP is both a routing problem and an assignment problem, different neighborhood and mutation operators have to be designed.
Here, we consider the following:
\begin{itemize}
	\item \emph{insert operator}:	adds an unvisited node $v_i$ in the cycle, the position where to insert $v_i$ is chosen in order to minimize the ring cost
	\item \emph{remove operator}:	removes a vertex $v_j$ of the ring
	\item \emph{2-opt operator}:	applies a 2-opt operator between two nodes of the cycle~$v_i$ and $v_j$, 
																{\itshape i.e.} the sequence of visited nodes located between~$v_i$ and $v_j$ is reversed.
\end{itemize}
For the LS, the neighbors of a solution are randomly explored, without considering any order between these three operators; and each neighbor is at most visited once.
Moreover, note that it is not necessary to completely re-evaluate a solution each time a neighborhood operator is applied.
Thus, after an \emph{insert} neighborhood operator, we just have to re-assign unvisited nodes in order to minimize the assignment cost.
After a \emph{remove} neighborhood operator, we just have to re-assign the nodes that were previously assigned to the one that has been removed.
And, after a \emph{2-opt} neighborhood operator, we just have to recompute the ring cost, the assignment cost being unchanged.
In the case of mutations, the operators are applied to randomly chosen vertexes.

\subsubsection{Crossover Operator}
The crossover operator is a quadratic crossover closely related to the one proposed in~\cite{RBL:04}.
Two randomly selected solutions $s_1$ and $s_2$ are divided in a particular position.
Then, the first part of $s_1$ is combined with the second part of~$s_2$ to build a first offspring, 
and the first part of~$s_2$ is combined with the second part of $s_1$ to build a second offspring.
Every node retains its random key so that it enables an easy reconstruction of the new individuals.
Figure~\ref{fig:crossover} illustrates a recombination between two solutions~($v_1$,$v_7$,$v_4$,$v_9$,$v_2$,$v_6$) and~($v_1$,$v_8$,$v_4$,$v_3$,$v_5$) 
after vertex~$v_6$, what gives rise to a couple of new solutions~($v_1$,$v_8$,$v_4$,$v_2$,$v_6$) and~($v_1$,$v_7$,$v_9$,$v_4$,$v_3$,$v_5$).
Thanks to the random keys encoding mechanism, solutions having a different ring size can easily be recombined, 
even if the initial ring structures are generally broken in the offspring solutions.
\begin{figure}[htb]
\footnotesize
\centering
\begin{tabular}{p{0.7in}l}
	\raggedright Parent 1
&
	\setlength{\tabcolsep}{2.5pt}
	\begin{tabular}{p{0.8in}|p{0.2in}|p{0.2in}|p{0.2in}|p{0.2in}|p{0.2in}|p{0.2in}||p{0.2in}|p{0.2in}|p{0.2in}|p{0.2in}p{0in}}
	Vertex				&	\centering $v_1$	& \centering $v_2$	&	\centering $v_3$	&	\centering $v_4$	&	\centering $v_5$
								&	\centering $v_6$	&	\centering $v_7$	& \centering $v_8$	&	\centering $v_9$	&	\centering $v_{10}$	&
								\\
	\hline
	Random key		&	\centering 0			& \centering 0.7		& \centering - 			&	\centering 0.3		&	\centering -
								&	\centering 0.8		&	\centering 0.2		&	\centering -			&	\centering 0.5		&	\centering -	 &
								\\
	\end{tabular}
\\
\\
	Parent 2
&
	\setlength{\tabcolsep}{2.5pt}
	\begin{tabular}{p{0.8in}|p{0.2in}|p{0.2in}|p{0.2in}|p{0.2in}|p{0.2in}|p{0.2in}||p{0.2in}|p{0.2in}|p{0.2in}|p{0.2in}p{0in}}
	Vertex				&	\centering $v_1$	& \centering $v_2$	&	\centering $v_3$	&	\centering $v_4$	&	\centering $v_5$
								&	\centering $v_6$	&	\centering $v_7$	& \centering $v_8$	&	\centering $v_9$	&	\centering $v_{10}$	&
								\\
	\hline
	Random key		&	\centering 0			& \centering -			& \centering 0.8		&	\centering 0.7		&	\centering 0.9
								&	\centering -			&	\centering -			&	\centering 0.2		&	\centering -			&	\centering -	 &
								\\
	\end{tabular}
\\
\\
\\
	Offspring 1
&
	\setlength{\tabcolsep}{2.5pt}
	\begin{tabular}{p{0.8in}|p{0.2in}|p{0.2in}|p{0.2in}|p{0.2in}|p{0.2in}|p{0.2in}||p{0.2in}|p{0.2in}|p{0.2in}|p{0.2in}p{0in}}
	Vertex				&	\centering $v_1$	& \centering $v_2$	&	\centering $v_3$	&	\centering $v_4$	&	\centering $v_5$
								&	\centering $v_6$	&	\centering $v_7$	& \centering $v_8$	&	\centering $v_9$	&	\centering $v_{10}$	&
								\\
	\hline
	Random key		&	\centering 0			& \centering 0.7		& \centering - 			&	\centering 0.3		&	\centering -
								&	\centering 0.8		&	\centering -			&	\centering 0.2		&	\centering -			&	\centering -	 &
								\\
	\end{tabular}
\\
\\
	Offspring 2
&
	\setlength{\tabcolsep}{2.5pt}
	\begin{tabular}{p{0.8in}|p{0.2in}|p{0.2in}|p{0.2in}|p{0.2in}|p{0.2in}|p{0.2in}||p{0.2in}|p{0.2in}|p{0.2in}|p{0.2in}p{0in}}
	Vertex				&	\centering $v_1$	& \centering $v_2$	&	\centering $v_3$	&	\centering $v_4$	&	\centering $v_5$
								&	\centering $v_6$	&	\centering $v_7$	& \centering $v_8$	&	\centering $v_9$	&	\centering $v_{10}$	&
								\\
	\hline
	Random key		&	\centering 0			& \centering -			& \centering 0.8		&	\centering 0.7		&	\centering 0.9
								&	\centering -			&	\centering 0.2		&	\centering -			&	\centering 0.5		&	\centering -	 &
								\\
	\end{tabular}
\\
\end{tabular}
\caption{Crossover operator.}
\label{fig:crossover}
\end{figure}

\noindent
Another crossover operator preserving a bigger part of the initial ring structures, has been experimented, 
but it was tending to reduce the number of nodes belonging to the cycle and then was causing a premature convergence with solutions having a small number of visited nodes.

\subsection{Experiments}

The metaheuristics described in the previous section have all been implemented using the ParadisEO-MOEO library%
\footnote{ParadisEO is available at \url{http://paradiseo.gforge.inria.fr}.}%
~\cite{LB+:07}.
ParadisEO-MOEO is a C++ white-box object-oriented framework dedicated to the reusable design of metaheuristics for multi-objective optimization.
All the algorithms share the same base components for a fair comparison between them.
Computational runs were performed on an Intel Core $2$ Duo $6600$ ($2 \times 2.40$ GHz) machine, with $2$ GB RAM.

\subsubsection{Experimental Protocol}
\label{sec:protocol}

\paragraph{Benchmark Test Instances.}
Search methods have all been tested on different benchmark instances taken from the TSPLIB%
\footnote{\url{http://www.iwr.uni-heidelberg.de/groups/comopt/software/TSPLIB95/}.}%
~\cite{Rei:91}.
These instances involve between~$50$ and $300$ nodes.
The number at the end of an instance's name represents the number of nodes for the instance under consideration.
Let $l_{ij}$ denote the distance between two nodes $v_i$ and $v_j$ of a TSPLIB file.
Then, the ring cost $c_{ij}$ and the assignment cost $d_{ij}$ have both been set to $l_{ij}$ for every pair of nodes $v_i$ and $v_j$.
Moreover, note that, for each resolution method proposed to tackle the bi-objective RSP, the search process stops after a certain amount of run time.
As shown in Table~\ref{tab:stop}, this stopping criteria has been arbitrary set according to the size of the instance under consideration.
\begin{table}[htb]
\centering
\caption{Stopping criteria: running time.}
\hspace{500pt} \small \setlength{\tabcolsep}{3pt}
\begin{tabular}{ccc}
\begin{tabular}{cl|cc}
&  \centering Instance  &  Running  &  \\
&                       &  time     &  \\
\hline
&  {\itshape eil51}    &   20'' &  \\
&  {\itshape st70}     &    1'  &  \\
&  {\itshape kroA100}  &    2'  &  \\
&  {\itshape bier127}  &    5'  &  \\
\end{tabular}
&&
\begin{tabular}{cl|cc}
&  \centering Instance  &  Running  &  \\
&                       &  time     &  \\
\hline
&  {\itshape kroA150}  &   10'  &  \\
&  {\itshape kroA200}  &   20'  &  \\
&  {\itshape pr264}    &   30'  &  \\
&  {\itshape pr299}    &   50'  &  \\
\end{tabular}
\end{tabular}
\label{tab:stop}
\end{table}

\begin{figure}[t]
\centering
\includegraphics[width=4.5in]{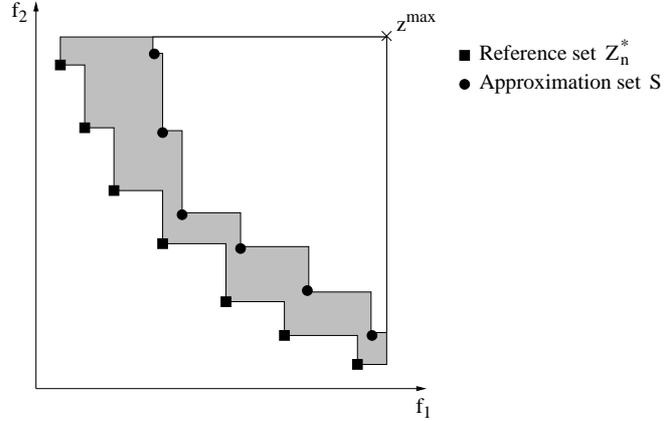}
\caption{Illustration of the hypervolume difference (I$^-_H$) between a reference set~$Z_N^\star$ and a non-dominated front approximation~$S$ (shaded area).}
\label{fig:hypervolume}
\end{figure}

\paragraph{Performance Assessment.}
For each one of the search method, a set of $20$ runs per instance, with different initial populations, has been performed.
In order to evaluate the quality of the non-dominated front approximations obtained for a specific test instance, 
we follow the protocol given by Knowles et al.~\cite{KTZ:06}.
For a specific instance, we first compute a reference set~$Z_N^\star$ of non-dominated points extracted from the union of all these fronts.
Second, we define $z^{max} = ( z^{max}_1 , z^{max}_2 )$, 
where $z^{max}_1$ (respectively $z^{max}_2$) denotes the upper bound of the first (respectively second) objective in the whole non-dominated front approximations.
Then, to measure the quality of an output set~$A$ in comparison to $Z_N^\star$, 
we compute the difference between these two sets by using the unary hypervolume metric~\cite{ZT:99}, 
$( 1.05 \times z^{max}_1 , 1.05 \times z^{max}_2 )$ being the reference point.
As illustrated in Figure~\ref{fig:hypervolume}, 
the hypervolume difference indicator (I$^-_H$) computes the portion of the objective space that is weakly dominated by $Z_N^\star$ and not by~$A$.
The more this measure is close to~$0$, the better is the approximation~$A$.
Furthermore, we also consider the additive $\epsilon$-indicator proposed in~\cite{ZT+:03}.
Contrary to the one proposed in Equation~\ref{eq:eps}, this indicator is used to compare non-dominated set approximations, and not solutions.
The unary additive $\epsilon$-indicator (I$^1_{\epsilon+}$) gives the minimum factor by which an approximation~$A$ 
has to be translated in the criterion space to weakly dominate the reference set~$Z_N^\star$.
I$^1_{\epsilon+}$ can be defined as follows:
\begin{equation}
I^1_{\epsilon+}(A) = I_{\epsilon+}(A,Z_N^\star) \ ,
\label{eq:eps2}
\end{equation}
where
\begin{equation}
I_{\epsilon+}(A,B) = \min_{\epsilon} \{ \forall z \in B, \exists z' \in A: z'_i - \epsilon \leq z_i, \forall 1 \leq i \leq n \} \ .
\label{eq:eps3}
\end{equation}
As a consequence, for each test instance, we obtain $20$~hypervolume differences and $20$~I$_{\epsilon+}$ measures, corresponding to the $20$ runs, per algorithm.
As suggested by Knowles et al.~\cite{KTZ:06}, once all these values are computed, 
we perform a statistical analysis on pairs of optimization methods for a comparison on a specific test instance.
To this end, we use the Mann-Whitney statistical test as described in~\cite{KTZ:06}, with a p-value lower than $5\%$.
Hence, for a specific test instance, and according to the p-value and to the metric under consideration, 
this statistical test reveals if the sample of approximation sets obtained by a search method~$M_1$ 
is significantly better than the sample of approximation sets obtained by a search method~$M_2$, or if there is no significant difference between both.
Note that all the performance assessment procedures have been achieved using the performance assessment tool suite provided in PISA%
\footnote{The package is available at \url{http://www.tik.ee.ethz.ch/pisa/assessment.html}.}%
~\cite{BL+:03}.

\subsubsection{Local Search Parameter Analysis}
\label{sec:ibmols_study}

\paragraph{Parameter Setting.}
For the IBMOLS algorithm, we set the population size according to the number of vertexes involved in the instance (see Table~\ref{tab:param}).
For each instance, a small~(S), a medium~(M) and a large~(L) population size have been experimented.
Likewise, the noise rate for the population re-initialization in the iterated version of IBMOLS (denoted by I-IBMOLS) is set to a fixed percentage of the instance's size.
We investigate three different values for this noise rate: $5\%$, $10\%$ and $20\%$.
Then, $0.05 \times n$, $0.1 \times n$ and $0.2 \times n$ random mutations are applied respectively for a problem with $n$~nodes.
Finally, the scaling factor $\kappa$ of the Equation~\ref{eq:fitness} is set to $0.05$.
\begin{table}[htb]
\centering
\caption{Population size parameter for the I-IBMOLS algorithm.}
\hspace{500pt} \small \setlength{\tabcolsep}{3pt}
\begin{tabular}{cr|p{0.4in}|p{0.4in}|p{0.4in}c}
&  \centering Instance  &  \multicolumn{3}{c}{Population size}                &  \\
&                       &  \centering S  &  \centering M   &  \centering L    &  \\
\hline
&  {\itshape eil51}    &  \centering 10  &  \centering 20  &  \centering 30   &  \\
&  {\itshape st70}     &  \centering 10  &  \centering 20  &  \centering 30   &  \\
&  {\itshape kroA100}  &  \centering 20  &  \centering 30  &  \centering 50   &  \\
&  {\itshape bier127}  &  \centering 20  &  \centering 30  &  \centering 50   &  \\
&  {\itshape kroA150}  &  \centering 20  &  \centering 30  &  \centering 50   &  \\
&  {\itshape kroA200}  &  \centering 20  &  \centering 30  &  \centering 50   &  \\
&  {\itshape pr264}    &  \centering 30  &  \centering 50  &  \centering 70   &  \\
&  {\itshape pr299}    &  \centering 30  &  \centering 50  &  \centering 70   &  \\
\end{tabular}
\label{tab:param}
\end{table}

\paragraph{Computational Results and Discussion.}
The results we obtained for the I$^-_H$ and the I$^1_{\epsilon+}$ metrics are respectively summarized in Table~\ref{tab:ibmols_hyp} and Table~\ref{tab:ibmols_eps}.
These tables give, for every instance, the number of times that a specific para\-metrization is significantly worse than another one for the I-IBMOLS algorithm.
So, a lower score is better, and {\itshape zero} means that the parametrization under consideration is never significantly outperformed.
According to the I$^-_H$ metric, a medium population size with a random noise of $5\%$ or $10\%$ seems to perform well on every test instance,
except for the {\itshape pr264} instance where a medium population size and a random noise of $10\%$ is significantly worse than another parametrization.
Moreover, as shown in Table~\ref{tab:ibmols_eps}, a medium population size with a random noise parametrized to $10\%$ 
generally seems to perform slightly better than with a random noise of $5\%$ according to the I$^1_{\epsilon+}$ metric.
The only instances for which such a parametrization is significantly worse than at least another one are the {\itshape kroA200} and the {\itshape pr264} problems.
Then, for all the other instances, we choose to set the number of individuals in the population to the medium size and the noise rate to $10\%$ 
for comparing I-IBMOLS to the other metaheuristics.
For the {\itshape kroA200} instance, the same parametrization will be used as it seems to be the best compromise between the two metrics.
Finally, for the {\itshape pr264} instance, we will use a small population size and a noise rate of $20\%$.
The value of all the parameters used by I-IBMOLS for the remainder of the paper are summarized in Table~\ref{tab:param_compare}.

\begin{table}[htbp]
\centering
\caption{
Performance comparison of different population size and noise rate parameter settings for I-IBMOLS according to the I$^-_H$ metric 
by using a Mann-Whitney statistical test with a significance level of $5\%$.
Each column gives the number of algorithms with another parameter setting by which it is statistically outperformed for the instance under consideration.
}
\hspace{500pt} \small \setlength{\tabcolsep}{3pt}
\begin{tabular}{cr|p{0.25in}|p{0.25in}|p{0.25in}|p{0.25in}|p{0.25in}|p{0.25in}|p{0.25in}|p{0.25in}|p{0.25in}c}
&  noise rate &  \multicolumn{3}{|c|}{5\%}                       &  \multicolumn{3}{|c|}{10\%}                      &  \multicolumn{3}{|c}{20\%}                       & \\
&  pop. size  &  \centering S  &  \centering M  &  \centering L  &  \centering S  &  \centering M  &  \centering L  &  \centering S  &  \centering M  &  \centering L  & \\
\hline
& {\itshape eil51}    &  \centering  3              &  \centering  {\bfseries 0}  &  \centering  3
											&  \centering  3              &  \centering  {\bfseries 0}  &  \centering  3
											&  \centering  8              &  \centering  1              &  \centering  3              &  \\
& {\itshape st70}     &  \centering  6              &  \centering  {\bfseries 0}  &  \centering  {\bfseries 0}
											&  \centering  7              &  \centering  {\bfseries 0}  &  \centering  {\bfseries 0}
											&  \centering  5              &  \centering  {\bfseries 0}  &  \centering  1              &  \\
& {\itshape kroA100}  &  \centering  1              &  \centering  {\bfseries 0}  &  \centering  6
											&  \centering  1              &  \centering  {\bfseries 0}  &  \centering  6
											&  \centering  2              &  \centering  2              &  \centering  6              &  \\
& {\itshape bier127}  &  \centering  1              &  \centering  {\bfseries 0}  &  \centering  6
											&  \centering  1              &  \centering  {\bfseries 0}  &  \centering  6
											&  \centering  5              &  \centering  1              &  \centering  6              &  \\
& {\itshape kroA150}  &  \centering  1              &  \centering  {\bfseries 0}  &  \centering  6
											&  \centering  4              &  \centering  {\bfseries 0}  &  \centering  6
											&  \centering  5              &  \centering  1              &  \centering  6              &  \\
& {\itshape kroA200}  &  \centering  2              &  \centering  {\bfseries 0}  &  \centering  1
											&  \centering  2              &  \centering  {\bfseries 0}  &  \centering  2
											&  \centering  8              &  \centering  2              &  \centering  2              &  \\
& {\itshape pr264}    &  \centering  5              &  \centering  {\bfseries 0}  &  \centering  2
											&  \centering  2              &  \centering  1              &  \centering  4
											&  \centering  {\bfseries 0}  &  \centering  1              &  \centering  5              &  \\
& {\itshape pr299}    &  \centering  6              &  \centering  {\bfseries 0}  &  \centering  {\bfseries 0}
											&  \centering  8              &  \centering  {\bfseries 0}  &  \centering  3
											&  \centering  3              &  \centering  {\bfseries 0}  &  \centering  4              &  \\
\end{tabular}
\label{tab:ibmols_hyp}
\end{table}

\begin{table}[htbp]
\centering
\caption{
Performance comparison of different population size and noise rate parameter settings for I-IBMOLS according to the I$^1_{\epsilon+}$ metric 
by using a Mann-Whitney statistical test with a significance level of $5\%$.
Each column gives the number of algorithms with another parameter setting by which it is statistically outperformed for the instance under consideration.
}
\hspace{500pt} \small \setlength{\tabcolsep}{3pt}
\begin{tabular}{cr|p{0.25in}|p{0.25in}|p{0.25in}|p{0.25in}|p{0.25in}|p{0.25in}|p{0.25in}|p{0.25in}|p{0.25in}c}
&  noise rate &  \multicolumn{3}{|c|}{5\%}                       &  \multicolumn{3}{|c|}{10\%}                      &  \multicolumn{3}{|c}{20\%}                       & \\
&  pop. size  &  \centering S  &  \centering M  &  \centering L  &  \centering S  &  \centering M  &  \centering L  &  \centering S  &  \centering M  &  \centering L  & \\
\hline
& {\itshape eil51}    &  \centering  7              &  \centering  {\bfseries 0}  &  \centering  3
											&  \centering  6              &  \centering  {\bfseries 0}  &  \centering  3
											&  \centering  5              &  \centering  {\bfseries 0}  &  \centering  3              &  \\
& {\itshape st70}     &  \centering  7              &  \centering  {\bfseries 0}  &  \centering  {\bfseries 0}
											&  \centering  7              &  \centering  {\bfseries 0}  &  \centering  {\bfseries 0}
											&  \centering  6              &  \centering  {\bfseries 0}  &  \centering  4              &  \\
& {\itshape kroA100}  &  \centering  {\bfseries 0}  &  \centering  {\bfseries 0}  &  \centering  6
											&  \centering  {\bfseries 0}  &  \centering  {\bfseries 0}  &  \centering  6
											&  \centering  1              &  \centering  2              &  \centering  6              &  \\
& {\itshape bier127}  &  \centering  {\bfseries 0}  &  \centering  {\bfseries 0}  &  \centering  6
											&  \centering  1              &  \centering  {\bfseries 0}  &  \centering  6
											&  \centering  1              &  \centering  1              &  \centering  6              &  \\
& {\itshape kroA150}  &  \centering  2              &  \centering  1              &  \centering  6
											&  \centering  3              &  \centering  {\bfseries 0}  &  \centering  6
											&  \centering  2              &  \centering  {\bfseries 0}  &  \centering  6              &  \\
& {\itshape kroA200}  &  \centering  8              &  \centering  1              &  \centering  1
											&  \centering  7              &  \centering  1              &  \centering  1
											&  \centering  1              &  \centering  {\bfseries 0}  &  \centering  1              &  \\
& {\itshape pr264}    &  \centering  2              &  \centering  2              &  \centering  6
											&  \centering  1              &  \centering  2              &  \centering  6
											&  \centering  {\bfseries 0}  &  \centering  2              &  \centering  6              &  \\
& {\itshape pr299}    &  \centering  8              &  \centering  3              &  \centering  3
											&  \centering  6              &  \centering  {\bfseries 0}  &  \centering  3
											&  \centering  {\bfseries 0}  &  \centering  {\bfseries 0}  &  \centering  6              &  \\
\end{tabular}
\label{tab:ibmols_eps}
\end{table}

\subsubsection{Evolutionary Algorithms Parameter Analysis}

\paragraph{Parameter Setting.}
For all the EAs, population sizes of $50$, $100$ and $200$ individuals have been experimented.
Note that the instance-dependent population size setting used for I-IBMOLS has also been investigated, but the results were overall less efficient.
The scaling factor~$\kappa$ used by IBEA is set to $0.05$.
At last, the crossover probability is set to $0.25$,
and the mutation probability to $1.00$, with a probability of $0.25$, $0.25$ and $0.50$ for the \emph{remove}, the \emph{insert} and the \emph{2-opt} operator, 
respectively.

\paragraph{Computational Results and Discussion.}
As pointed out in Table~\ref{tab:ibea}, the results obtained by IBEA are quite heterogeneous from a test instance to another, and no general conclusion can be drawn.
Thus, we will use a population size of $50$ for the {\itshape eil51}, the {\itshape st70}, the {\itshape pr264} and the {\itshape pr299} instances,
a population size of $100$ for the {\itshape kroA100} and the {\itshape bier127} instances,
and a population size of $200$ for the {\itshape kroA150} and the {\itshape kroA200} instances.
For NSGA-II, the best performances are usually obtained with a large population size.
Then, a population of $100$ individuals for instances with less than $100$ nodes and of $200$ individuals for instances with more than $100$ nodes will be used for the following.
These parameters are never significantly outperformed according to the metrics we used.
Finally, as shown in Table~\ref{tab:seea}, the population size does not appear to have a big influence on the outcome of SEEA, 
what is not surprising given the design of the algorithm (see Section~\ref{sec:seea}).
As a consequence, 
we will use a population of $100$ individuals as this is the only parametrization that is strictly never significantly outperformed for every test instance.
All the parameters used by the EAs for the next section are recapitulated in Table~\ref{tab:param_compare}.

\begin{table}[htbp]
\centering
\caption{
Performance comparison of different population size parameter settings for IBEA according to the I$^-_H$ and the I$^1_{\epsilon+}$ metrics 
by using a Mann-Whitney statistical test with a significance level of $5\%$.
Each column gives the number of algorithms with another parameter setting by which it is statistically outperformed for the instance under consideration.
}
\hspace{500pt} \small \setlength{\tabcolsep}{3pt}
\begin{tabular}{cr||p{0.25in}|p{0.25in}|p{0.25in}||p{0.25in}|p{0.25in}|p{0.25in}c}
&                 &  \multicolumn{3}{|c||}{I$^-_H$}                        &  \multicolumn{3}{|c}{I$^1_{\epsilon+}$}                 &  \\
&  pop. size      &  \centering 50  &  \centering 100   &  \centering 200  &  \centering 50  &  \centering 100  &  \centering 200  &  \\
\hline
& {\itshape eil51}    &  \centering  {\bfseries 0}  &  \centering  {\bfseries 0}  &  \centering  2
											&  \centering  {\bfseries 0}  &  \centering  {\bfseries 0}  &  \centering  2              &  \\
& {\itshape st70}     &  \centering  {\bfseries 0}  &  \centering  {\bfseries 0}  &  \centering  2
											&  \centering  {\bfseries 0}  &  \centering  1              &  \centering  1              &  \\
& {\itshape kroA100}  &  \centering  {\bfseries 0}  &  \centering  {\bfseries 0}  &  \centering  1
											&  \centering  {\bfseries 0}  &  \centering  {\bfseries 0}  &  \centering  {\bfseries 0}  &  \\
& {\itshape bier127}  &  \centering  {\bfseries 0}  &  \centering  {\bfseries 0}  &  \centering  {\bfseries 0}
											&  \centering  {\bfseries 0}  &  \centering  {\bfseries 0}  &  \centering  {\bfseries 0}  &  \\
& {\itshape kroA150}  &  \centering  1              &  \centering  {\bfseries 0}  &  \centering  {\bfseries 0}
											&  \centering  1              &  \centering  {\bfseries 0}  &  \centering  {\bfseries 0}  &  \\
& {\itshape kroA200}  &  \centering  {\bfseries 0}  &  \centering  2              &  \centering  {\bfseries 0}
											&  \centering  1              &  \centering  2              &  \centering  {\bfseries 0}  &  \\
& {\itshape pr264}    &  \centering  {\bfseries 0}  &  \centering  1              &  \centering  1
											&  \centering  {\bfseries 0}  &  \centering  1              &  \centering  1              &  \\
& {\itshape pr299}    &  \centering  {\bfseries 0}  &  \centering  1              &  \centering  1
											&  \centering  {\bfseries 0}  &  \centering  1              &  \centering  1              &  \\
\end{tabular}
\label{tab:ibea}
\end{table}

\begin{table}[htbp]
\centering
\caption{
Performance comparison of different population size parameter settings for NSGA-II according to the I$^-_H$ and the I$^1_{\epsilon+}$ metrics 
by using a Mann-Whitney statistical test with a significance level of $5\%$.
Each column gives the number of algorithms with another parameter setting by which it is statistically outperformed for the instance under consideration.
}
\hspace{500pt} \small \setlength{\tabcolsep}{3pt}
\begin{tabular}{cr||p{0.25in}|p{0.25in}|p{0.25in}||p{0.25in}|p{0.25in}|p{0.25in}c}
&                 &  \multicolumn{3}{|c||}{I$^-_H$}                        &  \multicolumn{3}{|c}{I$^1_{\epsilon+}$}                 &  \\
&  pop. size      &  \centering 50  &  \centering 100   &  \centering 200  &  \centering 50  &  \centering 100  &  \centering 200  &  \\
\hline
& {\itshape eil51}    &  \centering  1  &  \centering  {\bfseries 0}  &  \centering  2
											&  \centering  1  &  \centering  {\bfseries 0}  &  \centering  2              &  \\
& {\itshape st70}     &  \centering  2  &  \centering  {\bfseries 0}  &  \centering  {\bfseries 0}
											&  \centering  2  &  \centering  {\bfseries 0}  &  \centering  {\bfseries 0}  &  \\
& {\itshape kroA100}  &  \centering  2  &  \centering  1              &  \centering  {\bfseries 0}
											&  \centering  2  &  \centering  1              &  \centering  {\bfseries 0}  &  \\
& {\itshape bier127}  &  \centering  2  &  \centering  1              &  \centering  {\bfseries 0}
											&  \centering  2  &  \centering  1              &  \centering  {\bfseries 0}  &  \\
& {\itshape kroA150}  &  \centering  2  &  \centering  1              &  \centering  {\bfseries 0}
											&  \centering  2  &  \centering  1              &  \centering  {\bfseries 0}  &  \\
& {\itshape kroA200}  &  \centering  2  &  \centering  1              &  \centering  {\bfseries 0}
											&  \centering  2  &  \centering  1              &  \centering  {\bfseries 0}  &  \\
& {\itshape pr264}    &  \centering  2  &  \centering  1              &  \centering  {\bfseries 0}
											&  \centering  2  &  \centering  1              &  \centering  {\bfseries 0}  &  \\
& {\itshape pr299}    &  \centering  2  &  \centering  1              &  \centering  {\bfseries 0}
											&  \centering  2  &  \centering  1              &  \centering  {\bfseries 0}  &  \\
\end{tabular}
\label{tab:nsga2}
\end{table}

\begin{table}[htbp]
\centering
\caption{
Performance comparison of different population size parameter settings for SEEA according to the I$^-_H$ and the I$^1_{\epsilon+}$ metrics 
by using a Mann-Whitney statistical test with a significance level of $5\%$.
Each column gives the number of algorithms with another parameter setting by which it is statistically outperformed for the instance under consideration.
}
\hspace{500pt} \small \setlength{\tabcolsep}{3pt}
\begin{tabular}{cr||p{0.25in}|p{0.25in}|p{0.25in}||p{0.25in}|p{0.25in}|p{0.25in}c}
&                 &  \multicolumn{3}{|c||}{I$^-_H$}                        &  \multicolumn{3}{|c}{I$^1_{\epsilon+}$}                 &  \\
&  pop. size      &  \centering 50  &  \centering 100   &  \centering 200  &  \centering 50  &  \centering 100  &  \centering 200  &  \\
\hline
& {\itshape eil51}    &  \centering  {\bfseries 0}  &  \centering  {\bfseries 0}  &  \centering  {\bfseries 0}
											&  \centering  {\bfseries 0}  &  \centering  {\bfseries 0}  &  \centering  {\bfseries 0}  &  \\
& {\itshape st70}     &  \centering  {\bfseries 0}  &  \centering  {\bfseries 0}  &  \centering  {\bfseries 0}
											&  \centering  {\bfseries 0}  &  \centering  {\bfseries 0}  &  \centering  {\bfseries 0}  &  \\
& {\itshape kroA100}  &  \centering  {\bfseries 0}  &  \centering  {\bfseries 0}  &  \centering  {\bfseries 0}
											&  \centering  {\bfseries 0}  &  \centering  {\bfseries 0}  &  \centering  {\bfseries 0}  &  \\
& {\itshape bier127}  &  \centering  2              &  \centering  {\bfseries 0}  &  \centering  {\bfseries 0}
											&  \centering  2              &  \centering  {\bfseries 0}  &  \centering  {\bfseries 0}  &  \\
& {\itshape kroA150}  &  \centering  {\bfseries 0}  &  \centering  {\bfseries 0}  &  \centering  {\bfseries 0}
											&  \centering  {\bfseries 0}  &  \centering  {\bfseries 0}  &  \centering  {\bfseries 0}  &  \\
& {\itshape kroA200}  &  \centering  {\bfseries 0}  &  \centering  {\bfseries 0}  &  \centering  {\bfseries 0}
											&  \centering  {\bfseries 0}  &  \centering  {\bfseries 0}  &  \centering  {\bfseries 0}  &  \\
& {\itshape pr264}    &  \centering  {\bfseries 0}  &  \centering  {\bfseries 0}  &  \centering  {\bfseries 0}
											&  \centering  {\bfseries 0}  &  \centering  {\bfseries 0}  &  \centering  1              &  \\
& {\itshape pr299}    &  \centering  {\bfseries 0}  &  \centering  {\bfseries 0}  &  \centering  {\bfseries 0}
											&  \centering  {\bfseries 0}  &  \centering  {\bfseries 0}  &  \centering  {\bfseries 0}  &  \\
											\end{tabular}
\label{tab:seea}
\end{table}

\subsubsection{Comparison between Search Methods}
\begin{figure}[p]
\small 
\begin{minipage}[c]{0.1\linewidth}
Instance
\end{minipage}
\begin{minipage}[c]{0.44\linewidth}
\centering I$^-_H$
\end{minipage}
\hfill
\begin{minipage}[c]{0.44\linewidth}
\centering I$^1_{\epsilon+}$
\end{minipage}
\begin{minipage}[c]{0.1\linewidth}
\centering
{\itshape eil51}
\end{minipage}
\begin{minipage}[c]{0.44\linewidth}
\centering
\includegraphics[angle=270,width=2in]{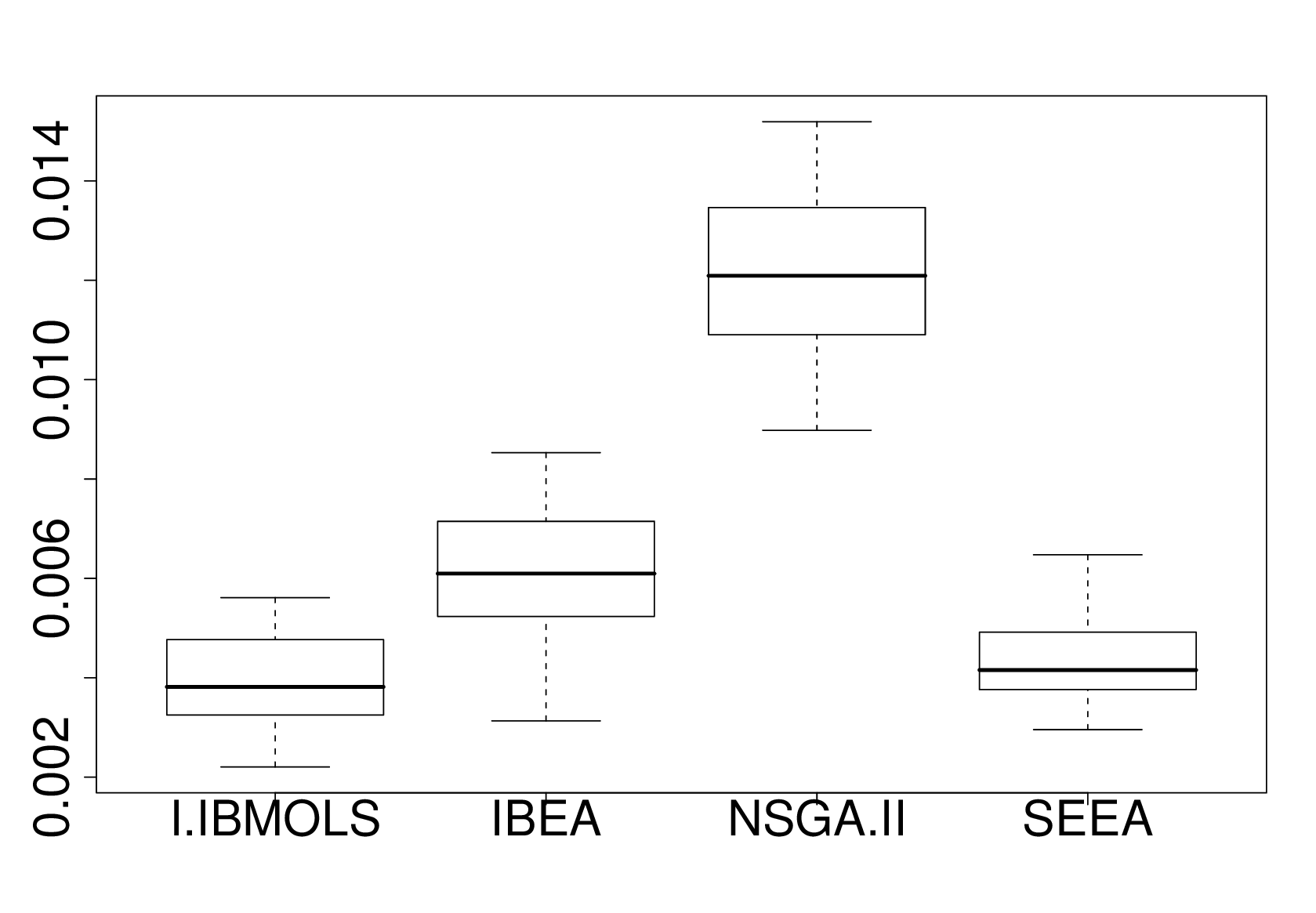}
\end{minipage}
\hfill
\begin{minipage}[c]{0.44\linewidth}
\centering
\includegraphics[angle=270,width=2in]{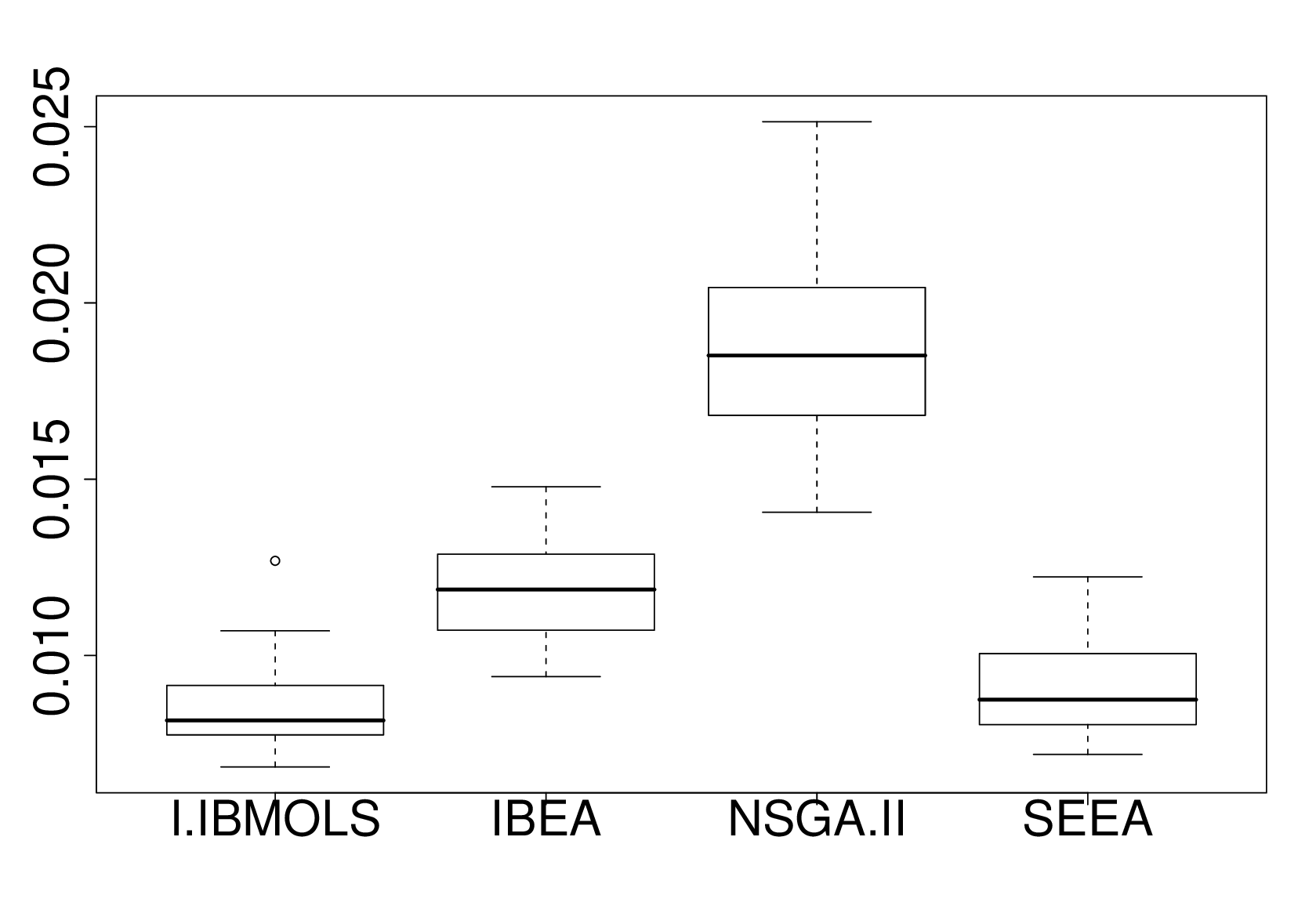}
\end{minipage}
\begin{minipage}[c]{0.1\linewidth}
\centering
{\itshape st70}
\end{minipage}
\begin{minipage}[c]{0.44\linewidth}
\centering
\includegraphics[angle=270,width=2in]{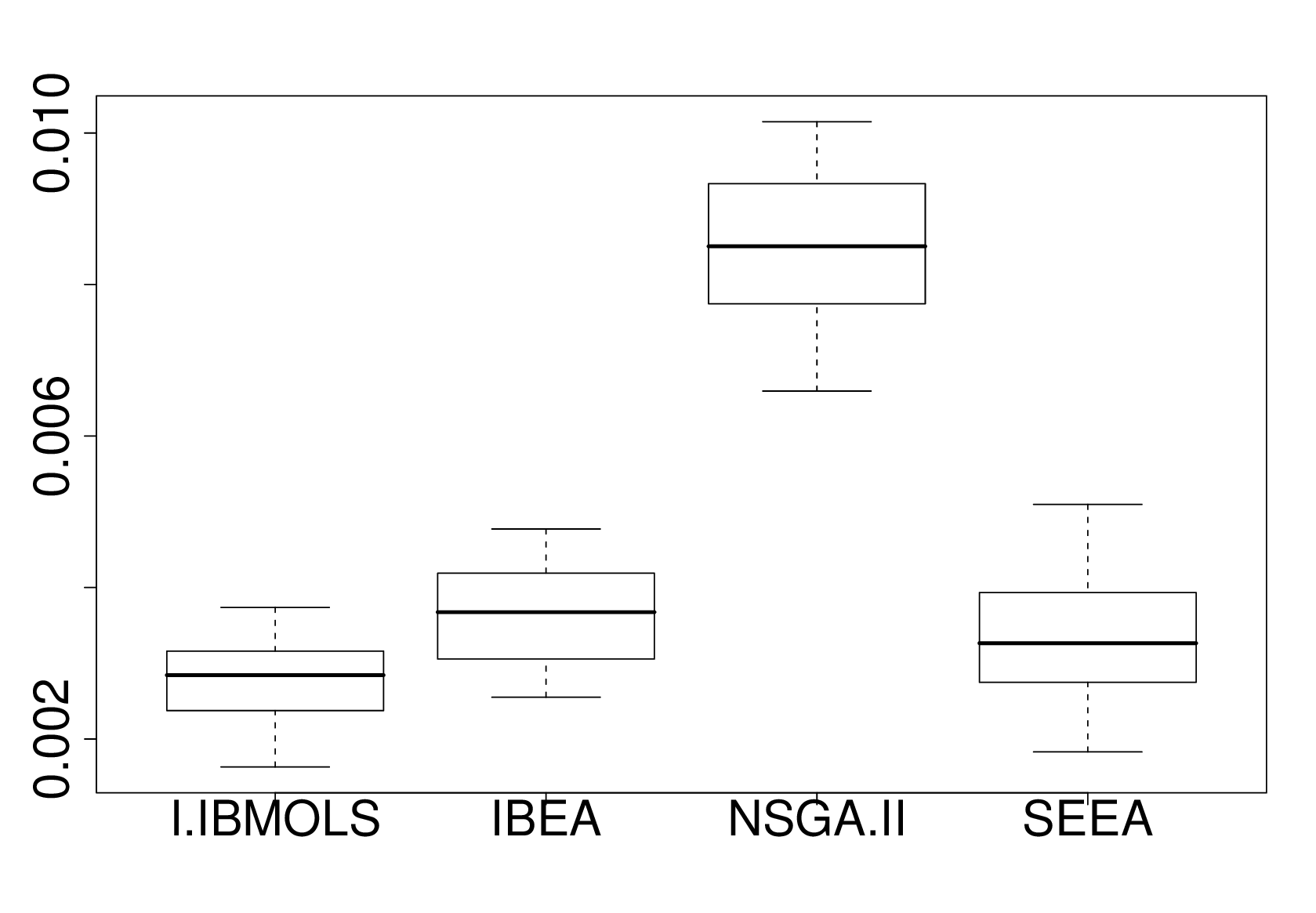}
\end{minipage}
\hfill
\begin{minipage}[c]{0.44\linewidth}
\centering
\includegraphics[angle=270,width=2in]{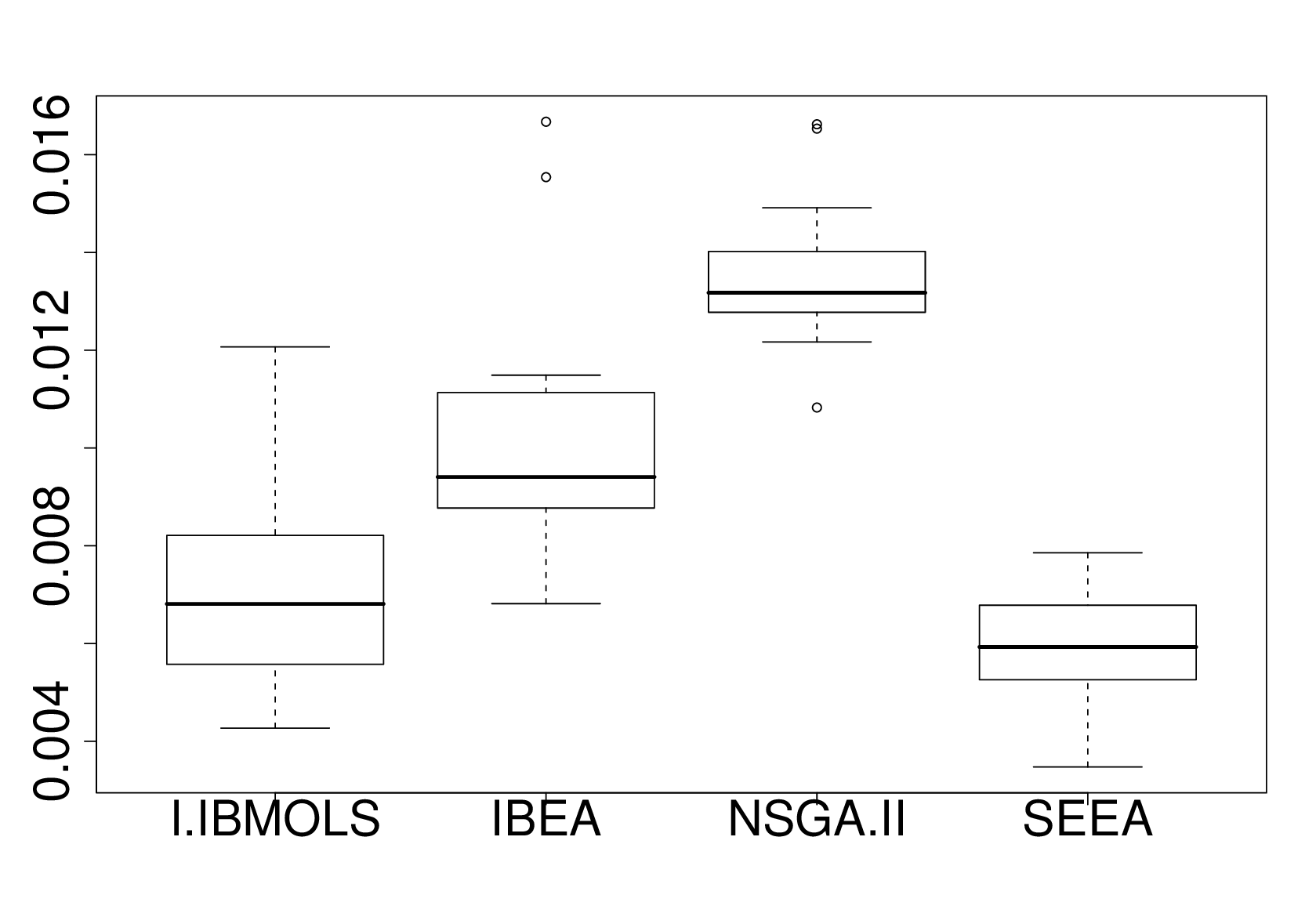}
\end{minipage}
\begin{minipage}[c]{0.1\linewidth}
\centering
{\itshape kroA100}
\end{minipage}
\begin{minipage}[c]{0.44\linewidth}
\centering
\includegraphics[angle=270,width=2in]{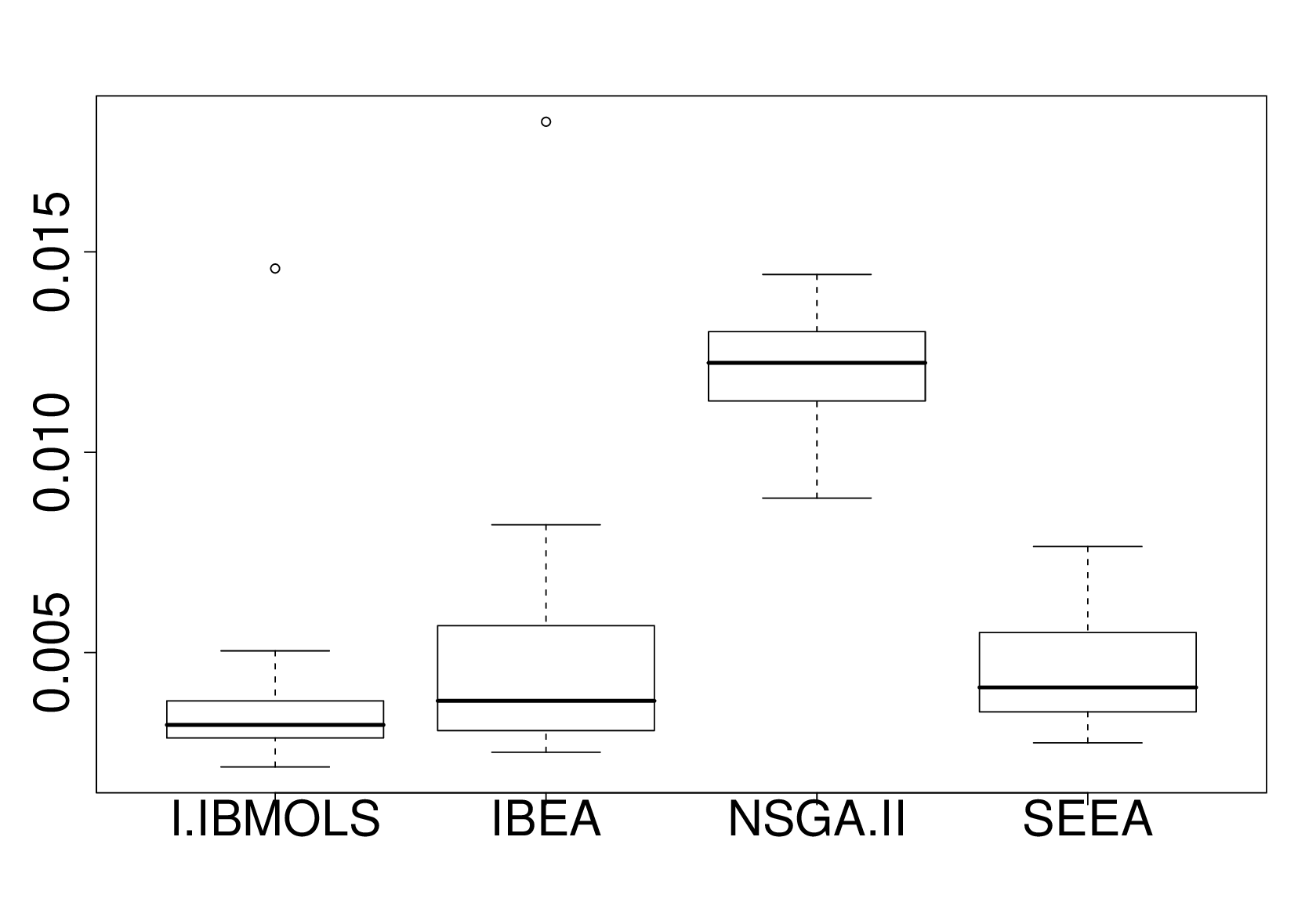}
\end{minipage}
\hfill
\begin{minipage}[c]{0.44\linewidth}
\centering
\includegraphics[angle=270,width=2in]{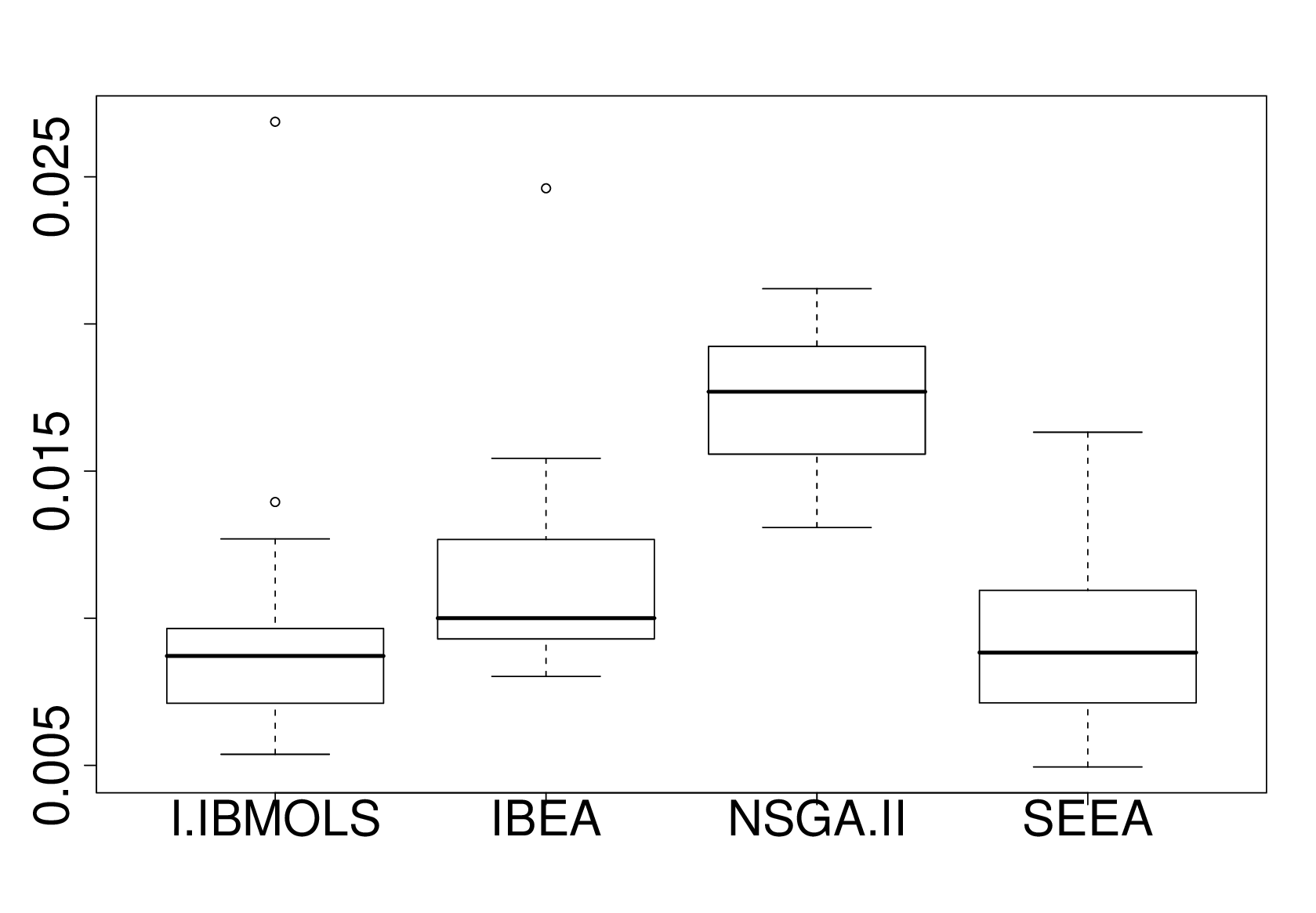}
\end{minipage}
\begin{minipage}[c]{0.1\linewidth}
\centering
{\itshape bier127}
\end{minipage}
\begin{minipage}[c]{0.44\linewidth}
\centering
\includegraphics[angle=270,width=2in]{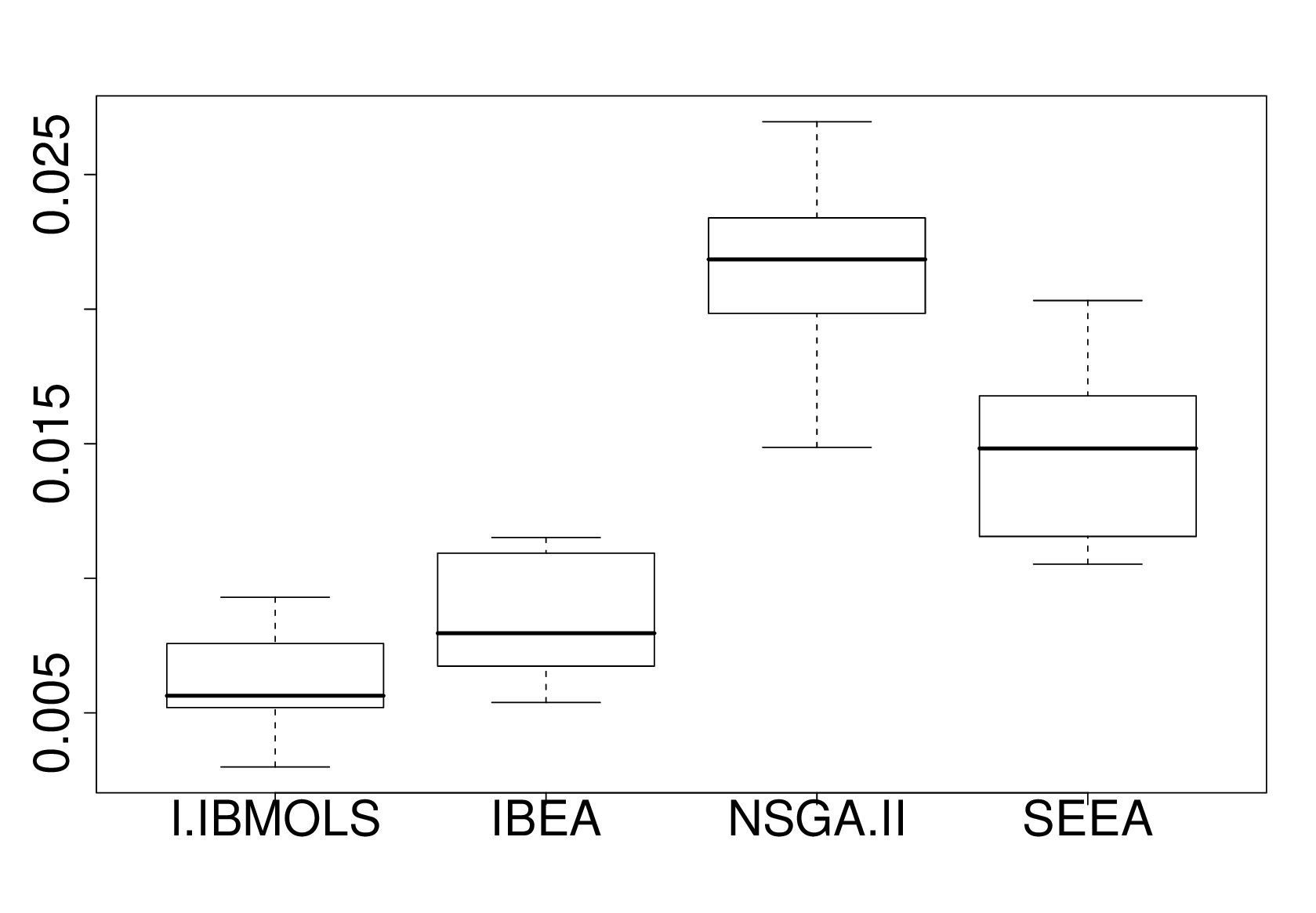}
\end{minipage}
\hfill
\begin{minipage}[c]{0.44\linewidth}
\centering
\includegraphics[angle=270,width=2in]{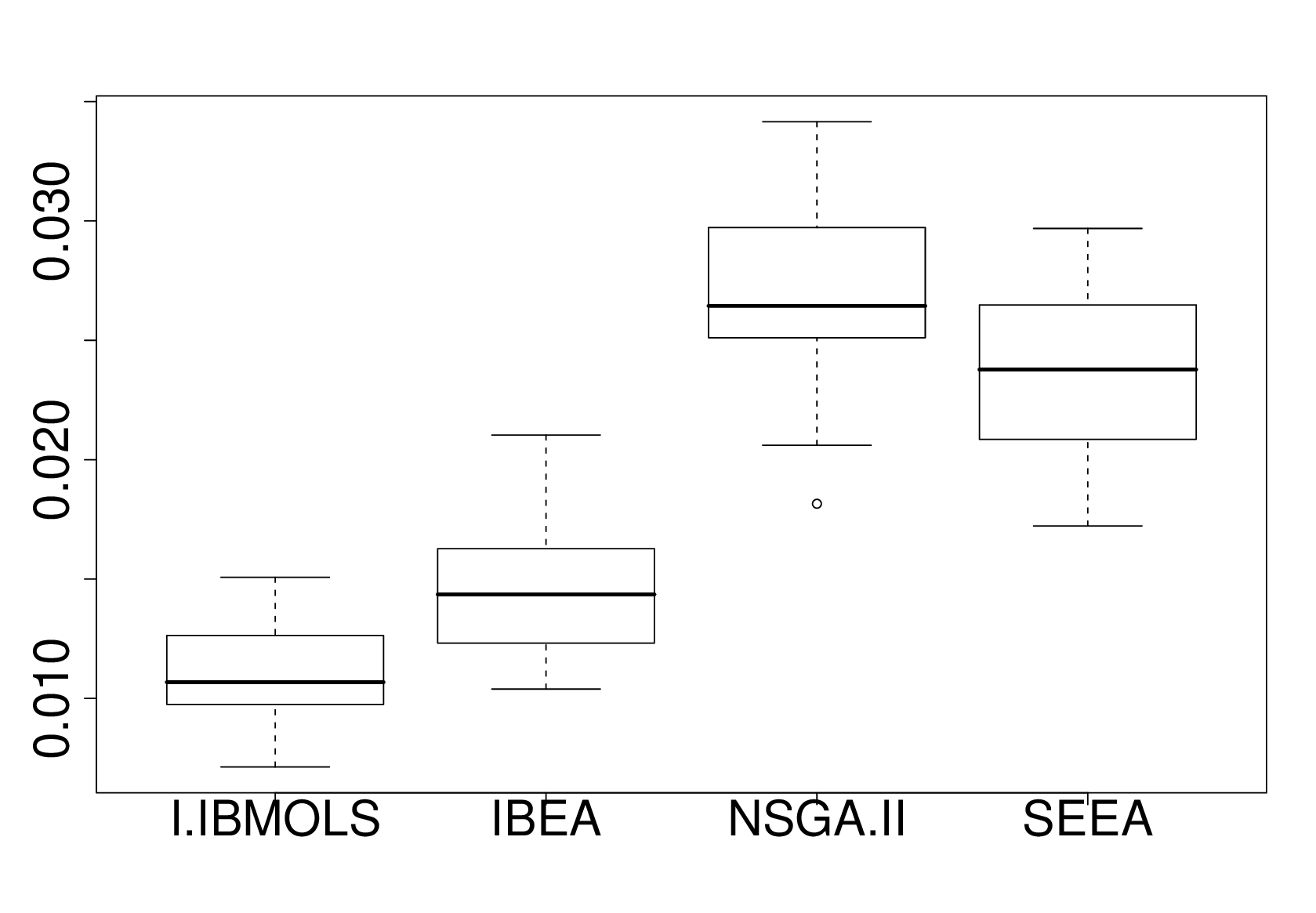}
\end{minipage}
\caption{Performance comparison for I-IBMOLS, IBEA, NSGA-II and SEEA according to the I$^-_H$ and the I$^1_{\epsilon+}$ metric (1).}
\label{fig:box_meta1}
\end{figure}

\begin{figure}[p]
\small 
\begin{minipage}[c]{0.1\linewidth}
Instance
\end{minipage}
\begin{minipage}[c]{0.44\linewidth}
\centering I$^-_H$
\end{minipage}
\hfill
\begin{minipage}[c]{0.44\linewidth}
\centering I$^1_{\epsilon+}$
\end{minipage}
\begin{minipage}[c]{0.1\linewidth}
\centering
{\itshape kroA150}
\end{minipage}
\begin{minipage}[c]{0.44\linewidth}
\centering
\includegraphics[angle=270,width=2in]{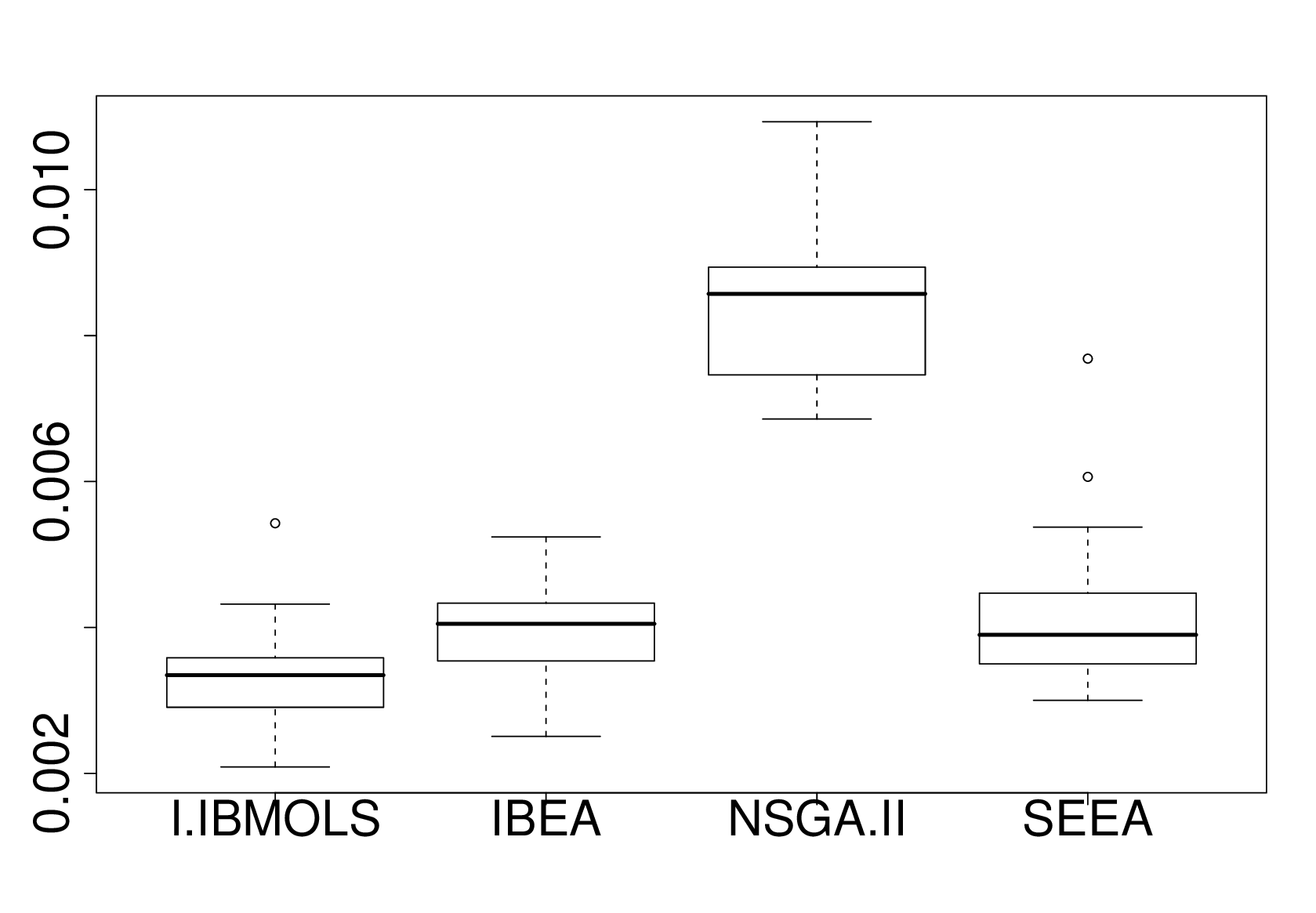}
\end{minipage}
\hfill
\begin{minipage}[c]{0.44\linewidth}
\centering
\includegraphics[angle=270,width=2in]{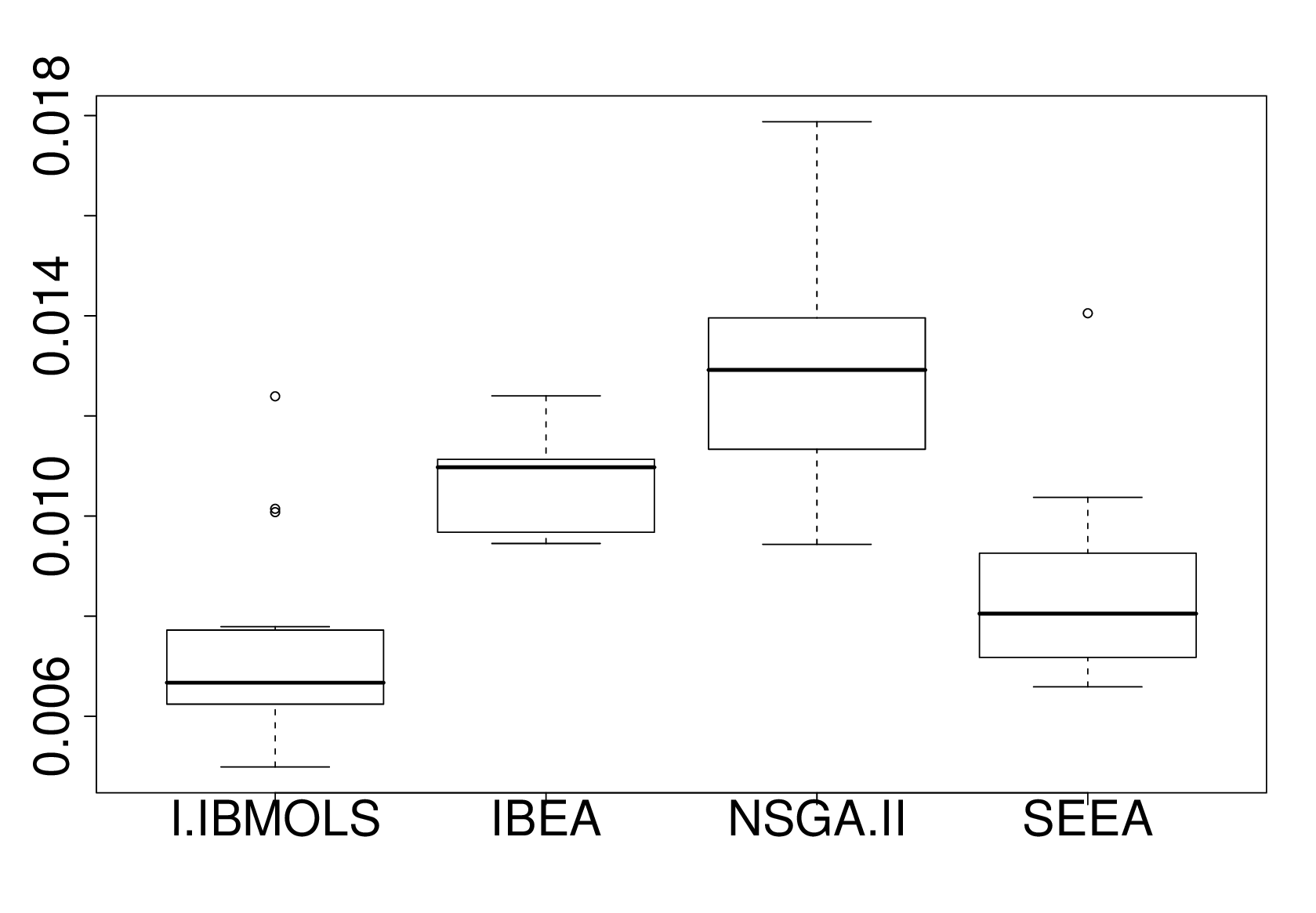}
\end{minipage}
\begin{minipage}[c]{0.1\linewidth}
\centering
{\itshape kroA200}
\end{minipage}
\begin{minipage}[c]{0.44\linewidth}
\centering
\includegraphics[angle=270,width=2in]{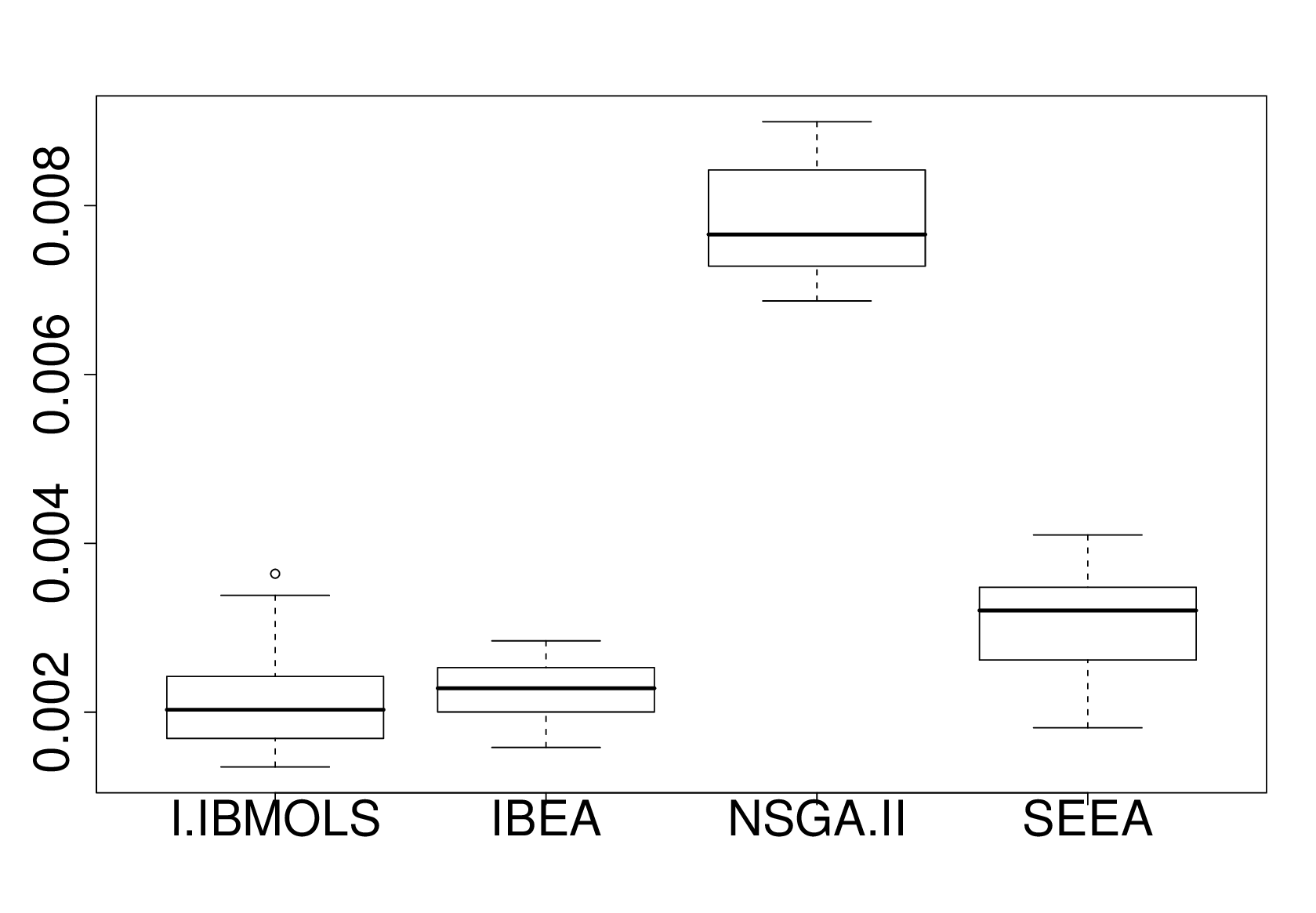}
\end{minipage}
\hfill
\begin{minipage}[c]{0.44\linewidth}
\centering
\includegraphics[angle=270,width=2in]{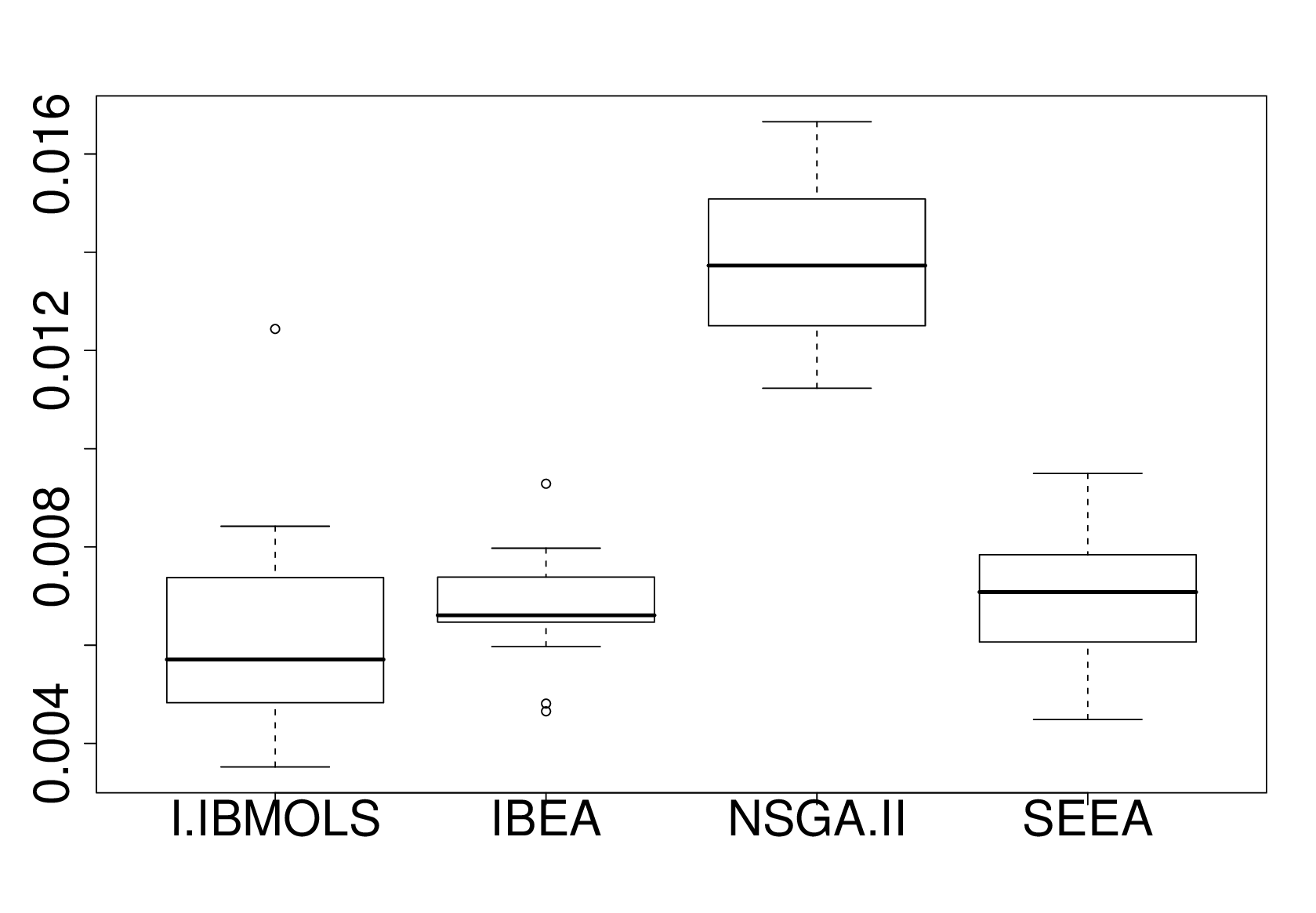}
\end{minipage}
\begin{minipage}[c]{0.1\linewidth}
\centering
{\itshape pr264}
\end{minipage}
\begin{minipage}[c]{0.44\linewidth}
\centering
\includegraphics[angle=270,width=2in]{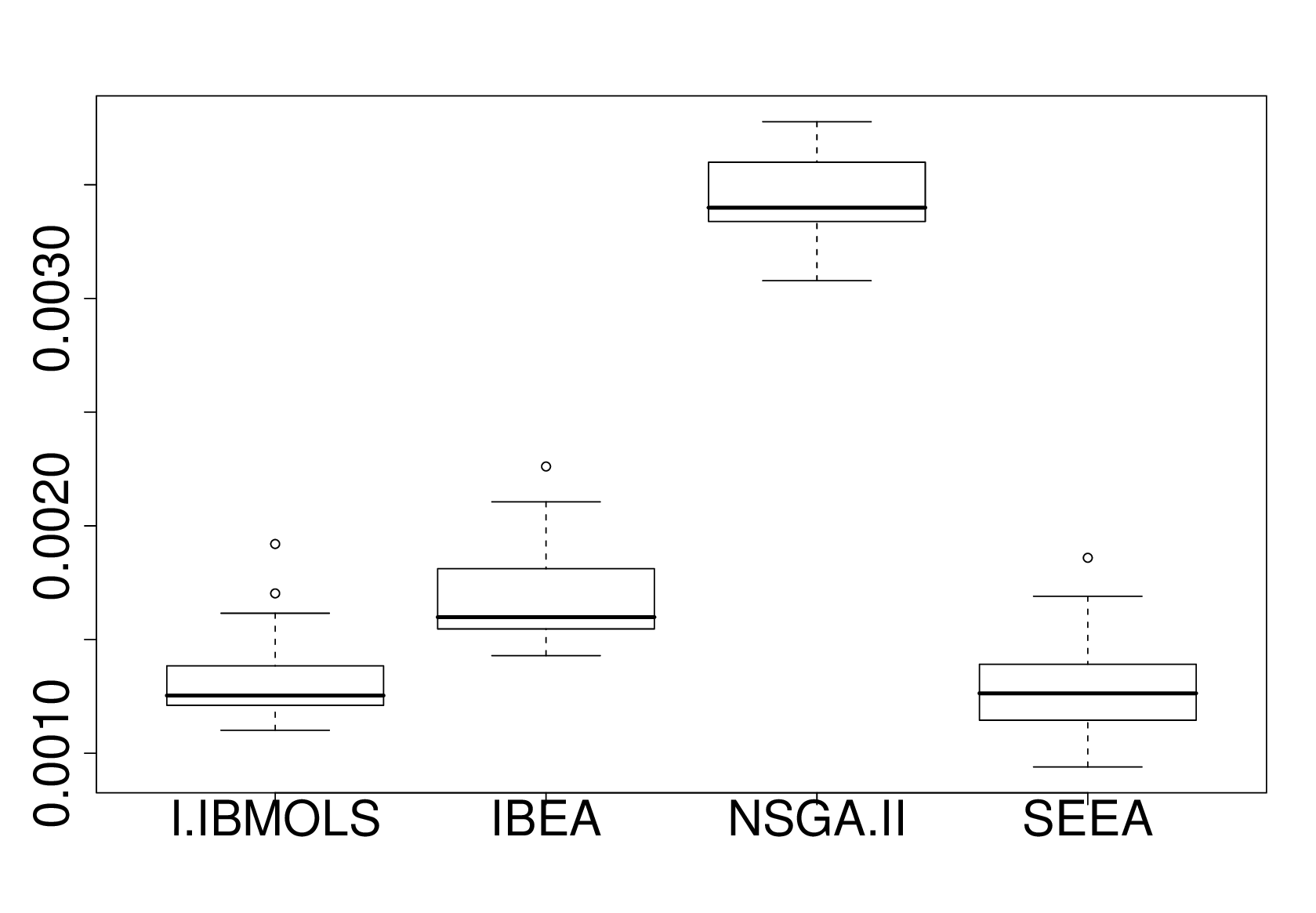}
\end{minipage}
\hfill
\begin{minipage}[c]{0.44\linewidth}
\centering
\includegraphics[angle=270,width=2in]{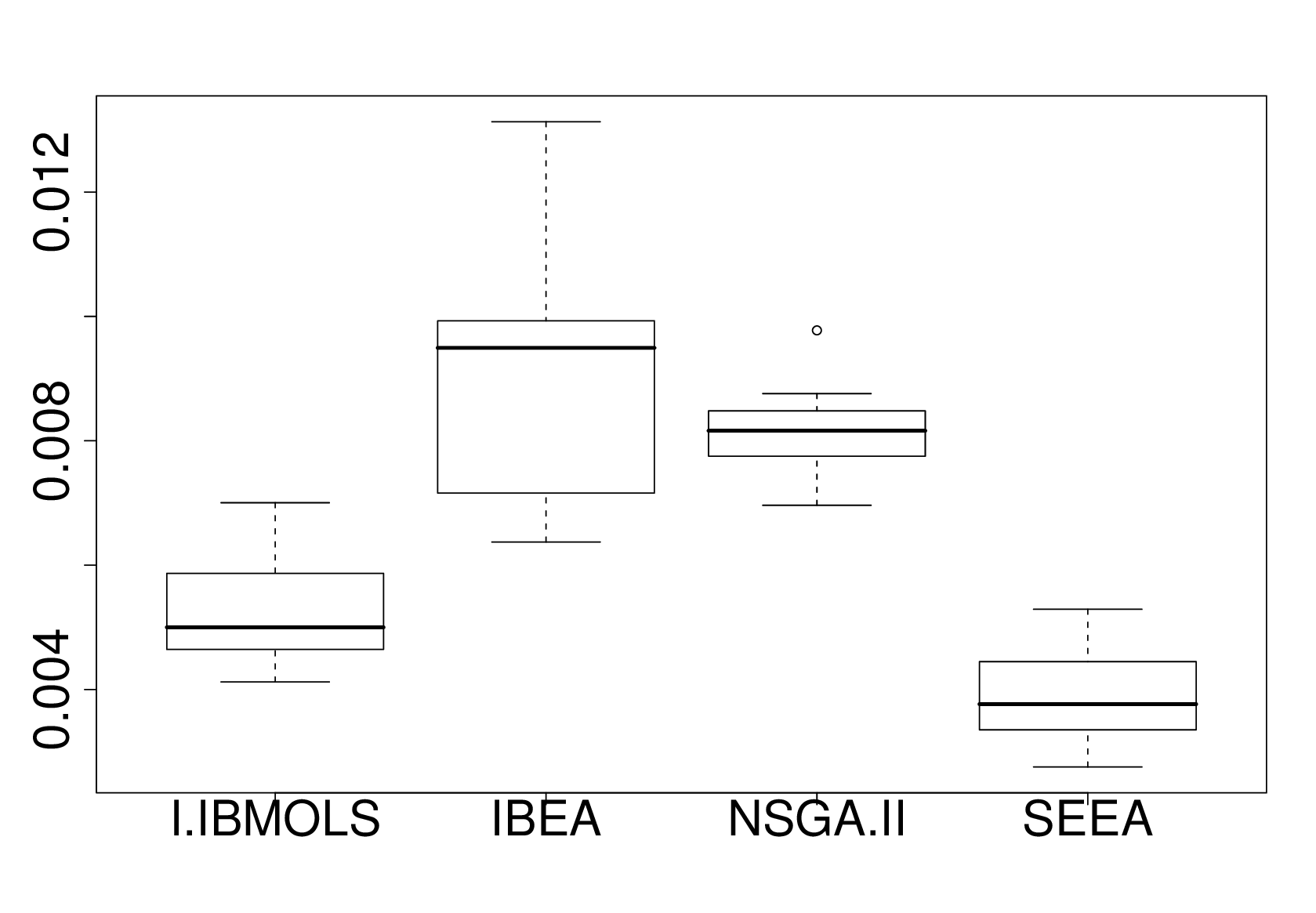}
\end{minipage}
\begin{minipage}[c]{0.1\linewidth}
\centering
{\itshape pr299}
\end{minipage}
\begin{minipage}[c]{0.44\linewidth}
\centering
\includegraphics[angle=270,width=2in]{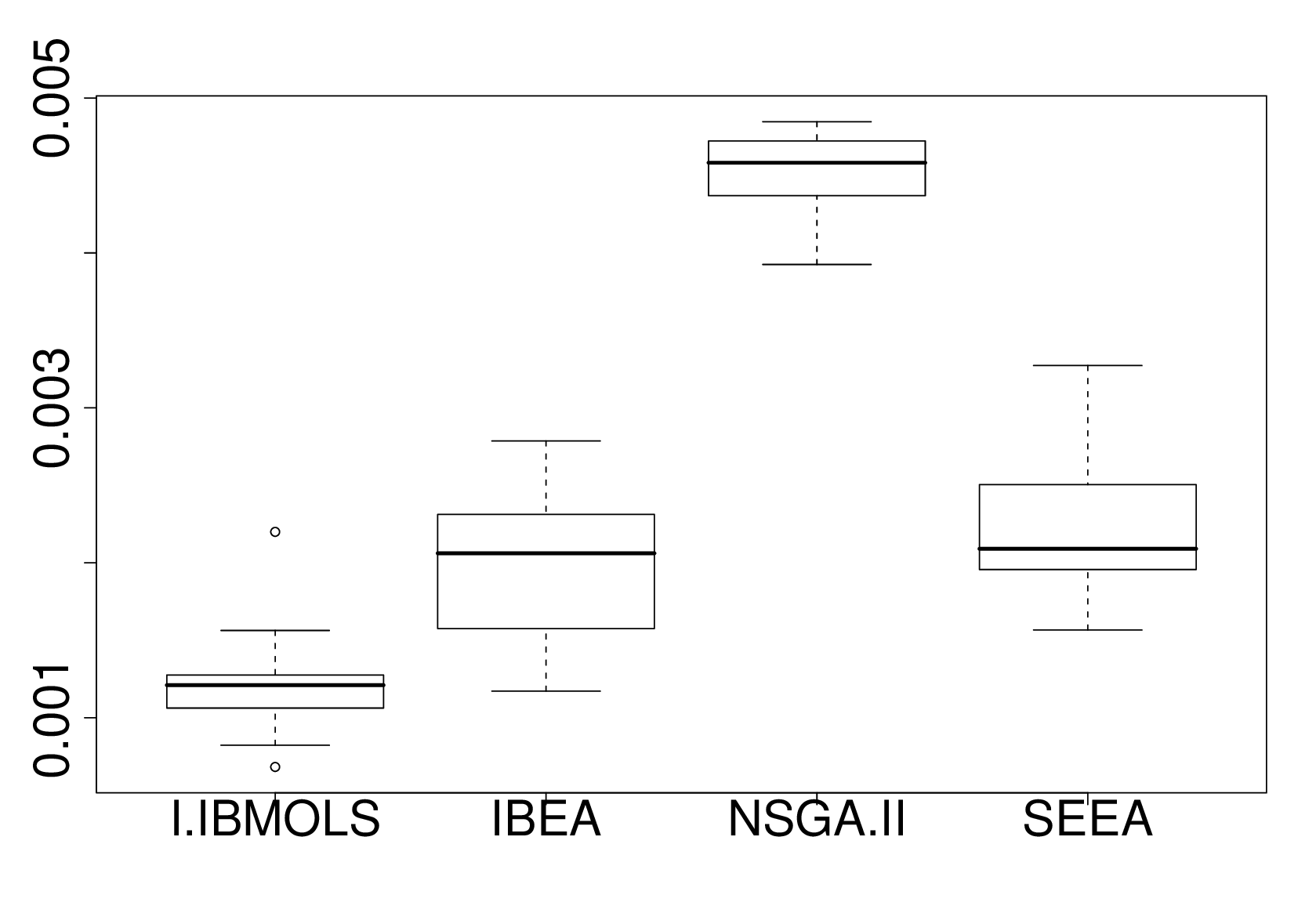}
\end{minipage}
\hfill
\begin{minipage}[c]{0.44\linewidth}
\centering
\includegraphics[angle=270,width=2in]{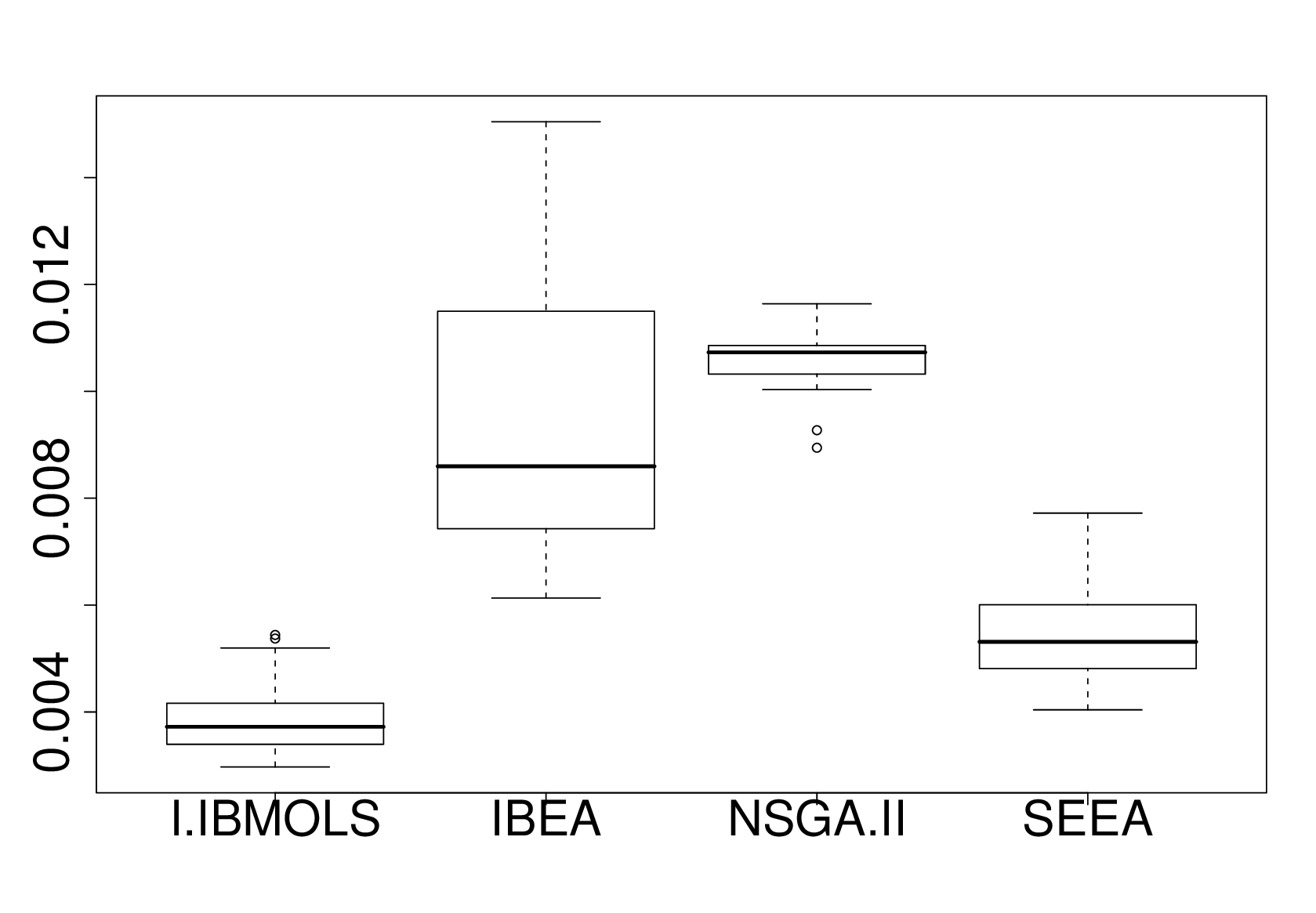}
\end{minipage}
\caption{Performance comparison for I-IBMOLS, IBEA, NSGA-II and SEEA according to the I$^-_H$ and the I$^1_{\epsilon+}$ metric (2).}
\label{fig:box_meta2}
\end{figure}

Previous experiments allowed us to determine satisfying parameters for each one of the proposed algorithms.
The parameter setting is summarized in Table~\ref{tab:param_compare}.
%
\begin{table}[htbp]
\centering
\caption{Parameter setting used to compare search methods.}
\hspace{500pt} \small \setlength{\tabcolsep}{3pt}
\begin{tabular}{cr||c|c||c||c||cc}
&                     &  \multicolumn{2}{c||}{I-IBMOLS}  &  IBEA       &  NSGA-II    &  SEEA      &  \\
&  Instance           &  pop. size  &  noise rate        &  pop. size  &  pop. size  &  pop.size  &  \\
\hline
& {\itshape eil51}    &  M = 20  &  10\%  &  100  &  100  &  100  &  \\
& {\itshape st70}     &  M = 20  &  10\%  &   50  &  100  &  100  &  \\
& {\itshape kroA100}  &  M = 30  &  10\%  &  100  &  200  &  100  &  \\
& {\itshape bier127}  &  M = 30  &  10\%  &  100  &  200  &  100  &  \\
& {\itshape kroA150}  &  M = 30  &  10\%  &  200  &  200  &  100  &  \\
& {\itshape kroA200}  &  M = 30  &  10\%  &  200  &  200  &  100  &  \\
& {\itshape pr264}    &  S = 30  &  20\%  &   50  &  200  &  100  &  \\
& {\itshape pr299}    &  M = 50  &  10\%  &   50  &  200  &  100  &  \\
\end{tabular}
\label{tab:param_compare}
\end{table}
In this section, experiments were carried out with the aim of comparing the results obtained by the different metaheuristics to solve the problem under consideration.
A further analysis of these results was conducted in comparison to the previous ones.
Then, Table~\ref{tab:meta_hyp} and Table~\ref{tab:meta_eps} present the results we respectively obtained for the I$^-_H$ and the I$^1_{\epsilon+}$ metrics.
These tables show, for every test instance, the outcome of the statistical test as well as the significance level given by the p-value.
With respect to the I$^-_H$ metric, we can see that I-IBMOLS is never outperformed by any other search method for all the instances we tested.
The only instance for which I-IBMOLS does not perform better than the other algorithms is {\itshape pr264}, 
where there is no significant difference between the results it obtained and the ones of SEEA.
Furthermore, according to the same metric, SEEA seems overall more efficient than IBEA and NSGA-II.
Indeed, it always performs significantly better than NSGA-II and is statistically outperformed by IBEA only on the {\itshape bier127} and the {\itshape kroA200} instances.
Lastly, IBEA performs significantly better than NSGA-II on every test instance.
Likewise, as shown in Table~\ref{tab:meta_eps}, the results we obtained for the I$^1_{\epsilon+}$ metric,
reveal that I-IBMOLS is globally more efficient than SEEA, IBEA and NSGA-II on most test instances.
Nevertheless, it is not the case on the {\itshape st70} and {\itshape kroA100} instances, where there is no significant difference between I-IBMOLS and SEEA,
and on the {\itshape pr264} instance, where I-IBMOLS is statistically outperformed by SEEA.
Moreover, the results of NSGA-II are always worst than the ones of the other methods for all the instances.
At last, SEEA provides better performance than IBEA on every problem, except for the {\itshape kroA200} instance, where the difference between both is not significant,
and for the {\itshape bier127} instance, where IBEA performs better.
Then, to summarize the results given in Table~\ref{tab:meta_hyp} and Table~\ref{tab:meta_eps} for the I$^-_H$ and the I$^1_{\epsilon+}$ metric respectively, 
I-IBMOLS seems to be the overall most efficient optimization method to solve the bi-objective ring star problem.
Then comes SEEA that performs better than any other EA.
Finally, IBEA performs better than NSGA-II, the later being always outperformed by the other search methods.
All these observations are confirmed by the box-plots given in Figure~\ref{fig:box_meta1} and Figure~\ref{fig:box_meta2}.
We can see that the measures obtained by NSGA-II and, to a lesser extent, IBEA are both generally quite far from the ones obtained by IBMOLS and SEEA,
except for some particular instances such as {\itshape bier127} and {\itshape kroA200}.
Moreover, the differences between the measures obtained by I-IBMOLS and SEEA are very thin for most of the cases, even if SEEA is generally outperformed, 
as pointed out above.

\begin{table}[p]
\centering
\caption{
Performance comparison for I-IBMOLS, IBEA, NSGA-II and SEEA according to the I$^-_H$ metric by using a Mann-Whitney statistical test.
The ``p-value'' columns give the p-value of the statistical test.
The ``T'' columns give the outcome of the statistical test for a significance level of $5\%$: 
either the results of the search method located at a specific row are significantly better than those of the search method located at a specific column~($\succ$), 
either they are worse~($\prec$), or there is no significant difference between both~($\equiv$).
}
\hspace{500pt} \small \setlength{\tabcolsep}{1.8pt}
\begin{tabular}{cr|p{0.65in}||p{0.65in}|c||p{0.65in}|c||p{0.65in}|cc}
&  \multicolumn{2}{c||}{}  &  \multicolumn{2}{|c||}{I-IBMOLS}  &  \multicolumn{2}{|c||}{IBEA}  &  \multicolumn{2}{|c}{NSGA-II}  &  \\
&  \multicolumn{2}{c||}{}  &  \centering p-value  &  T         &  \centering p-value  &  T     &  \centering p-value  &  T      &  \\
\hline
& {\itshape eil51}    &  IBEA      &  $4.578 \cdot 10^{-6}$  &  $\prec$   &                         &            &                         &            &  \\
&                     &  NSGA-II   &  $3.151 \cdot 10^{-8}$  &  $\prec$   &  $3.151 \cdot 10^{-8}$  &  $\prec$   &                         &            &  \\
&                     &  SEEA      &  $0.076$                &  $\prec$   &  $1.053 \cdot 10^{-4}$  &  $\succ$   &  $3.151 \cdot 10^{-8}$  &  $\succ$   &  \\
\hline
& {\itshape st70}     &  IBEA      &  $1.781 \cdot 10^{-4}$  &  $\prec$   &                         &            &                         &            &  \\
&                     &  NSGA-II   &  $3.151 \cdot 10^{-8}$  &  $\prec$   &  $3.151 \cdot 10^{-8}$  &  $\prec$   &                         &            &  \\
&                     &  SEEA      &  $0.011$                &  $\prec$   &  $0.111$                &  $\succ$   &  $3.151 \cdot 10^{-8}$  &  $\succ$   &  \\
\hline
& {\itshape kroA100}  &  IBEA      &  $0.042$                &  $\prec$   &                         &            &                         &            &  \\
&                     &  NSGA-II   &  $5.608 \cdot 10^{-7}$  &  $\prec$   &  $5.608 \cdot 10^{-7}$  &  $\prec$   &                         &            &  \\
&                     &  SEEA      &  $0.004$                &  $\prec$   &  $> 5 \%$               &  $\equiv$  &  $3.151 \cdot 10^{-8}$  &  $\succ$   &  \\
\hline
& {\itshape bier127}  &  IBEA      &  $4.832 \cdot 10^{-4}$  &  $\prec$   &                         &            &                         &            &  \\
&                     &  NSGA-II   &  $3.151 \cdot 10^{-8}$  &  $\prec$   &  $3.151 \cdot 10^{-8}$  &  $\prec$   &                         &            &  \\
&                     &  SEEA      &  $3.151 \cdot 10^{-8}$  &  $\prec$   &  $4.889 \cdot 10^{-7}$  &  $\prec$   &  $9.641 \cdot 10^{-7}$  &  $\succ$   &  \\
\hline
& {\itshape kroA150}  &  IBEA      &  $0.003$                &  $\prec$   &                         &            &                         &            &  \\
&                     &  NSGA-II   &  $3.151 \cdot 10^{-8}$  &  $\prec$   &  $3.151 \cdot 10^{-8}$  &  $\prec$   &                         &            &  \\
&                     &  SEEA      &  $0.003$                &  $\prec$   &  $> 5 \%$               &  $\equiv$  &  $8.914 \cdot 10^{-8}$  &  $\succ$   &  \\
\hline
& {\itshape kroA200}  &  IBEA      &  $0.112$                &  $\prec$   &                         &            &                         &            &  \\
&                     &  NSGA-II   &  $3.151 \cdot 10^{-8}$  &  $\prec$   &  $3.151 \cdot 10^{-8}$  &  $\prec$   &                         &            &  \\
&                     &  SEEA      &  $5.483 \cdot 10^{-5}$  &  $\prec$   &  $3.122 \cdot 10^{-5}$  &  $\prec$   &  $3.151 \cdot 10^{-8}$  &  $\succ$   &  \\
\hline
& {\itshape pr264}    &  IBEA      &  $1.084 \cdot 10^{-5}$  &  $\prec$   &                         &            &                         &            &  \\
&                     &  NSGA-II   &  $3.151 \cdot 10^{-8}$  &  $\prec$   &  $3.151 \cdot 10^{-8}$  &  $\prec$   &                         &            &  \\
&                     &  SEEA      &  $> 5 \%$               &  $\equiv$  &  $3.151 \cdot 10^{-8}$  &  $\succ$   &  $3.151 \cdot 10^{-8}$  &  $\succ$   &  \\
\hline
& {\itshape pr299}    &  IBEA      &  $5.877 \cdot 10^{-6}$  &  $\prec$   &                         &            &                         &            &  \\
&                     &  NSGA-II   &  $3.151 \cdot 10^{-8}$  &  $\prec$   &  $3.151 \cdot 10^{-8}$  &  $\prec$   &                         &            &  \\
&                     &  SEEA      &  $2.114 \cdot 10^{-7}$  &  $\prec$   &  $> 5 \%$               &  $\equiv$  &  $3.151 \cdot 10^{-8}$  &  $\succ$   &  \\
\end{tabular}
\label{tab:meta_hyp}
\end{table}

\begin{table}[p]
\centering
\caption{
Performance comparison for I-IBMOLS, IBEA, NSGA-II and SEEA according to the I$^1_{\epsilon+}$ metric by using a Mann-Whitney statistical test.
The ``p-value'' columns give the p-value of the statistical test.
The ``T'' columns give the outcome of the statistical test for a significance level of $5\%$: 
either the results of the search method located at a specific row are significantly better than those of the search method located at a specific column~($\succ$), 
either they are worse~($\prec$), or there is no significant difference between both~($\equiv$).
}
\hspace{500pt} \small \setlength{\tabcolsep}{1.8pt}
\begin{tabular}{cr|p{0.65in}||p{0.65in}|c||p{0.65in}|c||p{0.65in}|cc}
&  \multicolumn{2}{c||}{}  &  \multicolumn{2}{|c||}{I-IBMOLS}  &  \multicolumn{2}{|c||}{IBEA}  &  \multicolumn{2}{|c}{NSGA-II}  &  \\
&  \multicolumn{2}{c||}{}  &  \centering p-value  &  T         &  \centering p-value  &  T     &  \centering p-value  &  T      &  \\
\hline
& {\itshape eil51}    &  IBEA      &  $1.177 \cdot 10^{-6}$  &  $\prec$   &                         &            &                         &            &  \\
&                     &  NSGA-II   &  $3.151 \cdot 10^{-8}$  &  $\prec$   &  $3.663 \cdot 10^{-8}$  &  $\prec$   &                         &            &  \\
&                     &  SEEA      &  $0.117$                &  $\prec$   &  $9.604 \cdot 10^{-6}$  &  $\succ$   &  $3.151 \cdot 10^{-8}$  &  $\succ$   &  \\
\hline
& {\itshape st70}     &  IBEA      &  $1.084 \cdot 10^{-5}$  &  $\prec$   &                         &            &                         &            &  \\
&                     &  NSGA-II   &  $3.658 \cdot 10^{-8}$  &  $\prec$   &  $9.595 \cdot 10^{-6}$  &  $\prec$   &                         &            &  \\
&                     &  SEEA      &  $> 5 \%$               &  $\equiv$  &  $5.733 \cdot 10^{-8}$  &  $\succ$   &  $3.146 \cdot 10^{-8}$  &  $\succ$   &  \\
\hline
& {\itshape kroA100}  &  IBEA      &  $0.005$                &  $\prec$   &                         &            &                         &            &  \\
&                     &  NSGA-II   &  $7.363 \cdot 10^{-7}$  &  $\prec$   &  $1.846 \cdot 10^{-6}$  &  $\prec$   &                         &            &  \\
&                     &  SEEA      &  $> 5 \%$               &  $\equiv$  &  $0.027$                &  $\succ$   &  $1.192 \cdot 10^{-7}$  &  $\succ$   &  \\
\hline
& {\itshape bier127}  &  IBEA      &  $8.497 \cdot 10^{-5}$  &  $\prec$   &                         &            &                         &            &  \\
&                     &  NSGA-II   &  $3.151 \cdot 10^{-8}$  &  $\prec$   &  $4.941 \cdot 10^{-8}$  &  $\prec$   &                         &            &  \\
&                     &  SEEA      &  $3.151 \cdot 10^{-8}$  &  $\prec$   &  $1.377 \cdot 10^{-7}$  &  $\prec$   &  $0.008$                &  $\succ$   &  \\
\hline
& {\itshape kroA150}  &  IBEA      &  $2.947 \cdot 10^{-6}$  &  $\prec$   &                         &            &                         &            &  \\
&                     &  NSGA-II   &  $1.192 \cdot 10^{-7}$  &  $\prec$   &  $1.390 \cdot 10^{-4}$  &  $\prec$   &                         &            &  \\
&                     &  SEEA      &  $0.003$                &  $\prec$   &  $4.898 \cdot 10^{-6}$  &  $\succ$   &  $5.608 \cdot 10^{-7}$  &  $\succ$   &  \\
\hline
& {\itshape kroA200}  &  IBEA      &  $0.037$                &  $\prec$   &                         &            &                         &            &  \\
&                     &  NSGA-II   &  $5.733 \cdot 10^{-8}$  &  $\prec$   &  $3.146 \cdot 10^{-8}$  &  $\prec$   &                         &            &  \\
&                     &  SEEA      &  $0.037$                &  $\prec$   &  $> 5\%$                &  $\equiv$  &  $3.151 \cdot 10^{-8}$  &  $\succ$   &  \\
\hline
& {\itshape pr264}    &  IBEA      &  $8.902 \cdot 10^{-8}$  &  $\prec$   &                         &            &                         &            &  \\
&                     &  NSGA-II   &  $4.256 \cdot 10^{-8}$  &  $\prec$   &  $0.084$                &  $\succ$   &                         &            &  \\
&                     &  SEEA      &  $1.745 \cdot 10^{-5}$  &  $\succ$   &  $3.142 \cdot 10^{-8}$  &  $\succ$   &  $3.146 \cdot 10^{-8}$  &  $\succ$   &  \\
\hline
& {\itshape pr299}    &  IBEA      &  $3.151 \cdot 10^{-8}$  &  $\prec$   &                         &            &                         &            &  \\
&                     &  NSGA-II   &  $3.151 \cdot 10^{-8}$  &  $\prec$   &  $0.0495$               &  $\prec$   &                         &            &  \\
&                     &  SEEA      &  $5.877 \cdot 10^{-6}$  &  $\prec$   &  $3.225 \cdot 10^{-7}$  &  $\succ$   &  $3.151 \cdot 10^{-8}$  &  $\succ$   &  \\
\end{tabular}
\label{tab:meta_eps}
\end{table}

\begin{figure}[p]
\centering
\includegraphics[width=4.5in]{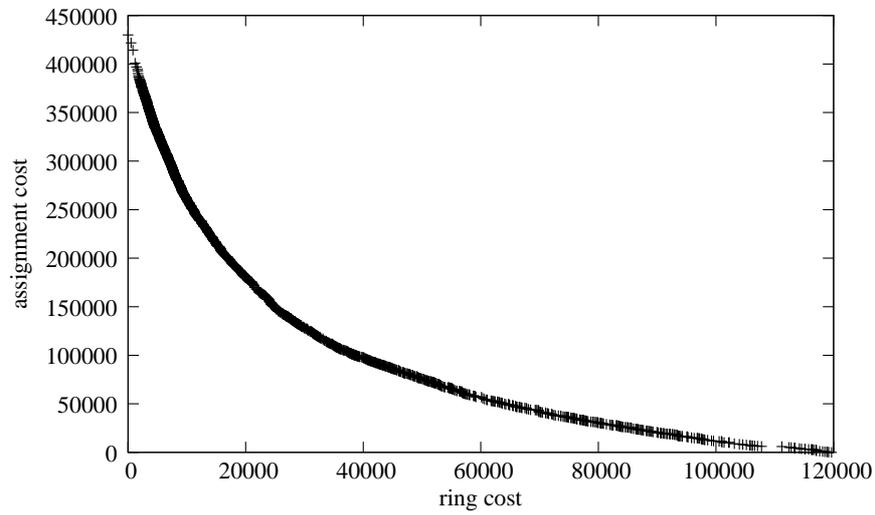}
\caption{Best set of non-dominated points found for the {\itshape bier127} test instance.}
\label{fig:front.bier127}
\end{figure}

\begin{figure}[p]
\centering
\includegraphics[width=4.5in]{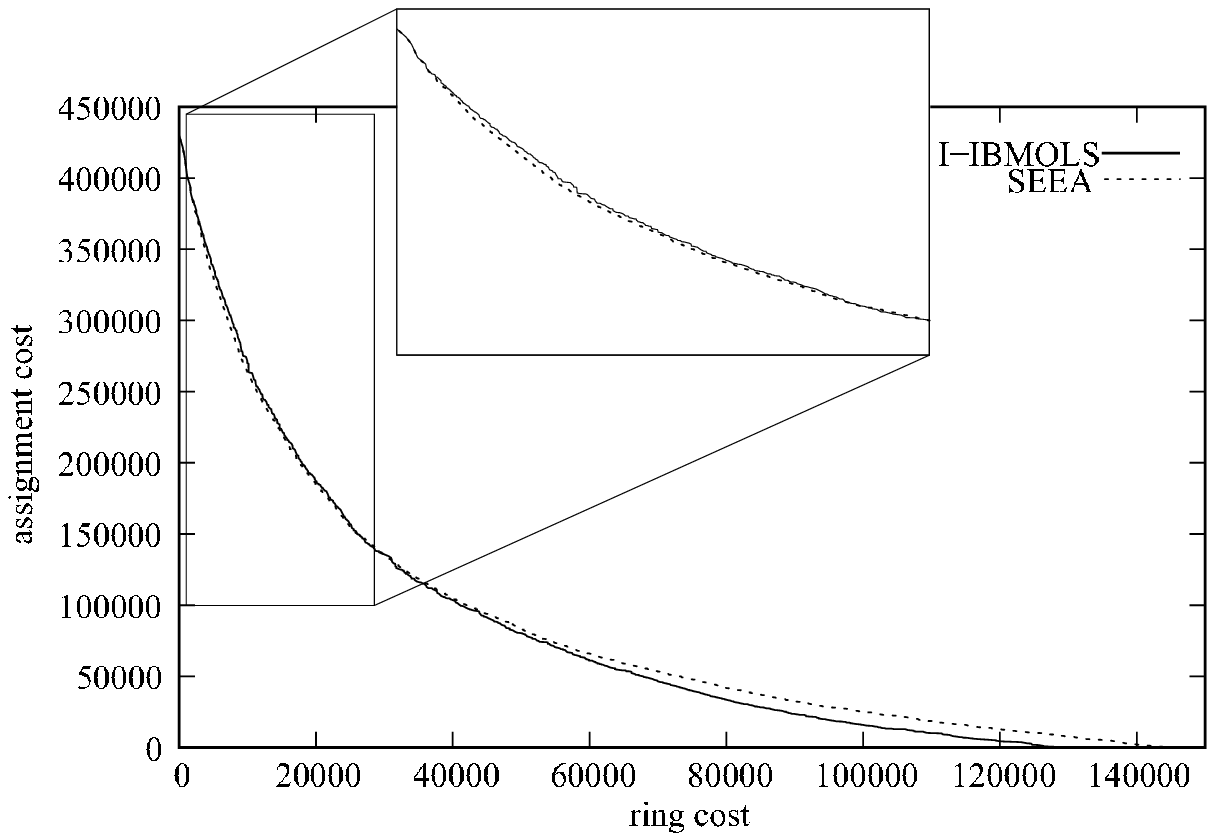}
\caption{$90\%$-attainment surface plot obtained by the approximation sets found by I-IBMOLS and SEEA for the {\itshape bier127} test instance.}
\label{fig:eaf.bier127}
\end{figure}

\begin{figure}[htbp]
\centering
\includegraphics[width=4.5in]{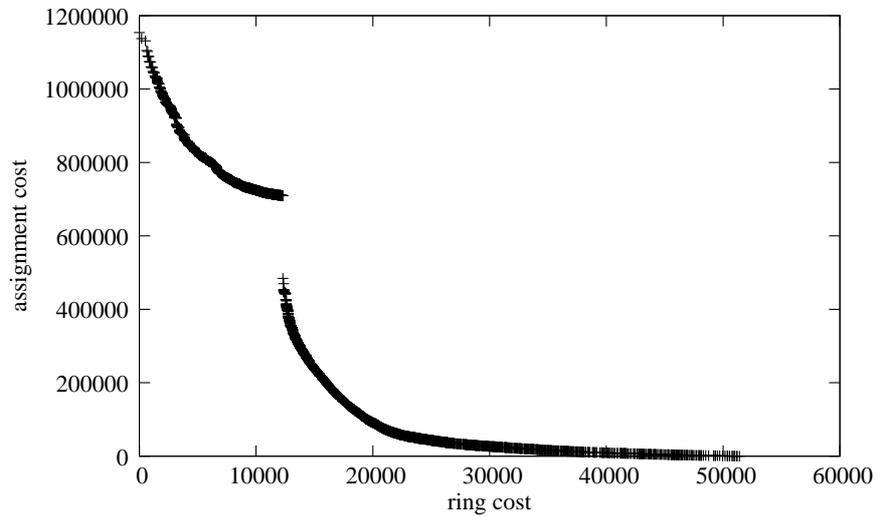}
\caption{Best set of non-dominated points found for the {\itshape pr264} test instance.}
\label{fig:front.pr264}
\end{figure}

\begin{figure}[htbp]
\centering
\includegraphics[width=4.5in]{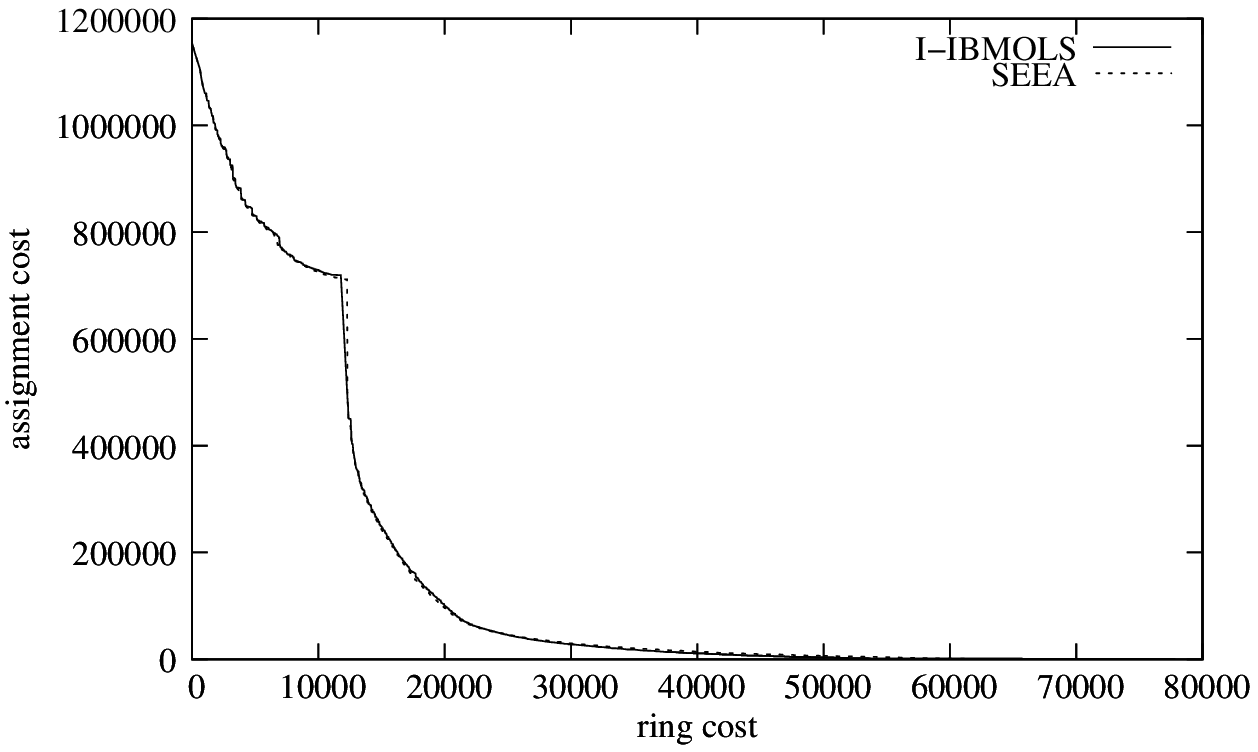}
\caption{$90\%$-attainment surface plot obtained by the approximation sets found by I-IBMOLS and SEEA for the {\itshape pr264} test instance.}
\label{fig:eaf.pr264}
\end{figure}

In order to study the differences between I-IBMOLS and SEEA more finely, 
examples of the empirical attainment function obtained by these two algorithms are represented in Figure~\ref{fig:eaf.bier127} and Figure~\ref{fig:eaf.pr264} 
for the {\itshape bier127} and the {\itshape pr264} instances respectively.
Theses figures illustrate the limit of the objective space that is attained by at least $90 \%$ of the runs by both methods.
The reader is referred to~\cite{FGP:05} for more details about the concept of empirical attainment function.
For the {\itshape bier127} instance, where the results obtained by I-IBMOLS are much better than the ones of SEEA with respect to both metrics, 
we can see that SEEA seems to have more trouble in finding solutions having a low ring cost, 
whereas it seems slightly better than I-IBMOLS for finding interesting solutions with a small assignment cost.
For the {\itshape pr264} instance, where SEEA is overall more efficient, it is less obvious to observe any superiority graphically.
Such a result mainly seems to be due to the discontinuity of the trade-off surface, for which a good approximation is given in Figure~\ref{fig:front.pr264}.
An example of a more often encountered Pareto front is illustrated in Figure~\ref{fig:front.bier127} for the {\itshape bier127} instance.

One of the main characteristics of the problem under consideration seems to be the high number of points located in the trade-off surface, 
as shown in Figure~\ref{fig:front.bier127} and Figure~\ref{fig:front.pr264}.
Then, after a certain number of iterations, a large part of the population involved in NSGA-II, IBEA and I-IBMOLS might map to potentially non-dominated points.
This could explain the low efficiency of NSGA-II.
Since the same rank is assigned to the major part of the population, only the crowding distance is used to compare solutions.
However, it is not the case in SEEA because all non-dominated solutions contained in the archive can potentially take part in the evolution engine.
Moreover, the indicator-based fitness assignment scheme is obviously much more suited to discriminate potentially efficient solutions than the single crowding distance.
Moreover, the high performance of I-IBMOLS in comparison to IBEA might depends on how close are the solutions which map to non-dominated points in the decision space.
If these solutions are close to each other, according to the neighborhood operators, 
a LS is known to be particularly well-suited to find additional interesting solutions by exploring the neighborhood of a potentially efficient solution.
On the contrary, an EA usually explores the decision space in a more random way.
Thus, a landscape analysis could be interesting to study the bi-objective RSP in more depth.
Besides, the results obtained by SEEA are quite satisfying in contrast to those obtained by the other EAs.
Indeed, NSGA-II and IBEA are two state-of-art search methods for multi-objective optimization and are widely used in the community.
Nevertheless, despite its simplicity, SEEA seems to be much more efficient for solving the B-RSP in a given amount of time.
It could then be interesting to experiment the ability of this method for the resolution of other kinds of multi-objective combinatorial optimization problems.

As a next step, we will try to solve problem instances involving an even bigger number of nodes to verify if our observations are still valid.
It could then be interesting to design a cooperation scheme between two different methods 
({\itshape i.e.} the local search procedure and an evolutionary algorithm) in order to benefit from the respective quality of each one of them.
The resulting hybrid metaheuristic could be particularly effective for solving large size problems.

\section{A Cooperative Approach}
\label{sec:hybrid}

This section presents a cooperative approach combining SEEA and the non-iterative version of IBMOLS in order to solve the bi-objective RSP. 
Two variants are proposed: a periodic one that operates a systematic cooperation and an auto-adaptive one that decides on-line when the cooperation must occur.

\subsection{Motivation}

Designing metaheuristics for solving combinatorial optimization problems is generally a matter of intensification and diversification.
This is even more pronounced for MOPs where the goal is to find a well-converged and well-diversified efficient set approximation.
However, LSs are known to be particularly efficient as intensifying methods 
whereas EAs are clearly powerful to explore the decision space thanks to their variation operators.
Instead of trying to improve one method in term of diversification or the other in term of intensification, 
a common approach is to hybridize both in order to make them cooperate and then to benefit of their respective behaviors.
Thus, hybrid metaheuristics have shown their efficiency to solve different kinds of optimization problems~\cite{Tal:02}, including MOPs~\cite{EG:08}.
They generally consist of a search method that cooperates with a second one in order to improve a solution or a set of solutions.

In the previous section, 
we saw that SEEA and I-IBMOLS were overall more efficient than NSGA-II and IBEA to approximate the efficient set for the problem under consideration.
Moreover, these two methods are quite different to each other and do not explore the search space in the same way.
SEEA has been conceived in order to find a rough approximation of the Pareto set in a very short amount of time 
whereas the non-iterative version of IBMOLS is able to improve an approximated set in a very efficient way.
It could then be interesting to design a cooperation scheme between these two algorithms.
The resulting hybrid metaheuristic could be particularly efficient for solving large size problems.
Furthermore, both methods maintain a secondary population (the archive) in parallel of the main population to store non-dominated solutions.
This archive is not only used as an external storage (what is the case, for instance, in IBEA and NSGA-II),
but also takes part in the evolution engine as it serves to build new solutions to explore.
Thus, each method can manage its own population and therefore use the archive as a single shared memory.

\subsection{Cooperative Schemes}

\begin{figure}[htbp]
\begin{minipage}[b]{0.49\linewidth}
\centering
\includegraphics[width=2.3in]{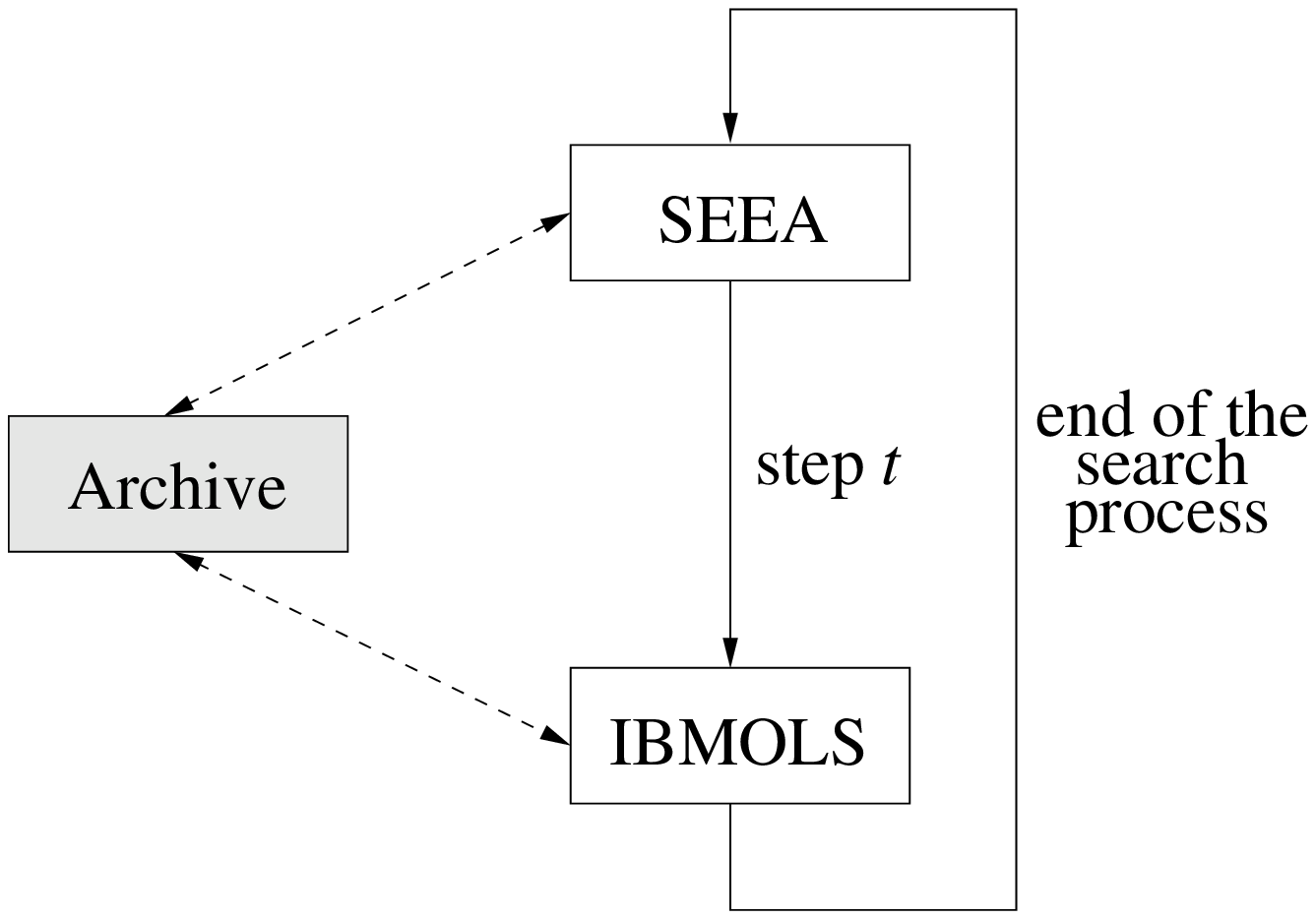}
\caption{Illustration of the Periodic Cooperative Search (PCS).}
\label{fig:pcs}
\end{minipage}
\hfill
\begin{minipage}[b]{0.49\linewidth}
\centering
\includegraphics[width=2.3in]{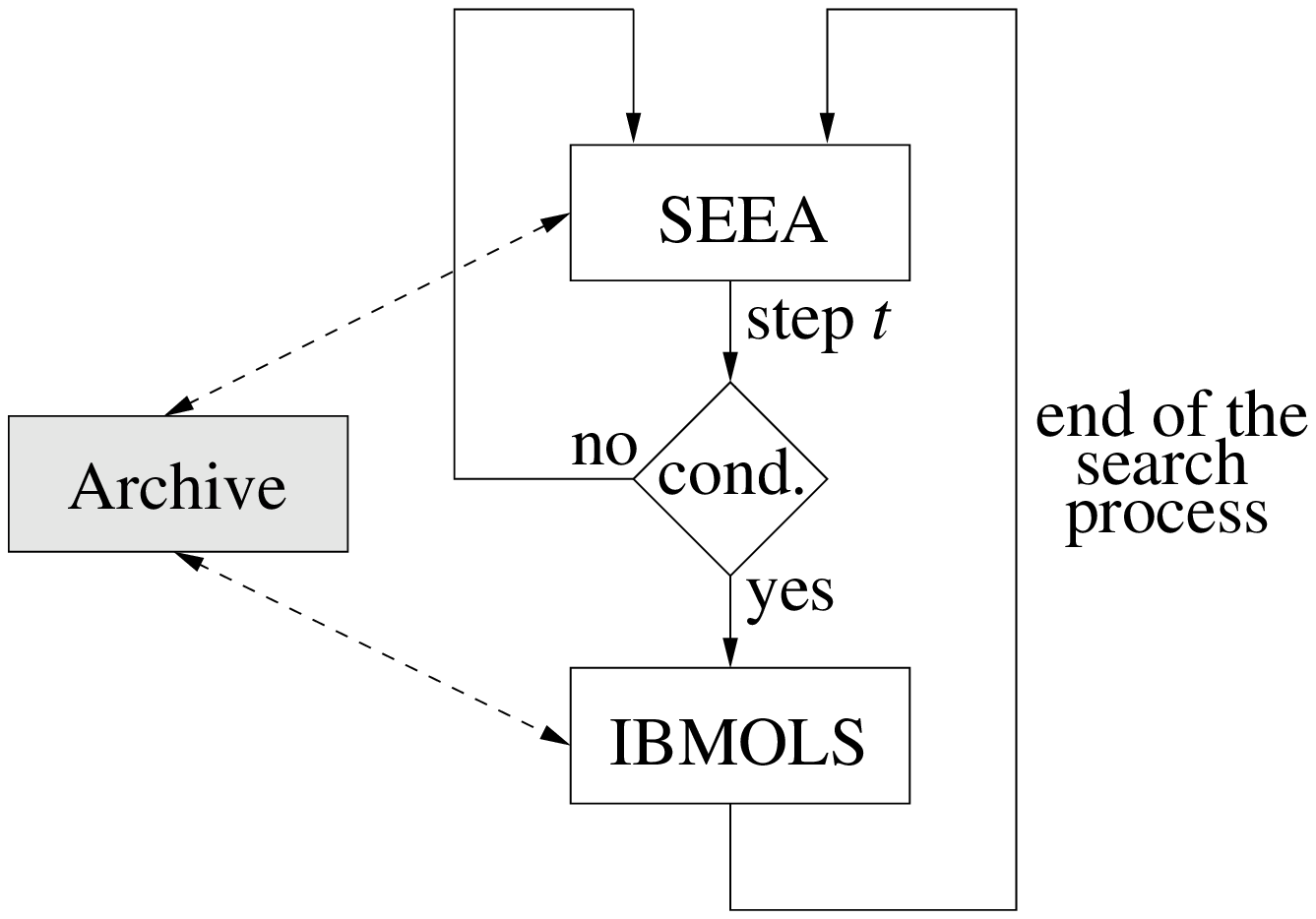}
\caption{Illustration of the Auto-adaptive Cooperative Search (ACS).}
\label{fig:acs}
\end{minipage}
\end{figure}

The general idea of our hybridization scheme is to run SEEA and to lunch IBMOLS regularly by using a subset of archive items as an initial population.
As the non-iterative version of IBMOLS stops when its own archive does not receive any new efficient solution anymore, 
we can restart the SEEA process until the next step of the hybrid algorithm.
A step of the hybrid metaheuristic can, for instance, be defined by a certain amount of time or by a certain number of generations.
Besides, as SEEA uses the non-dominated solutions found by IBMOLS to create new ones and vice versa,
the global archive is the only memory shared by the two search agents to exchange information.
Resulting from this, we can imagine two versions of the hybrid algorithm: 
($i$)~a \emph{periodic version}, in which IBMOLS is launched at each step, 
and ($ii$)~an \emph{auto-adaptive version}, in which IBMOLS is launched at a specific step only if a condition is verified.
These two approaches will be denoted by PCS (for \emph{Periodic Cooperative Search}) and ACS (for \emph{Auto-adaptive Cooperative Search}) in the remainder of the paper, 
and are respectively illustrated in Figure~\ref{fig:pcs} and Figure~\ref{fig:acs}.

The ACS method decides by itself, and on-line, if it would be interesting to launch IBMOLS at a specific instant of the search process.
The condition that we consider here to start IBMOLS during a particular step of the ACS 
is that the archive of potentially efficient solutions does not improve enough regarding to the optimization scenario.
One possibility is to measure the quality of the current archive~$A^t$ in comparison to the one of the previous step~$A^{t-1}$.
Different metrics exist to evaluate the convergence of an approximated efficient set regarding to another.
For instance, we could have computed the hypervolume~\cite{ZT:99} of $A^t$ and $A^{t-1}$ and measure the difference between both.
But, even if it evaluates the quality of an approximated efficient set both in term of convergence and diversity,
the hypervolume metric has the drawback of being time consuming.
Thus, let us introduce the \emph{contribution} metric $C$ proposed by Meunier et al.~\cite{MTR:00}.
This metric gives an idea of the quality of an approximated efficient set~$S_1$ compared to another one $S_2$ in term of convergence 
and can be computed in a reasonable running time.
If $S^\star$ denotes the set of non-dominated solutions given by $S_1 \cup S_2$,
the contribution of $S_1$ on $S_2$ evaluates the proportion of non-dominated solutions represented by $S_1$ in $S^\star$.
Then, if $f(S_1) = f(S_2)$, $C(S_1,S_2) = C(S_2,S_1) = 0.5$.
If every solution of $S_1$ is dominated by at least one solution of $S_2$, $C(S_1,S_2) = 0$. 
And, generally speaking, $C(S_1,S_2) + C(S_2,S_1) = 1$.
In our case, at each step~$t$ of the ACS, we compute the contribution of the current archive~$A^t$ on the archive of the previous step~$A^{t-1}$.
Thus, as non-dominated solutions are not lost between two steps, $A^t$ is at least as good as $A^{t+1}$ 
and the set of non-dominated solutions of $A^t \cup A^{t-1}$ is then $A^t$.
As a consequence, we know that $C(A^t,A^{t-1}) \in [0.5,1]$.
Assuming that the archive does not improve enough if the contribution of $A^t$ on $A^{t-1}$ is less than a pre-defined threshold $\delta \in [0.5,1]$,
we choose to launch IBMOLS only if $C(A^t,A^{t-1}) \leq \delta$.
Let us remark that an ACS with a $\delta = 0.5$ would result in SEEA, and that an ACS with a $\delta = 1.0$ would result in PCS, 
the time spent in calculating the different contribution values in less.


\subsection{Related Works}

Different schemes exist on how two search methods can be combined.
According to the taxonomy proposed in~\cite{Tal:02}, two levels (low and high) and two modes (relay and cooperative) of hybridization can be distinguished.
In \emph{low-level} hybridization, a given function of a metaheuristic is replaced by another metaheuristic. 
In \emph{high-level} hybridization, there is no direct relationship between the internal workings of a metaheuristic. 
A second hierarchical classification deals with the mode of cooperation.
The \emph{relay} mode consists of a method applied after another one in a pipeline way, the last one using the output of the previous one as its input.
On the contrary, \emph{teamwork} hybridization represents a cooperative optimization model in which search agents exchange information with the others.
Here, we consider that the optimization process alternates between two search algorithms: SEEA and the non-iterative version of IBMOLS.
Both methods are self-contained and cooperates with each other via a global archive.
Therefore, the resulting hybrid metaheuristic can be classified on the \emph{High-level Teamwork Hybrid} (HTH) class of the taxonomy introduced in~\cite{Tal:02} 
and can then be denoted by HTH(SEEA+IBMOLS).
In their survey on hybrid metaheuristics to solve multi-objective combinatorial optimization problems~\cite{EG:08}, 
Ehrgott and Gandibleux identify three categories of methods hybridizing an EA with a neighborhood search algorithm:
($i$)~an hybridization to make a method more aggressive, ($ii$)~an hybridization to drive a method and ($iii$)~an hybridization for exploiting complementary strengths.
The last one consists of alternating between both search methods, what is the case within our hybridization.
But, most of existing approaches occur in a pipeline way, like in the relay mode.
The few teamwork hybridization techniques of the literature to solve MOPs are most often compound of the same kind of metaheuristics, 
what is generally the case in the island model.
Moreover, this class of hybrid methods often uses a neighborhood search method combining a set of scalarizing functions.
The originality of the approach proposed in this paper is that search agents are based on different types of metaheuristics, are hybridized in a teamwork mode,
and do not use any scalar approach to convert the multiple objective functions into a single one.
Furthermore, one of the variant we propose, the auto-adaptive one, automatically detects when to start the neighborhood search according to the optimization scenario.

\subsection{Experiments}

In order to experiment the efficiency of our two cooperative approaches, 
we compare them to SEEA and to the iterative version of IBMOLS by using the same experimental protocol as the one defined in Section~\ref{sec:protocol}.
After having set the parameters for each hybrid search method, 
we give some computational results and we discuss the contribution of the hybridization for the resolution of the B-RSP.
Finally, we give an idea of the behavior of our approach compared to a state-of-art single-objective approach proposed in~\cite{LL+:04}.

\subsubsection{Parameter Setting}
For both cooperative search methods, the population size managed by SEEA is set to $100$, 
and the population size managed by IBMOLS is set according to the instance under consideration.
These sizes have been set on the same way that we did for the stand-alone iterative version of IBMOLS on the previous section (see Table~\ref{tab:param}).
But, for large-size instances, initial experiments were not satisfying since 
the hybrid algorithms were generally not able to launch IBMOLS more than once during the search process as it was too much time consuming.
For this reason, we bounded the IBMOLS population size to $30$. 
Then, a IBMOLS population of $20$ individuals has been set for instances with less than $100$ nodes, 
and a IBMOLS population of $30$ individuals has been set for instances with $100$ nodes and more.
The step $t$ of the hybrid algorithms has been set to $0.5\%$ of the total run time for the instance under consideration, rounded to the nearest upper integer.
Finally, we investigated different $\delta$ values for the ACS method: $0.6$, $0.7$, $0.8$ and $0.9$.
As pointed out above, an ACS with a $\delta$ value of $0.5$ and $1.0$ is similar to PCS and SEEA, respectively.

\subsubsection{Computational Results and Discussion}

\paragraph{ACS Parameter Analysis.}
As a first step, we analyze the influence of the $\delta$ parameter on the results of the ACS method.
Results are summarized in Table~\ref{tab:acs}.
A value of $0.6$ is generally less efficient than at least one other setting on most test instances according to both metrics.
For the other values, the results are roughly comparable, 
even if a $\delta = 0.8$ is never statistically outperformed by any other value according to I$^-_H$ and I$^1_{\epsilon+}$.
Therefore, we will use this parameter setting for the next experiments.
\begin{table}[htbp]
\centering
\caption{
Performance comparison of different $\delta$ parameter settings for the ACS method according to the I$^-_H$ and the I$^1_{\epsilon+}$ metrics 
by using a Mann-Whitney statistical test with a significance level of $5\%$.
Each column gives the number of algorithms with another parameter setting by which it is statistically outperformed for the instance under consideration.
}
\hspace{500pt} \small \setlength{\tabcolsep}{3pt}
\begin{tabular}{cr||p{0.25in}|p{0.25in}|p{0.25in}|p{0.25in}||p{0.25in}|p{0.25in}|p{0.25in}|p{0.25in}c}
&            &  \multicolumn{4}{|c||}{I$^-_H$}                                       &  \multicolumn{4}{|c}{I$^1_{\epsilon+}$}                              &  \\
&  $\delta$  &  \centering 0.60  &  \centering 0.70  &  \centering 0.80  &  \centering 0.90  &  \centering 0.60  &  \centering 0.70 &  \centering 0.80  &  \centering 0.90  &  \\
\hline
& {\itshape eil51}    &  \centering  {\bfseries 0}  &  \centering  {\bfseries 0}  &  \centering  {\bfseries 0}  &  \centering  {\bfseries 0}
											&  \centering  {\bfseries 0}  &  \centering  {\bfseries 0}  &  \centering  {\bfseries 0}  &  \centering  {\bfseries 0}  &  \\
& {\itshape st70}     &  \centering  {\bfseries 0}  &  \centering  {\bfseries 0}  &  \centering  {\bfseries 0}  &  \centering  {\bfseries 0}
											&  \centering  {\bfseries 0}  &  \centering  {\bfseries 0}  &  \centering  {\bfseries 0}  &  \centering  {\bfseries 0}  &  \\
& {\itshape kroA100}  &  \centering  1              &  \centering  {\bfseries 0}  &  \centering  {\bfseries 0}  &  \centering  {\bfseries 0}
											&  \centering  2              &  \centering  1              &  \centering  {\bfseries 0}  &  \centering  {\bfseries 0}  &  \\
& {\itshape bier127}  &  \centering  {\bfseries 0}  &  \centering  {\bfseries 0}  &  \centering  {\bfseries 0}  &  \centering  1
											&  \centering  3              &  \centering  {\bfseries 0}  &  \centering  {\bfseries 0}  &  \centering  {\bfseries 0}  &  \\
& {\itshape kroA150}  &  \centering  2              &  \centering  1              &  \centering  {\bfseries 0}  &  \centering  {\bfseries 0}
											&  \centering  3              &  \centering  2              &  \centering  {\bfseries 0}  &  \centering  {\bfseries 0}  &  \\
& {\itshape kroA200}  &  \centering  3              &  \centering  {\bfseries 0}  &  \centering  {\bfseries 0}  &  \centering  {\bfseries 0}
											&  \centering  3              &  \centering  2              &  \centering  {\bfseries 0}  &  \centering  {\bfseries 0}  &  \\
& {\itshape pr264}    &  \centering  {\bfseries 0}  &  \centering  {\bfseries 0}  &  \centering  {\bfseries 0}  &  \centering  {\bfseries 0}
											&  \centering  {\bfseries 0}  &  \centering  {\bfseries 0}  &  \centering  {\bfseries 0}  &  \centering  {\bfseries 0}  &  \\
& {\itshape pr299}    &  \centering  3              &  \centering  {\bfseries 0}  &  \centering  {\bfseries 0}  &  \centering  {\bfseries 0}
											&  \centering  3              &  \centering  {\bfseries 0}  &  \centering  {\bfseries 0}  &  \centering  {\bfseries 0}  &  \\
\end{tabular}
\label{tab:acs}
\end{table}

\paragraph{Comparison between Search Methods.}
In addition to the benchmark test instances we investigated previously, we experimented two bigger ones: {\itshape pr439} and {\itshape pr1002}.
The maximum running time for these new instances has been set to $70$ and $100$ minutes, respectively.
Following the general trend identified in Section~\ref{sec:ibmols_study}, we set the instance-specific parameters for the iterative version of IBMOLS as follows:
a population of $70$ and $100$ individuals respectively, and a noise rate of $10\%$.
The other parameters were set in the same way than for other instances.

\begin{figure}[p]
\small 
\begin{minipage}[c]{0.1\linewidth}
Instance
\end{minipage}
\begin{minipage}[c]{0.44\linewidth}
\centering I$^-_H$
\end{minipage}
\hfill
\begin{minipage}[c]{0.44\linewidth}
\centering I$^1_{\epsilon+}$
\end{minipage}
\begin{minipage}[c]{0.1\linewidth}
\centering
{\itshape eil51}
\end{minipage}
\begin{minipage}[c]{0.44\linewidth}
\centering
\includegraphics[angle=270,width=2in]{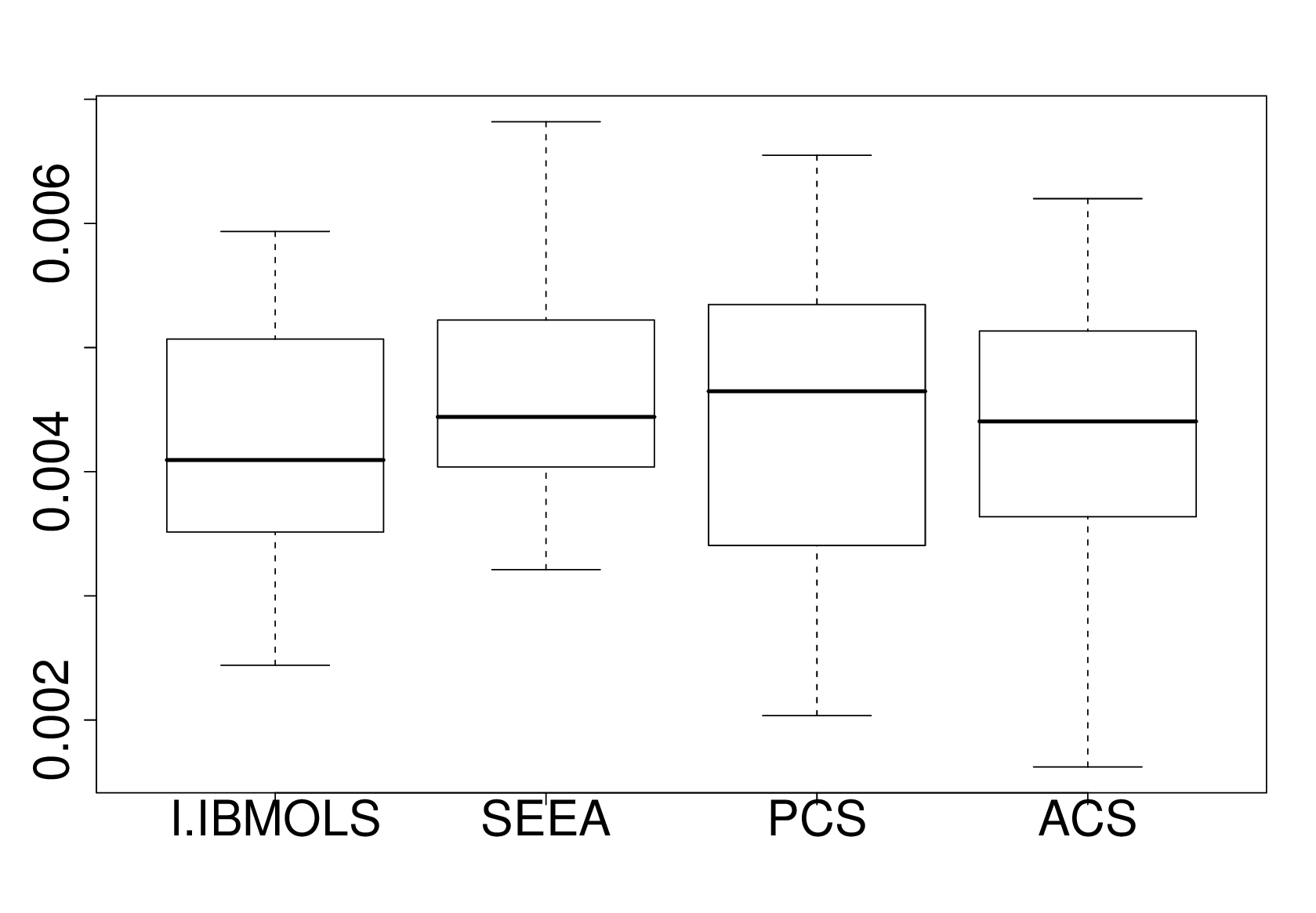}
\end{minipage}
\hfill
\begin{minipage}[c]{0.44\linewidth}
\centering
\includegraphics[angle=270,width=2in]{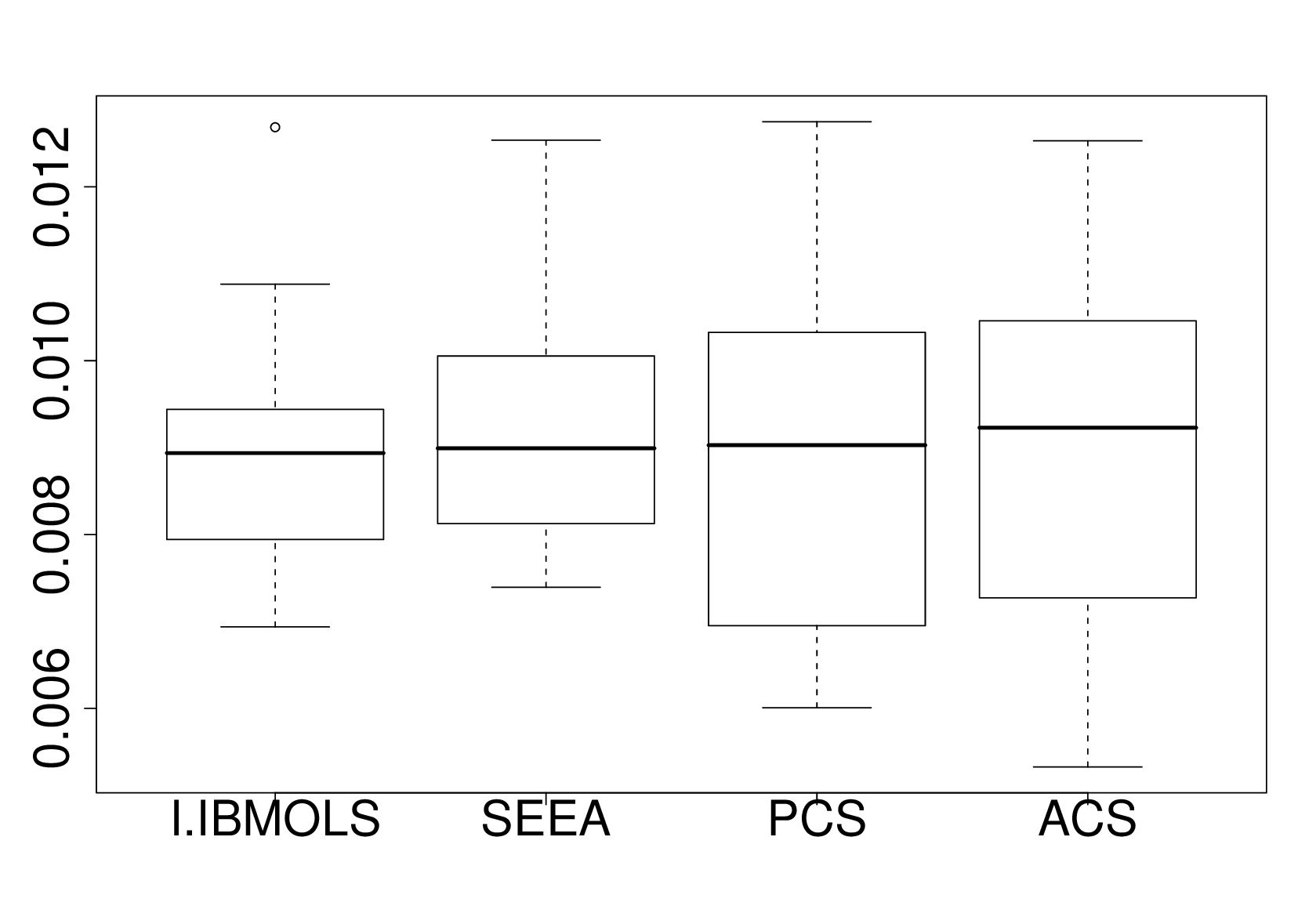}
\end{minipage}
\begin{minipage}[c]{0.1\linewidth}
\centering
{\itshape st70}
\end{minipage}
\begin{minipage}[c]{0.44\linewidth}
\centering
\includegraphics[angle=270,width=2in]{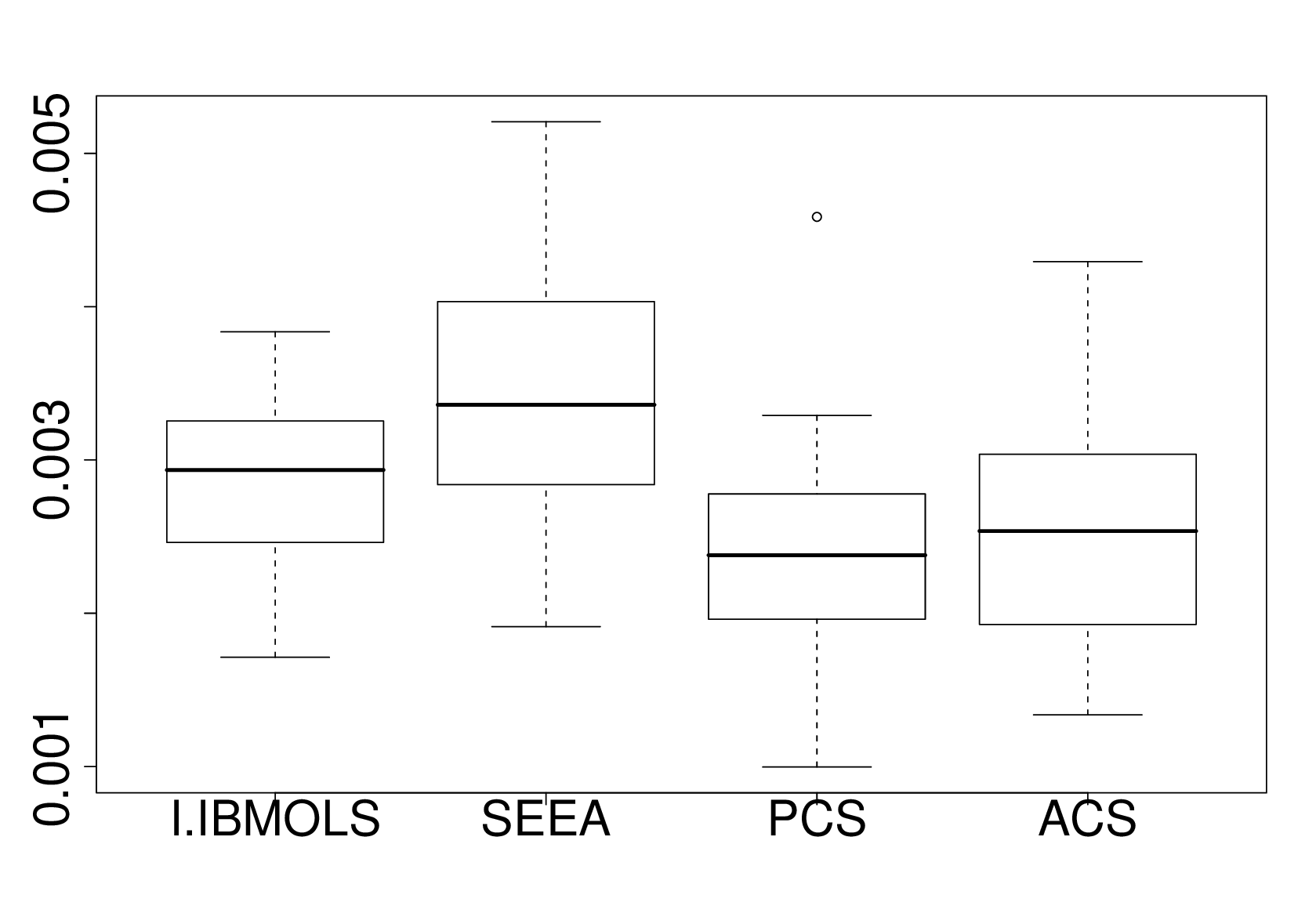}
\end{minipage}
\hfill
\begin{minipage}[c]{0.44\linewidth}
\centering
\includegraphics[angle=270,width=2in]{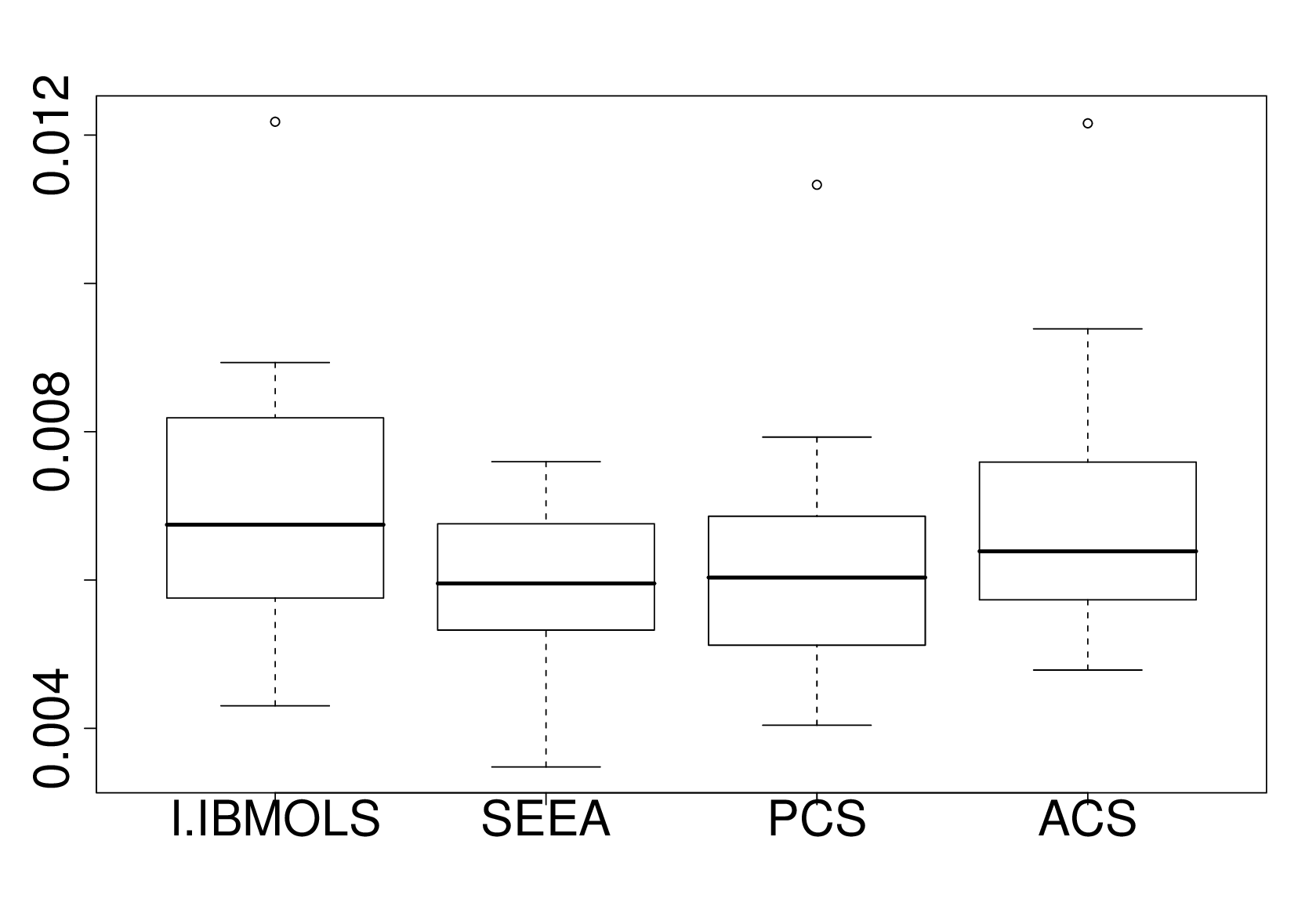}
\end{minipage}
\begin{minipage}[c]{0.1\linewidth}
\centering
{\itshape kroA100}
\end{minipage}
\begin{minipage}[c]{0.44\linewidth}
\centering
\includegraphics[angle=270,width=2in]{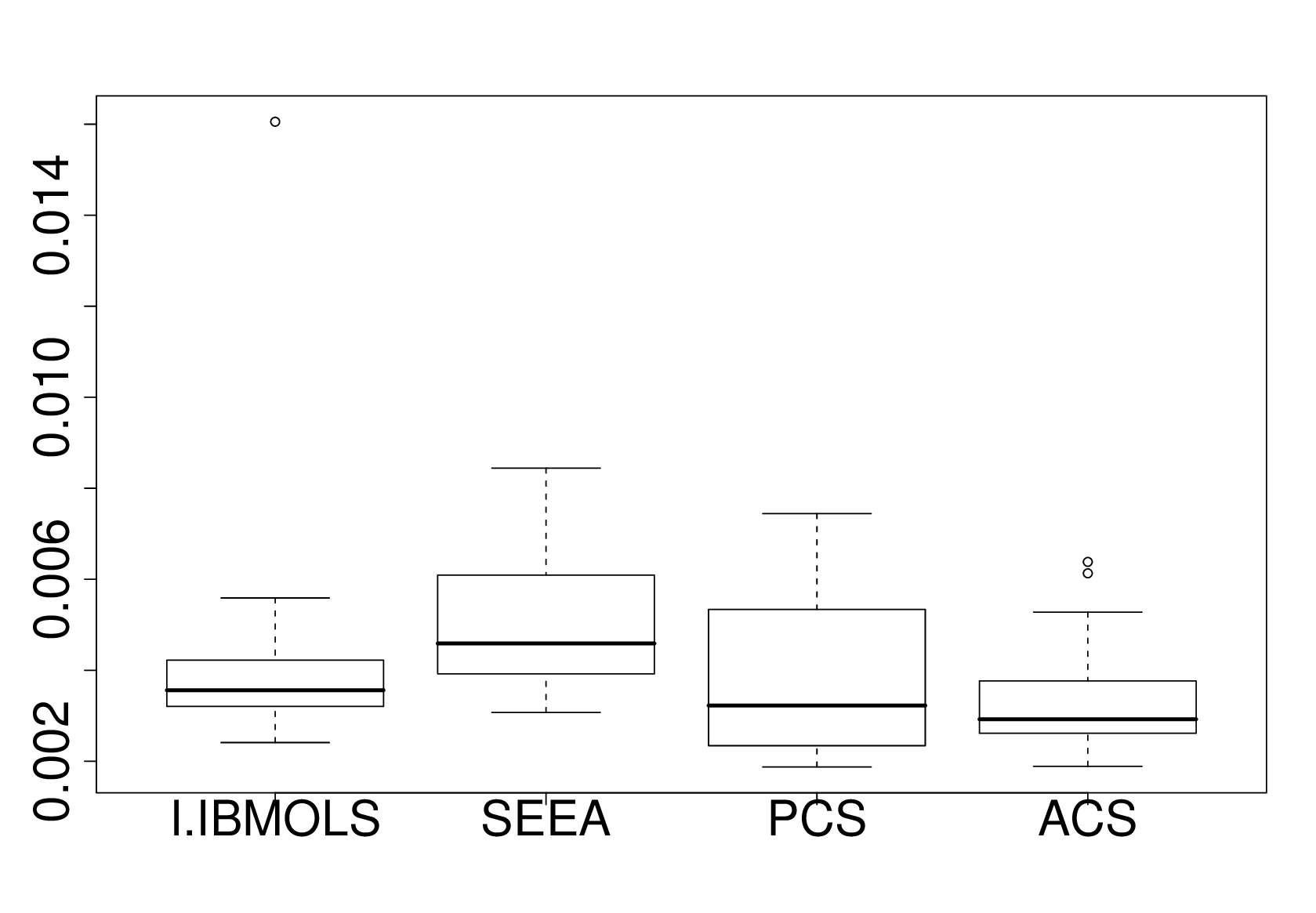}
\end{minipage}
\hfill
\begin{minipage}[c]{0.44\linewidth}
\centering
\includegraphics[angle=270,width=2in]{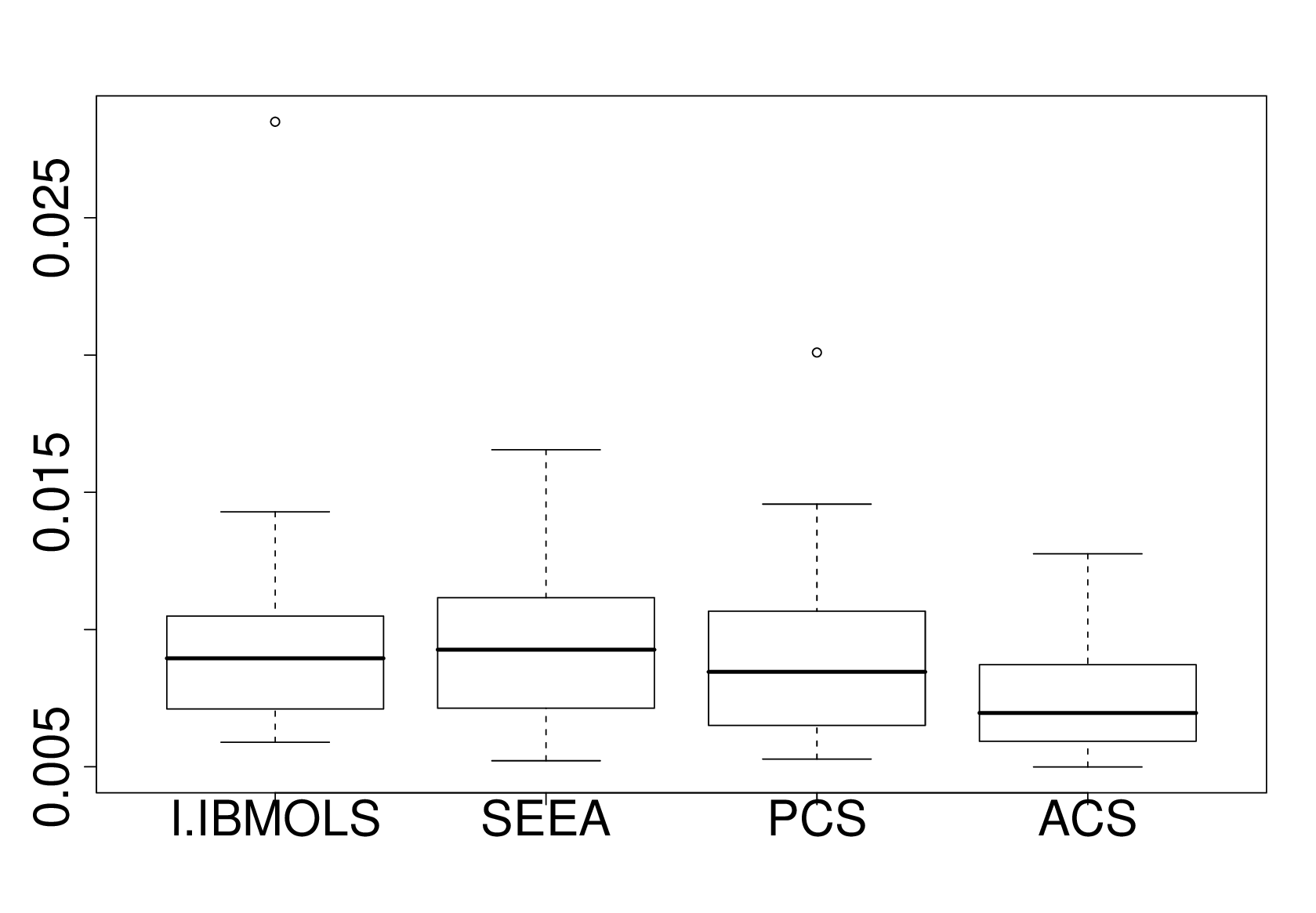}
\end{minipage}
\begin{minipage}[c]{0.1\linewidth}
\centering
{\itshape bier127}
\end{minipage}
\begin{minipage}[c]{0.44\linewidth}
\centering
\includegraphics[angle=270,width=2in]{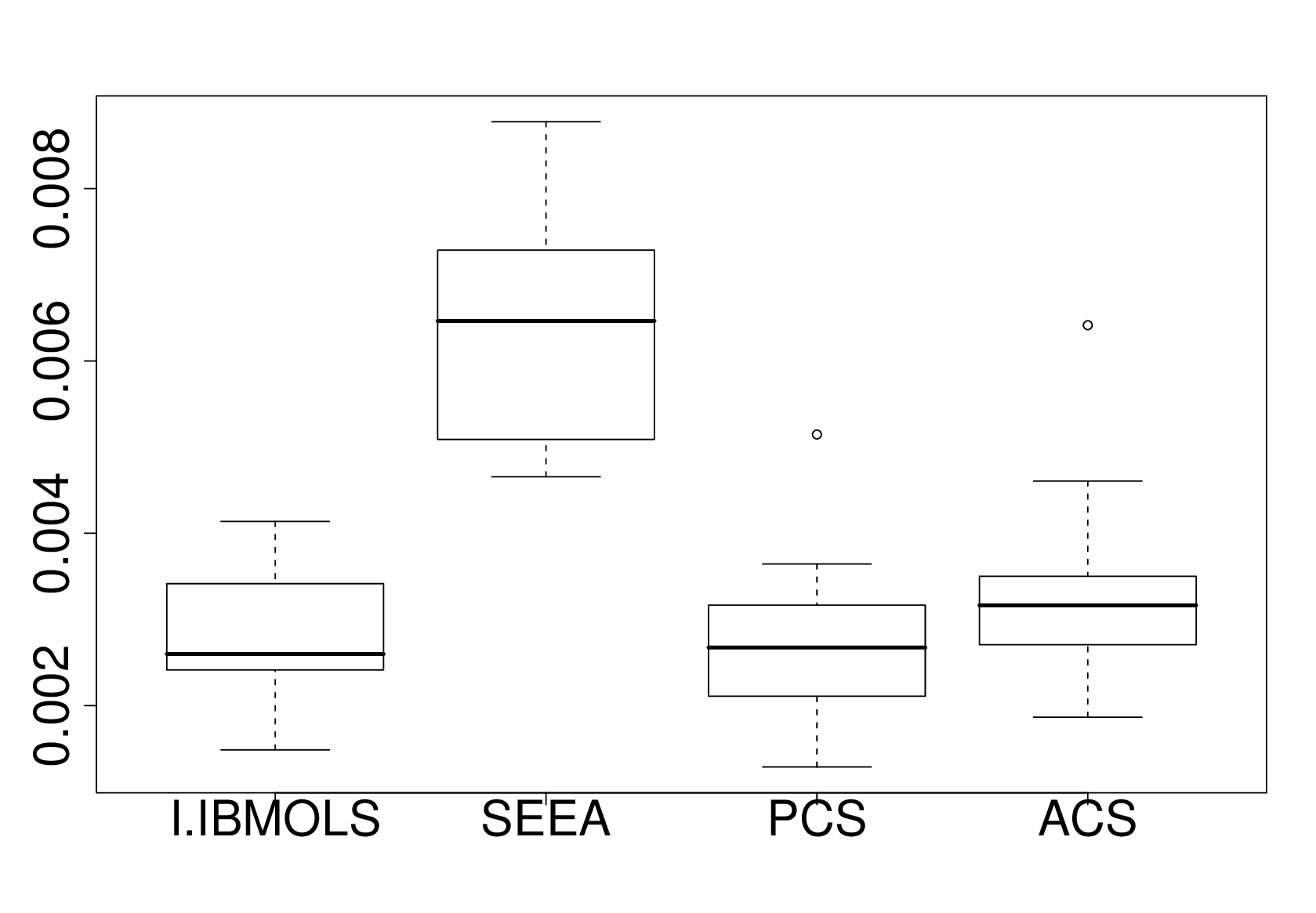}
\end{minipage}
\hfill
\begin{minipage}[c]{0.44\linewidth}
\centering
\includegraphics[angle=270,width=2in]{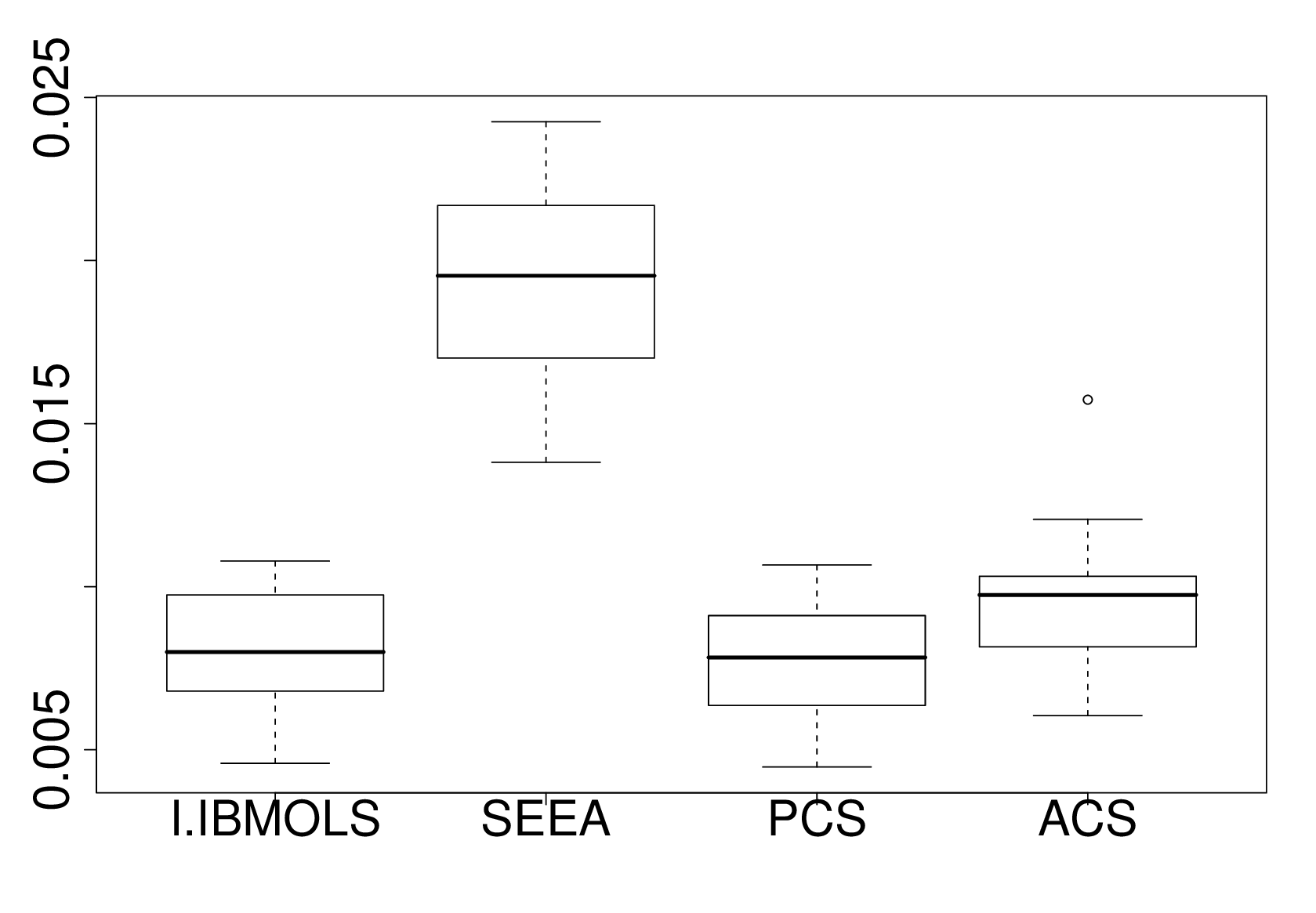}
\end{minipage}
\begin{minipage}[c]{0.1\linewidth}
\centering
{\itshape kroA150}
\end{minipage}
\begin{minipage}[c]{0.44\linewidth}
\centering
\includegraphics[angle=270,width=2in]{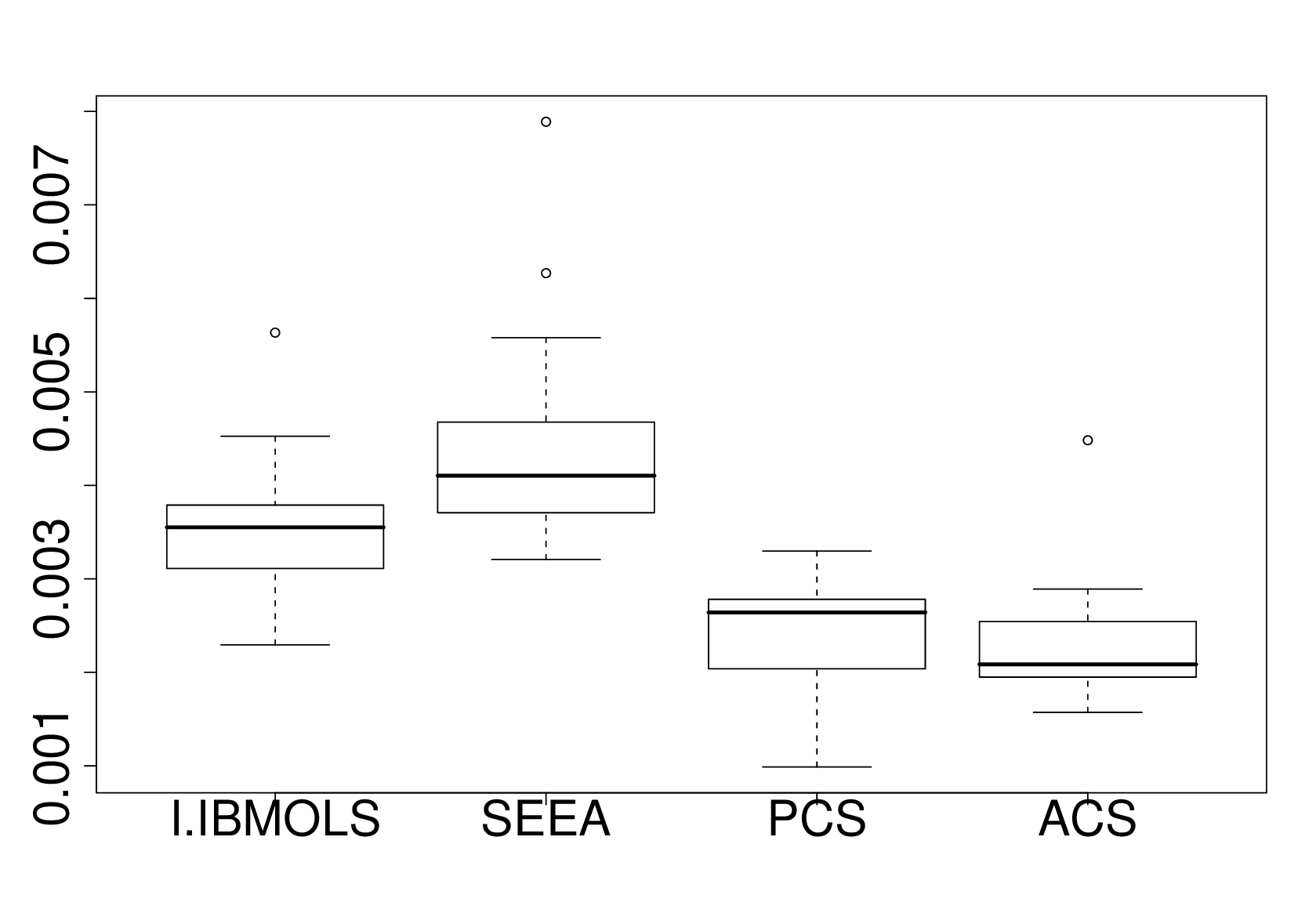}
\end{minipage}
\hfill
\begin{minipage}[c]{0.44\linewidth}
\centering
\includegraphics[angle=270,width=2in]{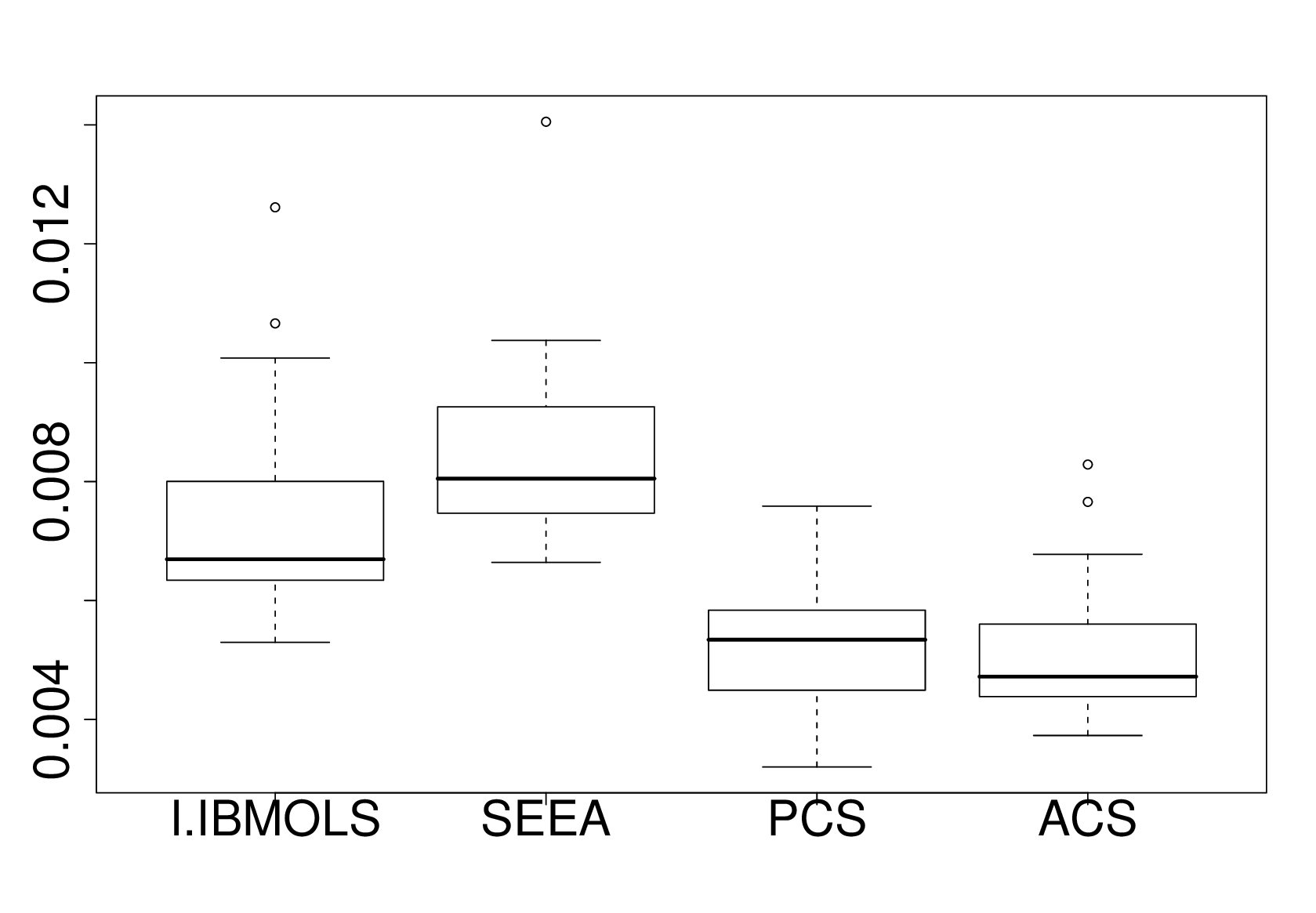}
\end{minipage}
\caption{Performance comparison for I-IBMOLS, SEEA, PCS and ACS according to the I$^-_H$ and the I$^1_{\epsilon+}$ metric (1).}
\label{fig:box_hybrid1}
\end{figure}

\begin{figure}[p]
\small 
\begin{minipage}[c]{0.1\linewidth}
Instance
\end{minipage}
\begin{minipage}[c]{0.44\linewidth}
\centering I$^-_H$
\end{minipage}
\hfill
\begin{minipage}[c]{0.44\linewidth}
\centering I$^1_{\epsilon+}$
\end{minipage}
\begin{minipage}[c]{0.1\linewidth}
\centering
{\itshape kroA200}
\end{minipage}
\begin{minipage}[c]{0.44\linewidth}
\centering
\includegraphics[angle=270,width=2in]{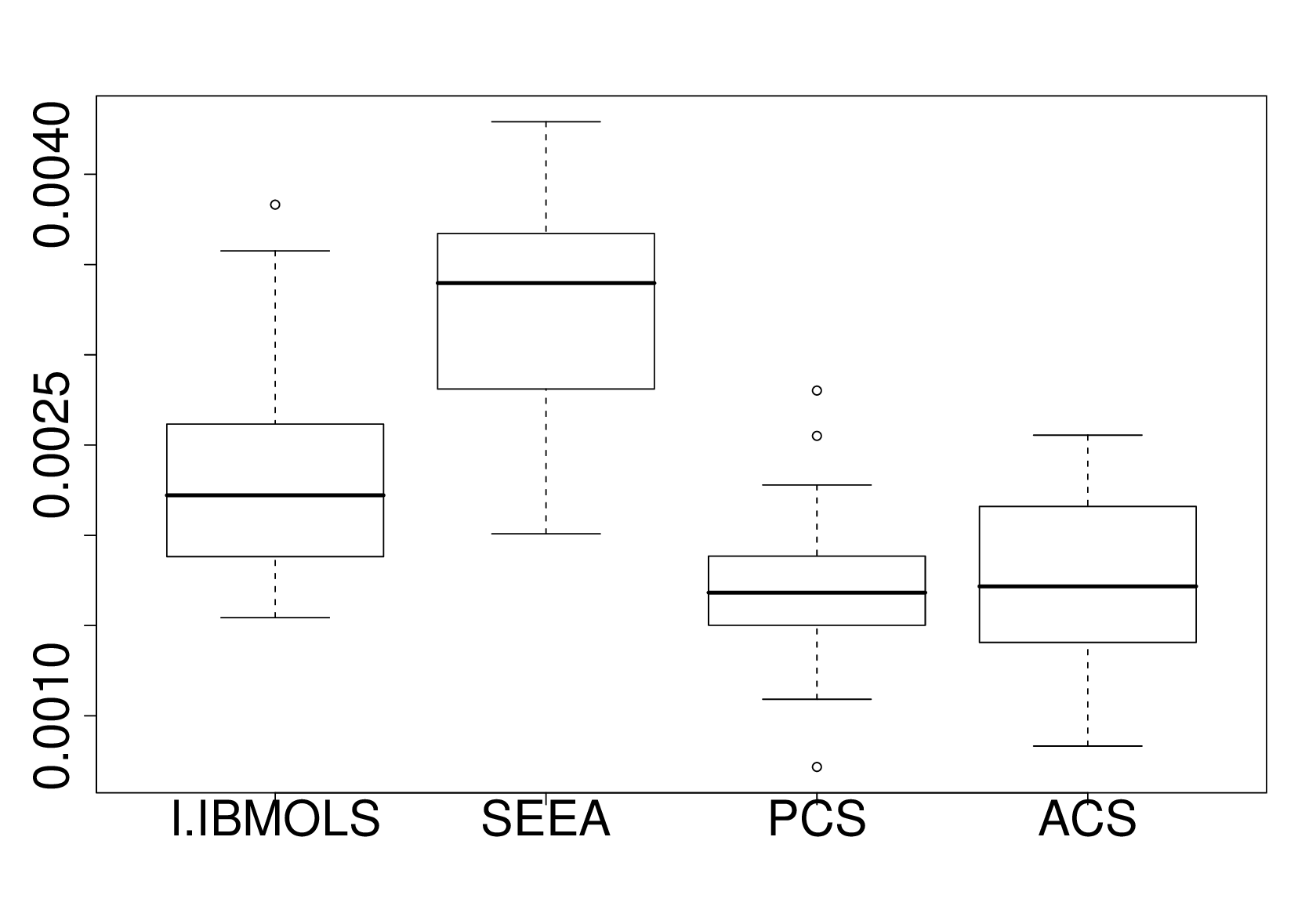}
\end{minipage}
\hfill
\begin{minipage}[c]{0.44\linewidth}
\centering
\includegraphics[angle=270,width=2in]{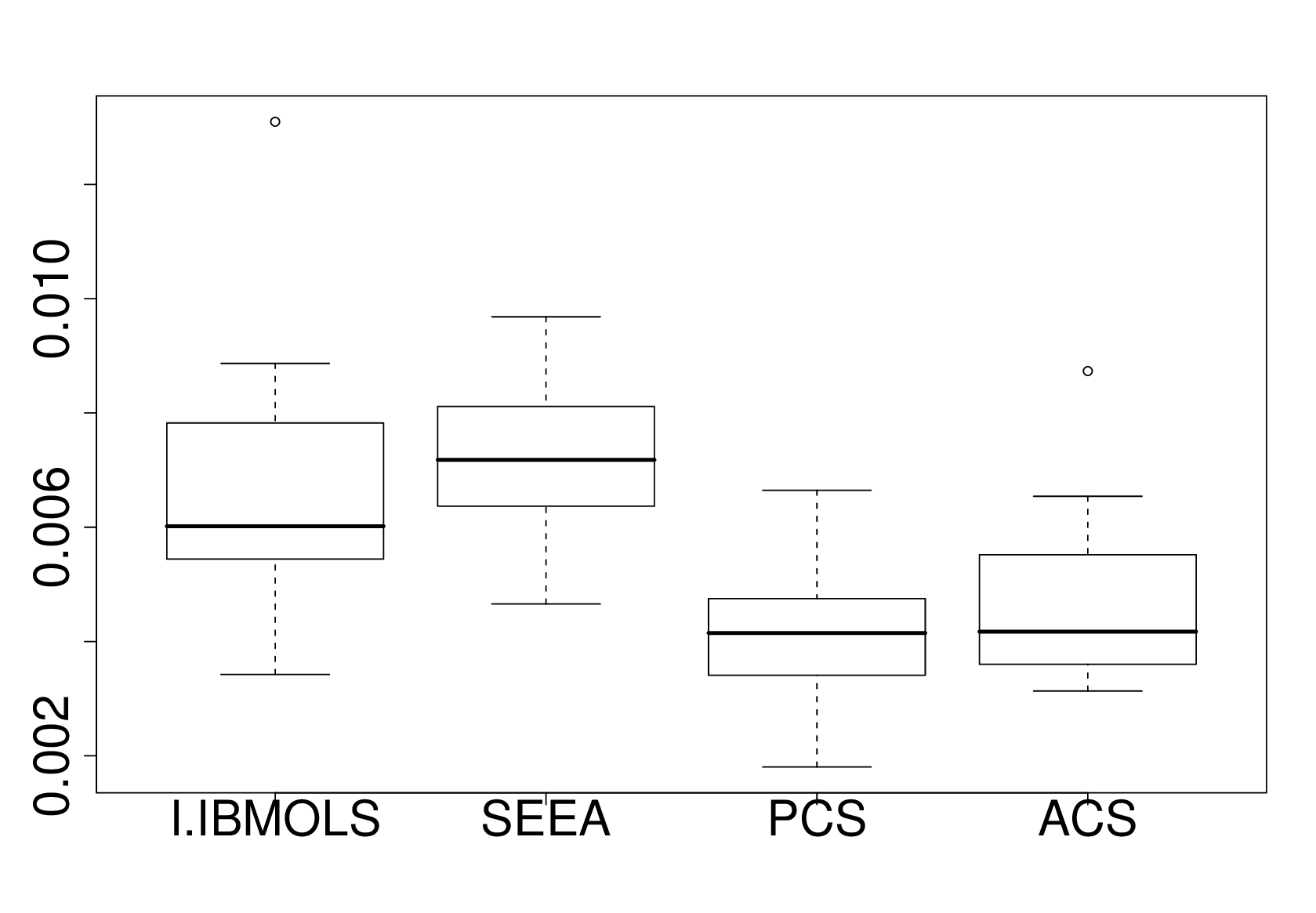}
\end{minipage}
\begin{minipage}[c]{0.1\linewidth}
\centering
{\itshape pr264}
\end{minipage}
\begin{minipage}[c]{0.44\linewidth}
\centering
\includegraphics[angle=270,width=2in]{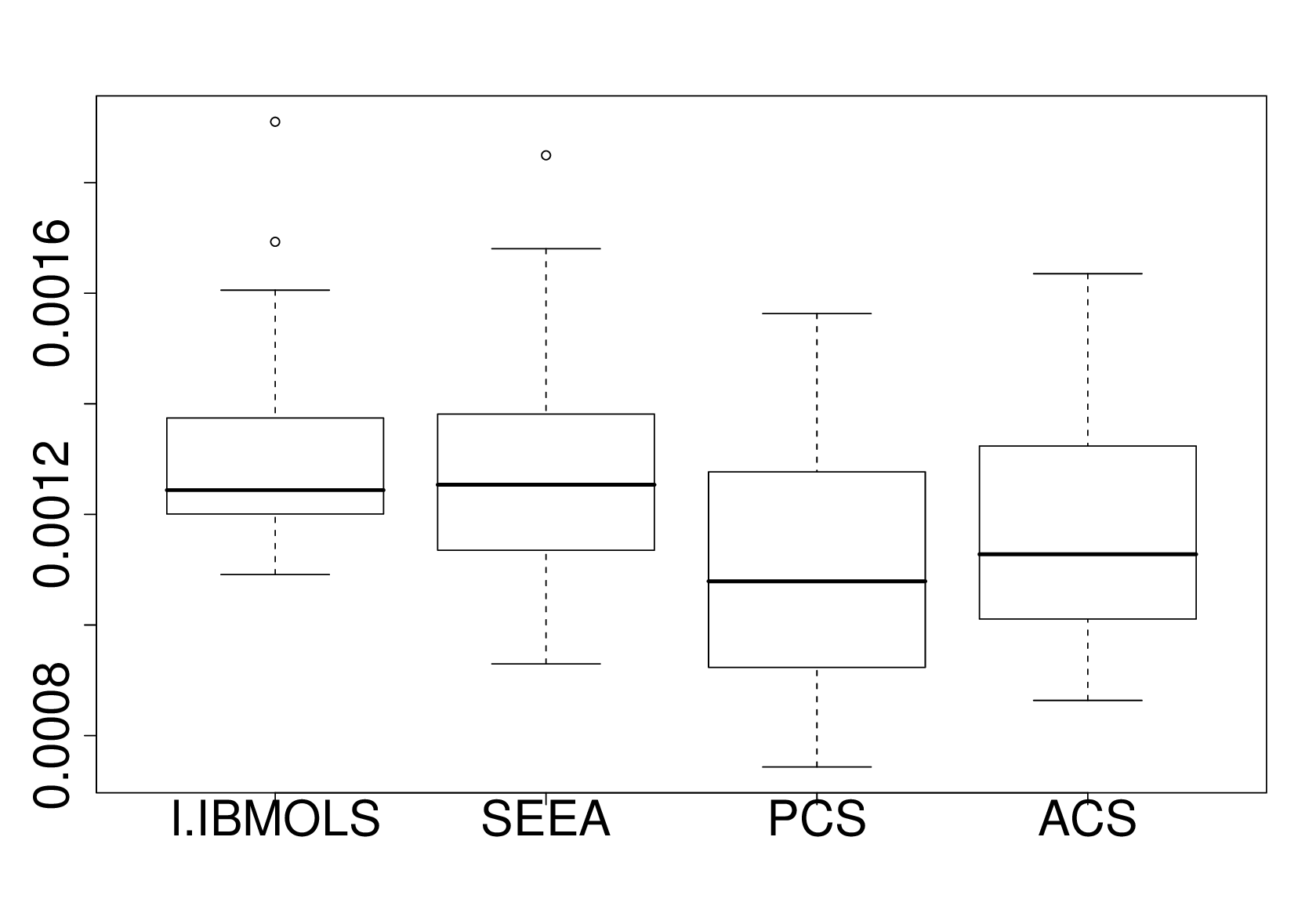}
\end{minipage}
\hfill
\begin{minipage}[c]{0.44\linewidth}
\centering
\includegraphics[angle=270,width=2in]{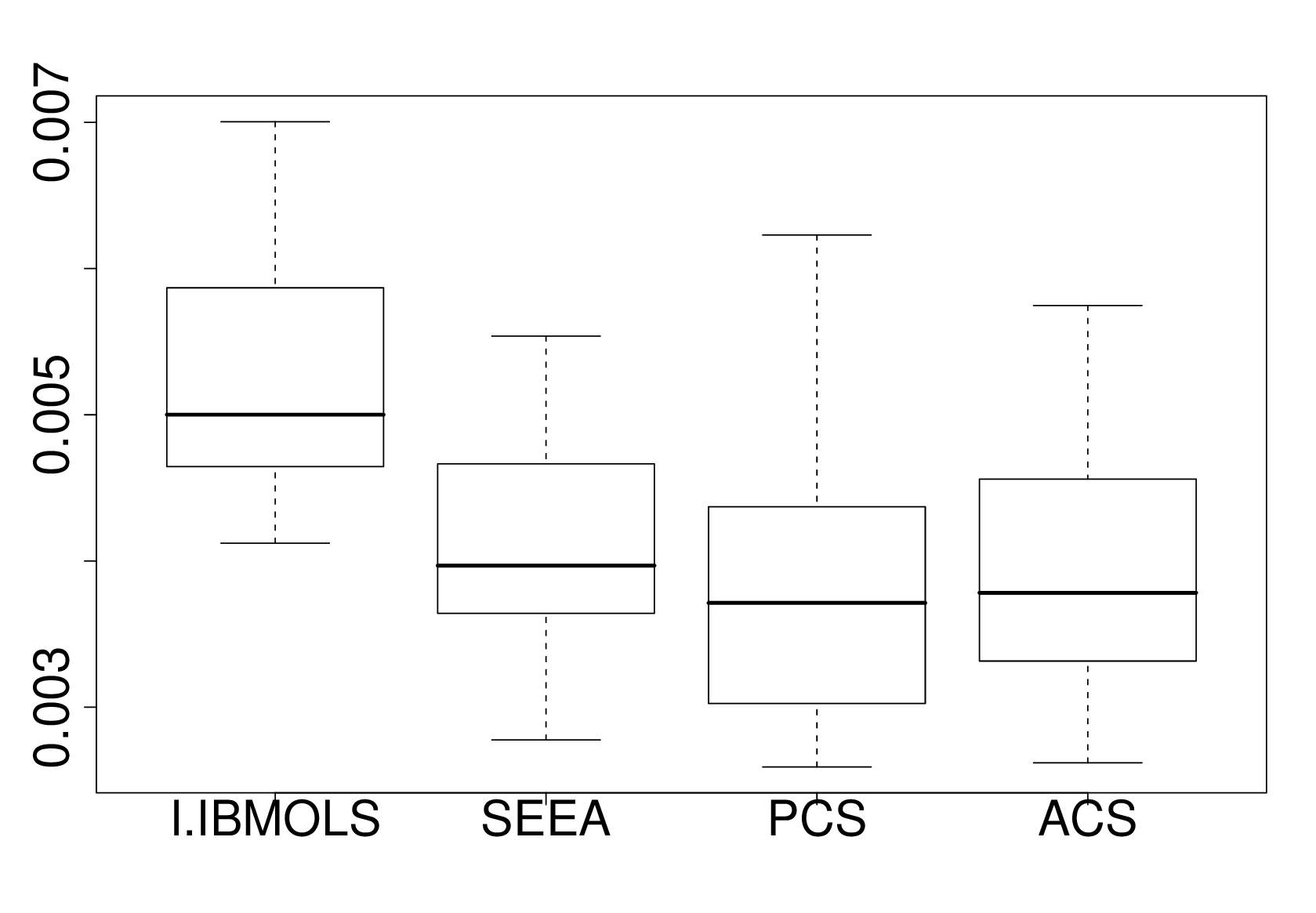}
\end{minipage}
\begin{minipage}[c]{0.1\linewidth}
\centering
{\itshape pr299}
\end{minipage}
\begin{minipage}[c]{0.44\linewidth}
\centering
\includegraphics[angle=270,width=2in]{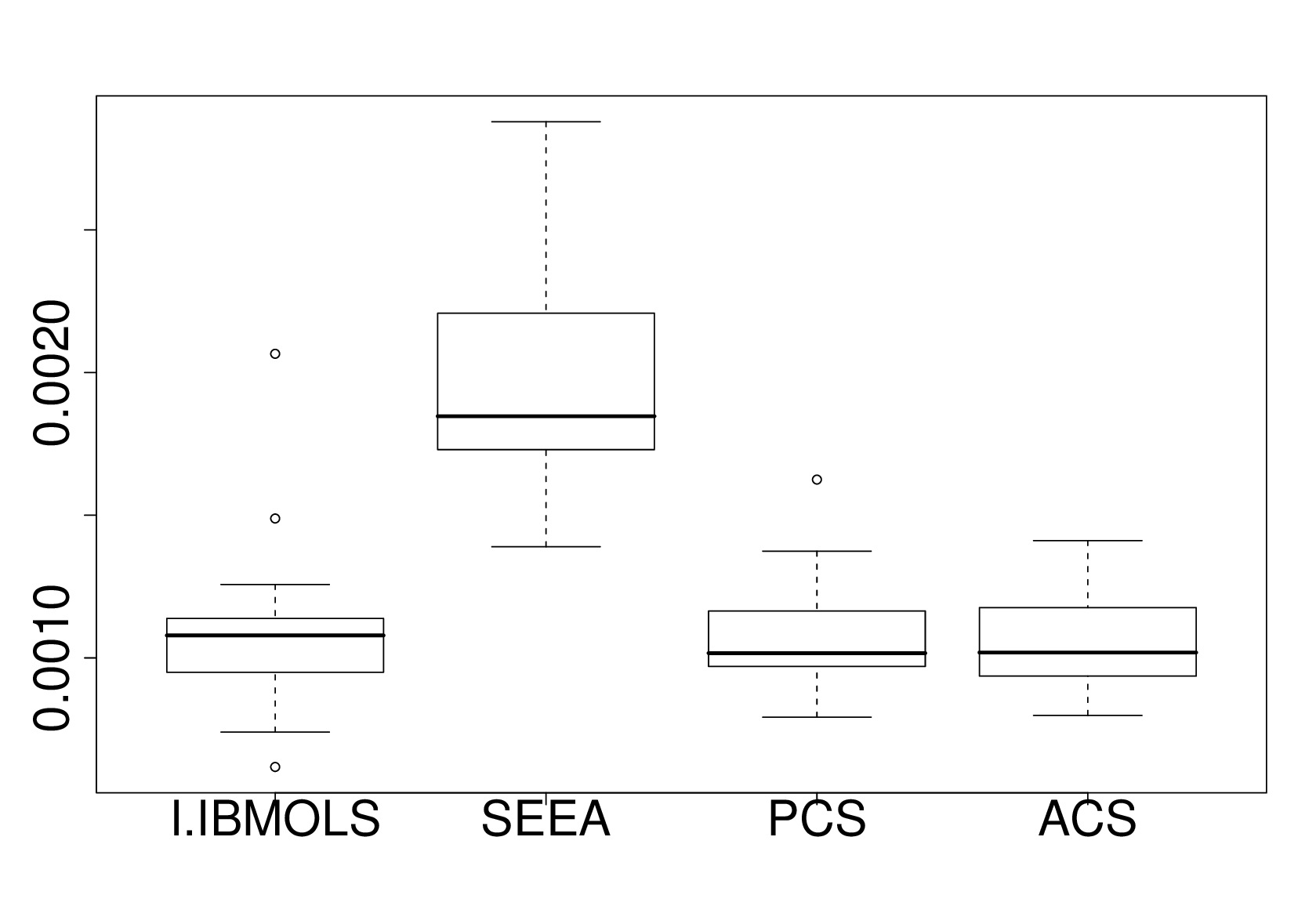}
\end{minipage}
\hfill
\begin{minipage}[c]{0.44\linewidth}
\centering
\includegraphics[angle=270,width=2in]{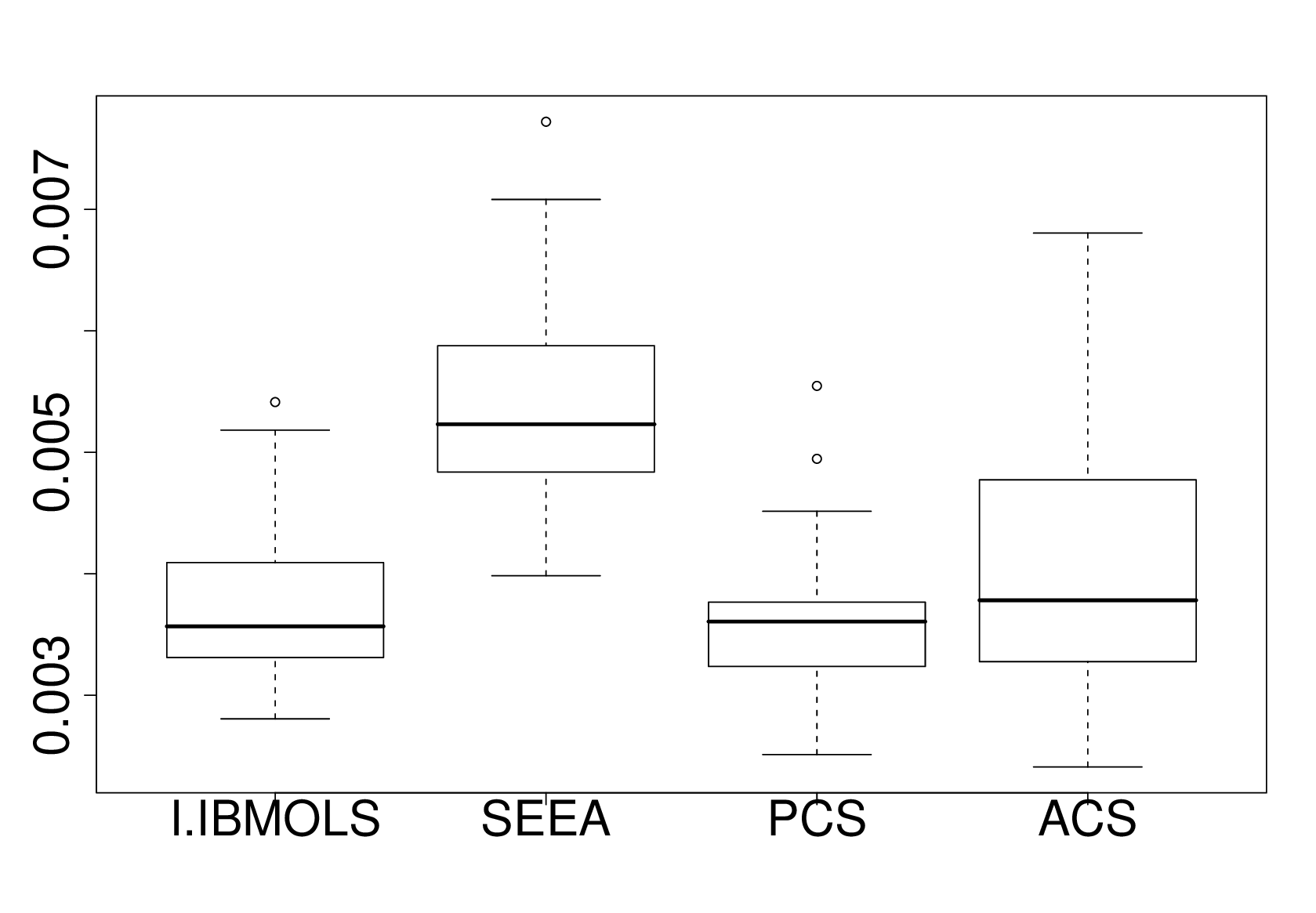}
\end{minipage}
\begin{minipage}[c]{0.1\linewidth}
\centering
{\itshape pr439}
\end{minipage}
\begin{minipage}[c]{0.44\linewidth}
\centering
\includegraphics[angle=270,width=2in]{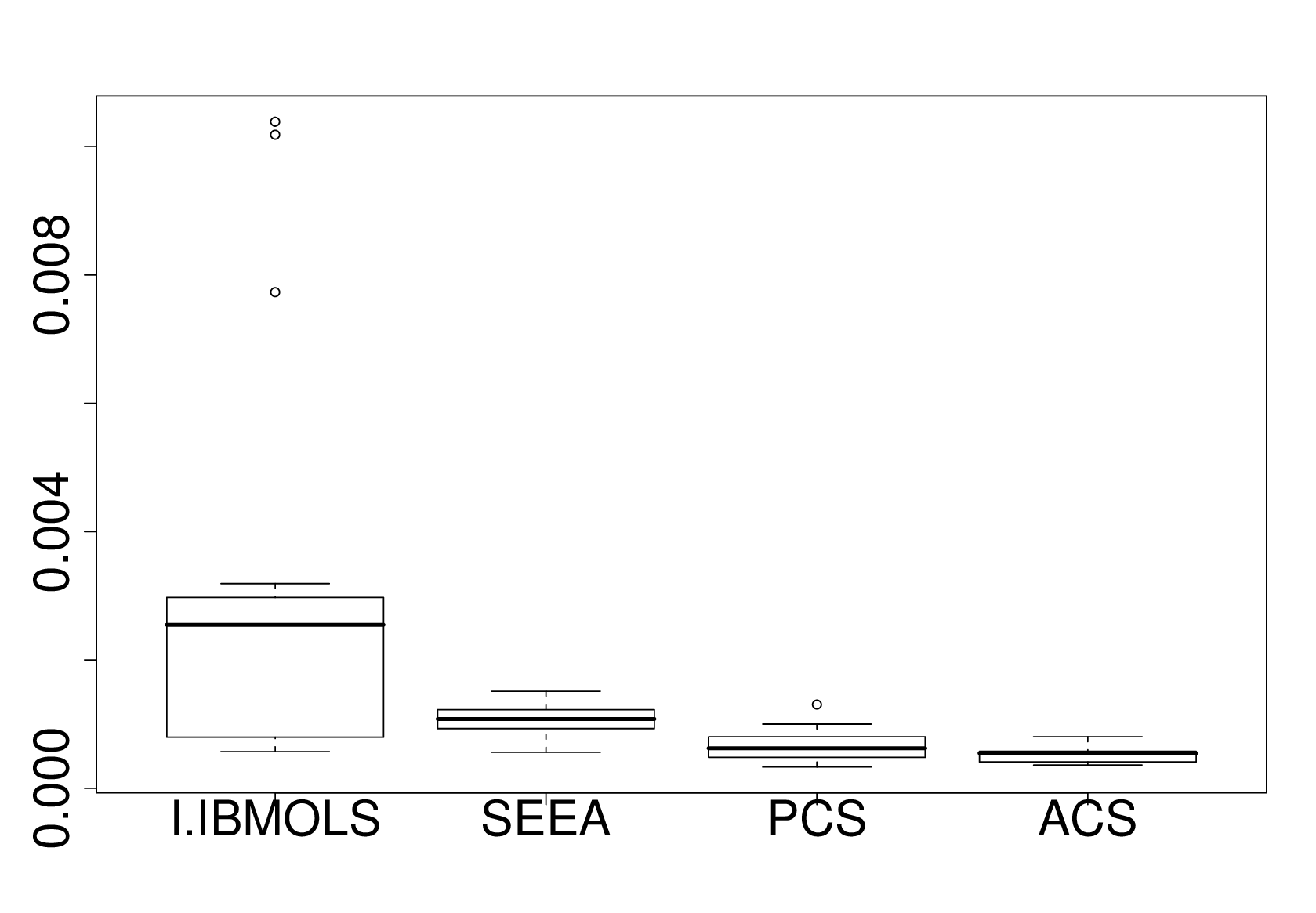}
\end{minipage}
\hfill
\begin{minipage}[c]{0.44\linewidth}
\centering
\includegraphics[angle=270,width=2in]{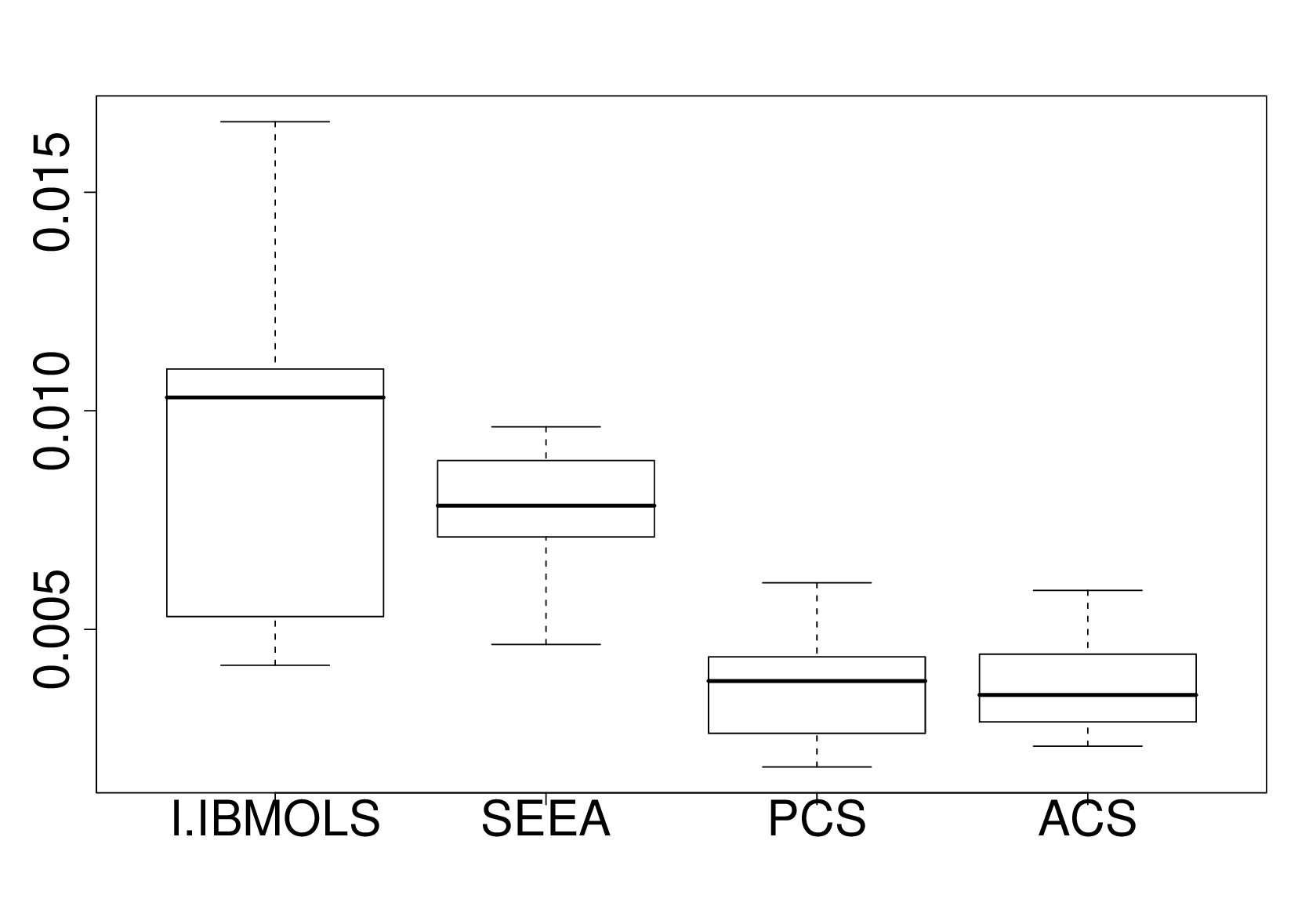}
\end{minipage}
\begin{minipage}[c]{0.1\linewidth}
\centering
{\itshape pr1002}
\end{minipage}
\begin{minipage}[c]{0.44\linewidth}
\centering
\includegraphics[angle=270,width=2in]{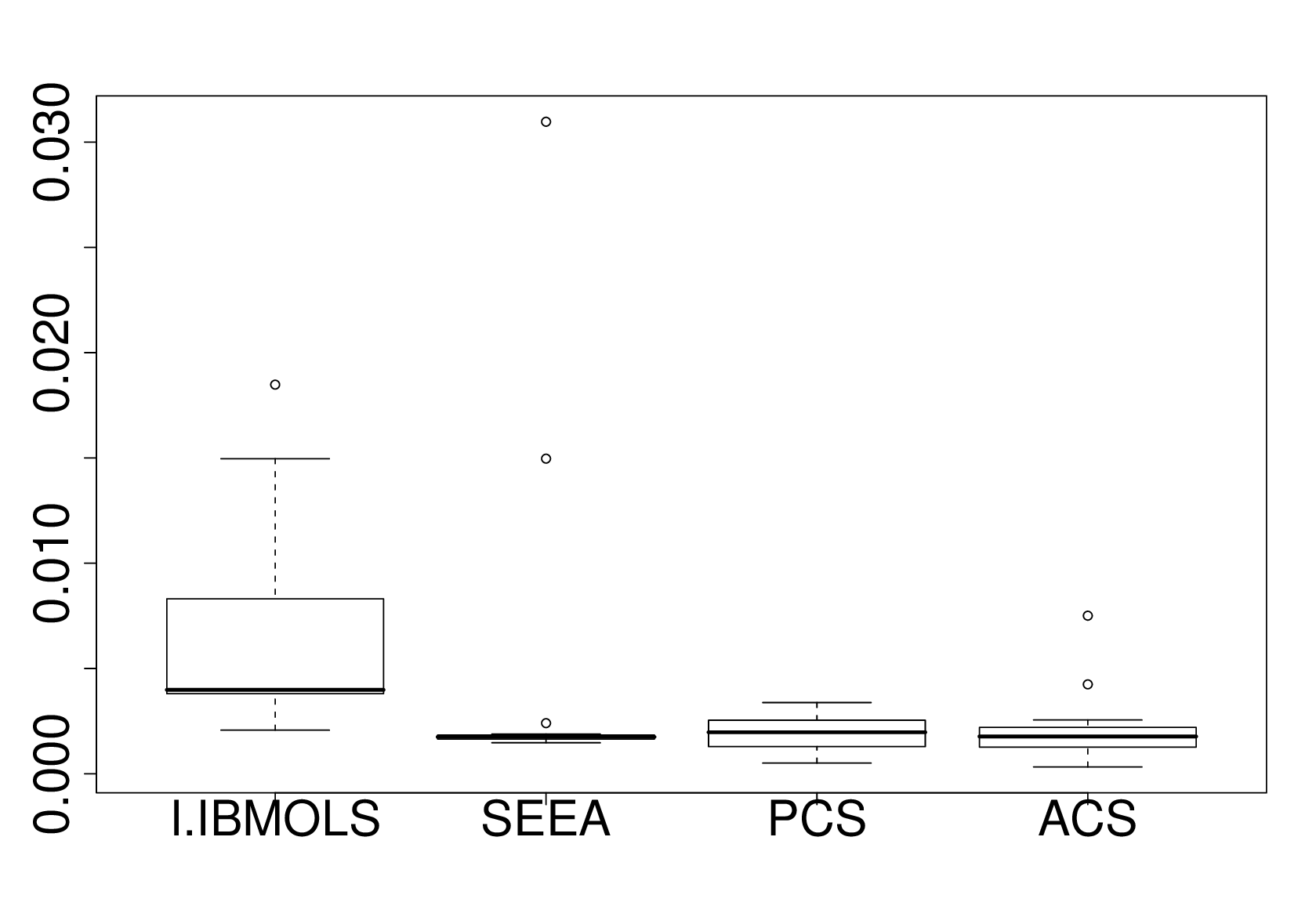}
\end{minipage}
\hfill
\begin{minipage}[c]{0.44\linewidth}
\centering
\includegraphics[angle=270,width=2in]{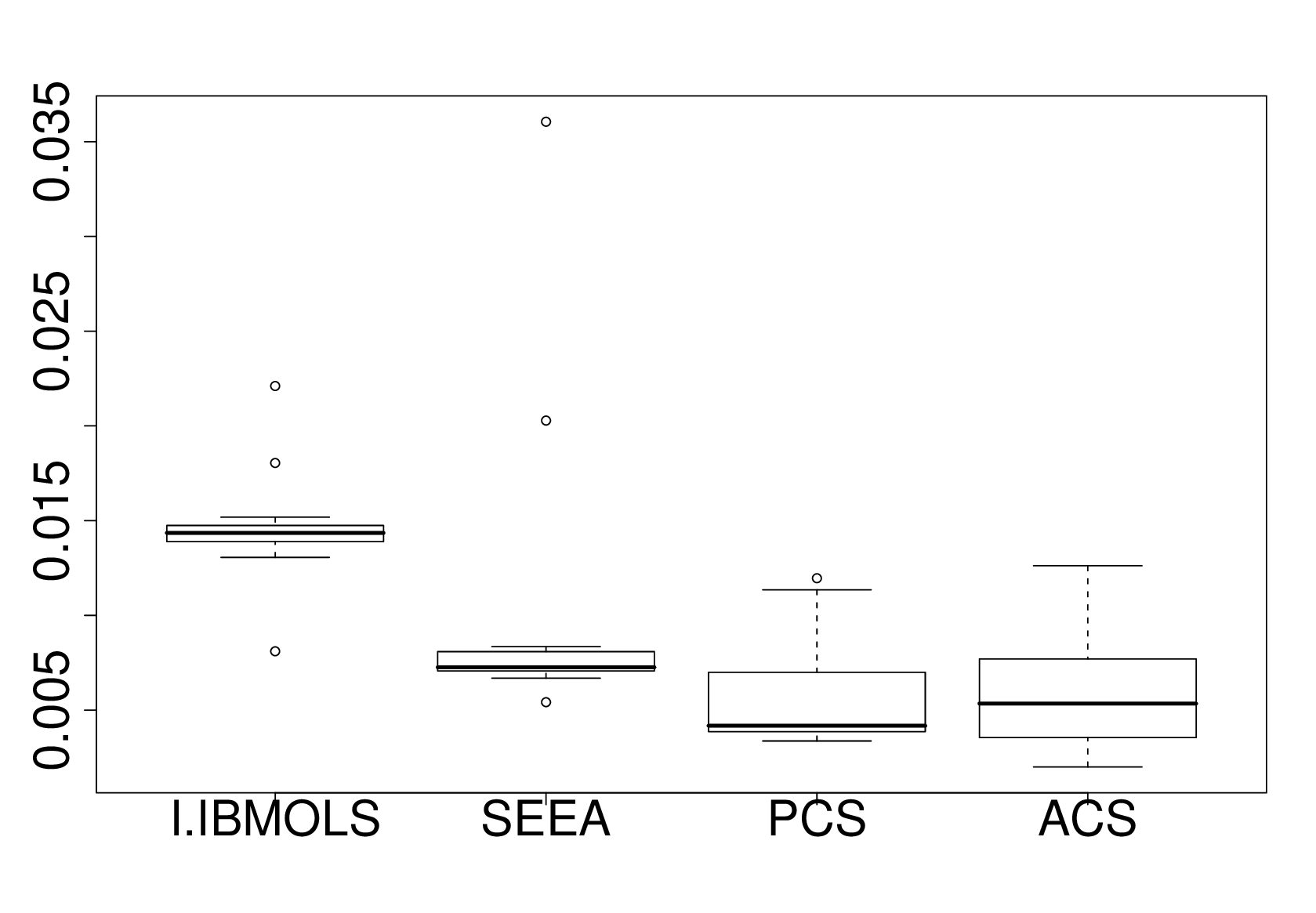}
\end{minipage}
\caption{Performance comparison for I-IBMOLS, SEEA, PCS and ACS according to the I$^-_H$ and the I$^1_{\epsilon+}$ metric (2).}
\label{fig:box_hybrid2}
\end{figure}

\begin{figure}[htbp]
\centering
\includegraphics[width=4.5in]{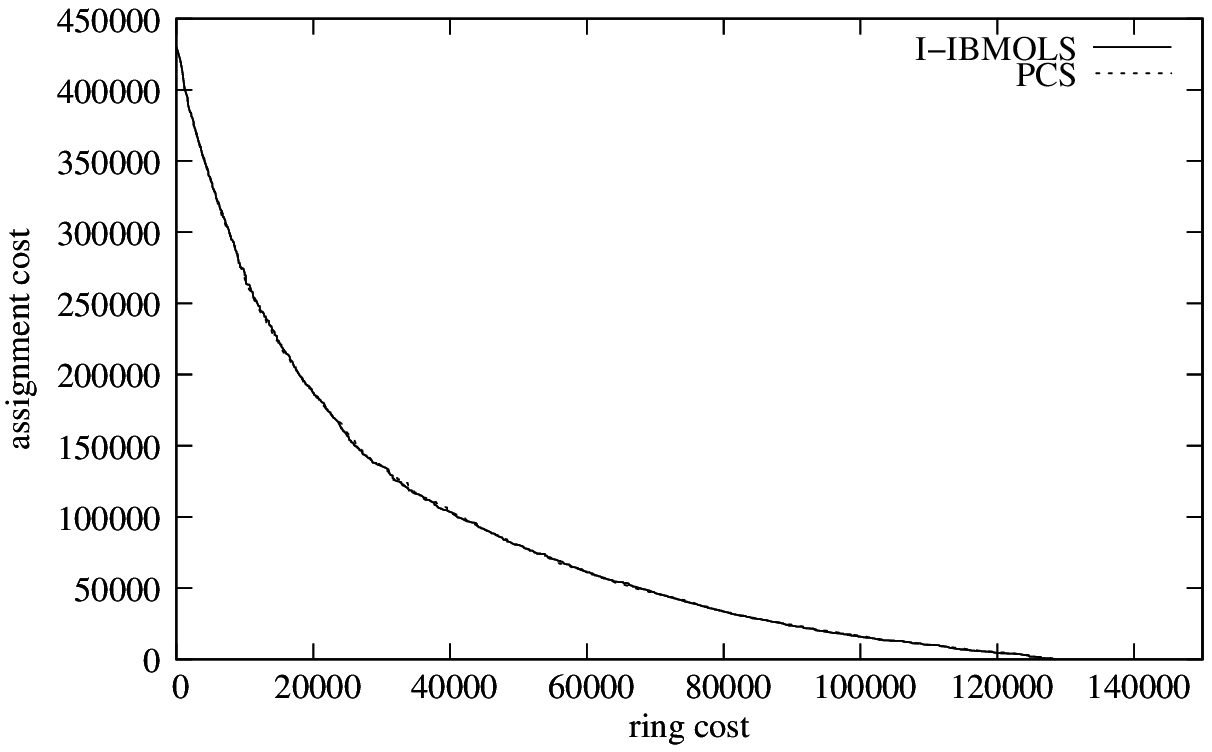}
\caption{$90\%$-attainment surface plot obtained by the approximation sets found by I-IBMOLS and PCS for the {\itshape bier127} test instance.}
\label{fig:eaf.bier127.hybrid}
\end{figure}

Table~\ref{tab:hybrid_hyp} and Table~\ref{tab:hybrid_eps} respectively summarize the results we obtained by comparing I-IBMOLS, SEEA, PCS and ACS 
according to the I$^-_H$ and to the I$^1_{\epsilon+}$ metrics.
For the $10$ benchmark test instances we experimented, we can see that both cooperative search methods are never statistically outperformed by SEEA.
Moreover, they perform significantly better according to at least one metric on all instances except {\itshape eil51}.
Similarly, PCS and ACS are at least as good as I-IBMOLS on every test instance.
While there is generally no significant difference for small-size instances, 
the efficient set approximations found by the cooperative approaches are statistically better than I-IBMOLS for problems with $150$ nodes and more 
(except for the {\itshape pr299} instance).
Thus, contrary to what has been pointed out in the previous section about SEEA, 
the cooperative methods do not seem to have trouble in finding solution having a low ring cost in comparison to I-IBMOLS, 
as illustrated in the empirical attainment function given in Figure~\ref{fig:eaf.bier127.hybrid}.
At last, in most cases, the differences between PCS and ACS are not statistically significant, 
except for the {\itshape bier127} instance where PCS performs better, 
and for the {\itshape kroA150} and the {\itshape pr439} instances where ACS performs better according to the I$^-_H$ metric.
Box-plots for I$^-_H$ and I$^1_{\epsilon+}$ are given in Figure~\ref{fig:box_hybrid1} and Figure~\ref{fig:box_hybrid2}.
In general, these figures graphically confirm the superiority of the cooperative search methods for test instances having a high number of nodes.
Thus, to summarize the results obtained by the hybrid search methods, the benefit of the cooperation scheme is being felt in most large-size test instances, 
but the addition of an auto-adaptive aspect does not seem to have a major influence on the results.
Nevertheless, as pointed out in Table~\ref{tab:nb_ibmols}, 
the difference between the average number of times that IBMOLS is launched during the search process of PCS and ACS is relatively thin.
This can be explained by the fact that 
($i$) IBMOLS puts more time to find non-dominated solutions by starting with a population of poorer quality, 
what is the case in PCS comparing to ACS (at least for the first launch),
and ($ii$) a part of the run-time allocated to the algorithm is used to compute a contribution value at every step of ACS, 
whereas PCS devotes all of its run-time in the search process.
These two aspects lead to the fact that the number of times that IBMOLS is launched is, in the end, more or less balanced between both cooperative methods.

\begin{table}[p]
\centering
\caption{
Performance comparison for I-IBMOLS, SEEA, PCS and ACS according to the I$^-_H$ metric by using a Mann-Whitney statistical test.
The ``p-value'' columns give the p-value of the statistical test.
The ``T'' columns give the outcome of the statistical test for a significance level of $5\%$: 
either the results of the search method located at a specific row are significantly better than those of the search method located at a specific column~($\succ$), 
either they are worse~($\prec$), or there is no significant difference between both~($\equiv$).
}
\hspace{500pt} \small \setlength{\tabcolsep}{1.8pt}
\begin{tabular}{cr|p{0.65in}||p{0.65in}|c||p{0.65in}|c||p{0.65in}|cc}
&  \multicolumn{2}{c||}{}  &  \multicolumn{2}{|c||}{I-IBMOLS}  &  \multicolumn{2}{|c||}{SEEA}  &  \multicolumn{2}{|c}{PCS}  &  \\
&  \multicolumn{2}{c||}{}  &  \centering p-value  &  T         &  \centering p-value  &  T     &  \centering p-value  &  T  &  \\
\hline
& {\itshape eil51}    &  SEEA  &  $> 5\%$                &  $\equiv$  &                         &            &                         &            &  \\
&                     &  PCS   &  $> 5\%$                &  $\equiv$  &  $> 5\%$                &  $\equiv$  &                         &            &  \\
&                     &  ACS   &  $> 5\%$                &  $\equiv$  &  $> 5\%$                &  $\equiv$  &  $> 5\%$                &  $\equiv$  &  \\
\hline
& {\itshape st70}     &  SEEA  &  $0.011$                &  $\prec$   &                         &            &                         &            &  \\
&                     &  PCS   &  $0.013$                &  $\succ$   &  $2.677 \cdot 10^{-4}$  &  $\succ$   &                         &            &  \\
&                     &  ACS   &  $> 5\%$                &  $\equiv$  &  $0.001$                &  $\succ$   &  $> 5\%$                &  $\equiv$  &  \\
\hline
& {\itshape kroA100}  &  SEEA  &  $0.004$                &  $\prec$   &                         &            &                         &            &  \\
&                     &  PCS   &  $> 5\%$                &  $\equiv$  &  $0.006$                &  $\succ$   &                         &            &  \\
&                     &  ACS   &  $0.015$                &  $\succ$   &  $1.053 \cdot 10^{-4}$  &  $\succ$   &  $> 5\%$                &  $\equiv$  &  \\
\hline
& {\itshape bier127}  &  SEEA  &  $3.151 \cdot 10^{-8}$  &  $\prec$   &                         &            &                         &            &  \\
&                     &  PCS   &  $> 5\%$                &  $\equiv$  &  $7.700 \cdot 10^{-8}$  &  $\succ$   &                         &            &  \\
&                     &  ACS   &  $> 5\%$                &  $\equiv$  &  $1.377 \cdot 10^{-7}$  &  $\succ$   &  $0.012$                &  $\prec$   &  \\
\hline
& {\itshape kroA150}  &  SEEA  &  $0.003$                &  $\prec$   &                         &            &                         &            &  \\
&                     &  PCS   &  $1.639 \cdot 10^{-6}$  &  $\succ$   &  $3.663 \cdot 10^{-8}$  &  $\succ$   &                         &            &  \\
&                     &  ACS   &  $1.259 \cdot 10^{-6}$  &  $\succ$   &  $2.114\cdot 10^{-7}$   &  $\succ$   &  $0.037$                &  $\succ$   &  \\
\hline
& {\itshape kroA200}  &  SEEA  &  $5.483 \cdot 10^{-5}$  &  $\prec$   &                         &            &                         &            &  \\
&                     &  PCS   &  $4.387 \cdot 10^{-4}$  &  $\succ$   &  $1.192 \cdot 10^{-7}$  &  $\succ$   &                         &            &  \\
&                     &  ACS   &  $9.330 \cdot 10^{-4}$  &  $\succ$   &  $8.914 \cdot 10^{-8}$  &  $\succ$   &  $> 5\%$                &  $\equiv$  &  \\
\hline
& {\itshape pr264}    &  SEEA  &  $> 5\%$                &  $\equiv$  &                         &            &                         &            &  \\
&                     &  PCS   &  $0.003$                &  $\succ$   &  $0.015$                &  $\succ$   &                         &            &  \\
&                     &  ACS   &  $0.012$                &  $\succ$   &  $0.058$                &  $\succ$   &  $> 5\%$                &  $\equiv$  &  \\
\hline
& {\itshape pr299}    &  SEEA  &  $3.225 \cdot 10^{-7}$  &  $\prec$   &                         &            &                         &            &  \\
&                     &  PCS   &  $> 5\%$                &  $\equiv$  &  $4.941 \cdot 10^{-8}$  &  $\succ$   &                         &            &  \\
&                     &  ACS   &  $> 5\%$                &  $\equiv$  &  $4.941 \cdot 10^{-8}$  &  $\succ$   &  $> 5\%$                &  $\equiv$  &  \\
\hline
& {\itshape pr439}    &  SEEA  &  $0.009$                &  $\succ$   &                         &            &                         &            &  \\
&                     &  PCS   &  $8.503 \cdot 10^{-6}$  &  $\succ$   &  $5.877 \cdot 10^{-6}$  &  $\succ$   &                         &            &  \\
&                     &  ACS   &  $1.377 \cdot 10^{-7}$  &  $\succ$   &  $1.031 \cdot 10^{-7}$  &  $\succ$   &  $0.033$                &  $\succ$   &  \\
\hline
& {\itshape pr1002}   &  SEEA  &  $7.523 \cdot 10^{-6}$  &  $\succ$   &                         &            &                         &            &  \\
&                     &  PCS   &  $1.192 \cdot 10^{-7}$  &  $\succ$   &  $> 5\%$                &  $\equiv$  &                         &            &  \\
&                     &  ACS   &  $1.868 \cdot 10^{-6}$  &  $\succ$   &  $> 5\%$                &  $\equiv$  &  $> 5\%$                &  $\equiv$  &  \\
\end{tabular}
\label{tab:hybrid_hyp}
\end{table}

\begin{table}[htbp]
\centering
\caption{
Performance comparison for I-IBMOLS, SEEA, PCS and ACS according to the I$^1_{\epsilon+}$ metric by using a Mann-Whitney statistical test.
The ``p-value'' columns give the p-value of the statistical test.
The ``T'' columns give the outcome of the statistical test for a significance level of $5\%$: 
either the results of the search method located at a specific row are significantly better than those of the search method located at a specific column~($\succ$), 
either they are worse~($\prec$), or there is no significant difference between both~($\equiv$).
}
\hspace{500pt} \small \setlength{\tabcolsep}{1.8pt}
\begin{tabular}{cr|p{0.65in}||p{0.65in}|c||p{0.65in}|c||p{0.65in}|cc}
&  \multicolumn{2}{c||}{}  &  \multicolumn{2}{|c||}{I-IBMOLS}  &  \multicolumn{2}{|c||}{SEEA}  &  \multicolumn{2}{|c}{PCS}  &  \\
&  \multicolumn{2}{c||}{}  &  \centering p-value  &  T         &  \centering p-value  &  T     &  \centering p-value  &  T  &  \\
\hline
& {\itshape eil51}    &  SEEA  &  $> 5\%$                &  $\equiv$  &                         &            &                         &            &  \\
&                     &  PCS   &  $> 5\%$                &  $\equiv$  &  $> 5\%$                &  $\equiv$  &                         &            &  \\
&                     &  ACS   &  $> 5\%$                &  $\equiv$  &  $> 5\%$                &  $\equiv$  &  $> 5\%$                &  $\equiv$  &  \\
\hline
& {\itshape st70}     &  SEEA  &  $> 5\%$                &  $\equiv$  &                         &            &                         &            &  \\
&                     &  PCS   &  $> 5\%$                &  $\equiv$  &  $> 5\%$                &  $\equiv$  &                         &            &  \\
&                     &  ACS   &  $> 5\%$                &  $\equiv$  &  $> 5\%$                &  $\equiv$  &  $> 5\%$                &  $\equiv$  &  \\
\hline
& {\itshape kroA100}  &  SEEA  &  $> 5\%$                &  $\equiv$  &                         &            &                         &            &  \\
&                     &  PCS   &  $> 5\%$                &  $\equiv$  &  $> 5\%$                &  $\equiv$  &                         &            &  \\
&                     &  ACS   &  $0.007$                &  $\succ$   &  $0.006$                &  $\succ$   &  $> 5\%$                &  $\equiv$  &  \\
\hline
& {\itshape bier127}  &  SEEA  &  $3.146 \cdot 10^{-8}$  &  $\prec$   &                         &            &                         &            &  \\
&                     &  PCS   &  $> 5\%$                &  $\equiv$  &  $3.151 \cdot 10^{-8}$  &  $\succ$   &                         &            &  \\
&                     &  ACS   &  $> 5\%$                &  $\equiv$  &  $4.935 \cdot 10^{-8}$  &  $\succ$   &  $> 5\%$                &  $\equiv$  &  \\
\hline
& {\itshape kroA150}  &  SEEA  &  $0.004$                &  $\prec$   &                         &            &                         &            &  \\
&                     &  PCS   &  $1.378 \cdot 10^{-5}$  &  $\succ$   &  $1.192 \cdot 10^{-7}$  &  $\succ$   &                         &            &  \\
&                     &  ACS   &  $1.747 \cdot 10^{-5}$  &  $\succ$   &  $3.708 \cdot 10^{-7}$  &  $\succ$   &  $> 5\%$                &  $\equiv$  &  \\
\hline
& {\itshape kroA200}  &  SEEA  &  $0.047$                &  $\prec$   &                         &            &                         &            &  \\
&                     &  PCS   &  $6.124 \cdot 10^{-5}$  &  $\succ$   &  $6.428 \cdot 10^{-7}$  &  $\succ$   &                         &            &  \\
&                     &  ACS   &  $7.067 \cdot 10^{-4}$  &  $\succ$   &  $8.503 \cdot 10^{-6}$  &  $\succ$   &  $> 5\%$                &  $\equiv$  &  \\
\hline
& {\itshape pr264}    &  SEEA  &  $4.386 \cdot 10^{-5}$  &  $\succ$   &                         &            &                         &            &  \\
&                     &  PCS   &  $2.208 \cdot 10^{-5}$  &  $\succ$   &  $> 5\%$                &  $\equiv$  &                         &            &  \\
&                     &  ACS   &  $1.963 \cdot 10^{-5}$  &  $\succ$   &  $> 5\%$                &  $\equiv$  &  $> 5\%$                &  $\equiv$  &  \\
\hline
& {\itshape pr299}    &  SEEA  &  $3.556 \cdot 10^{-6}$  &  $\prec$   &                         &            &                         &            &  \\
&                     &  PCS   &  $> 5\%$                &  $\equiv$  &  $1.259 \cdot 10^{-6}$  &  $\succ$   &                         &            &  \\
&                     &  ACS   &  $> 5\%$                &  $\equiv$  &  $1.974 \cdot 10^{-4}$  &  $\succ$   &  $> 5\%$                &  $\equiv$  &  \\
\hline
& {\itshape pr439}    &  SEEA  &  $0.018$                &  $\succ$   &                         &            &                         &            &  \\
&                     &  PCS   &  $7.363 \cdot 10^{-7}$  &  $\succ$   &  $5.733 \cdot 10^{-8}$  &  $\succ$   &                         &            &  \\
&                     &  ACS   &  $4.889 \cdot 10^{-7}$  &  $\succ$   &  $5.733 \cdot 10^{-8}$  &  $\succ$   &  $> 5\%$                &  $\equiv$  &  \\
\hline
& {\itshape pr1002}   &  SEEA  &  $9.604 \cdot 10^{-6}$  &  $\succ$   &                         &            &                         &            &  \\
&                     &  PCS   &  $4.256 \cdot 10^{-8}$  &  $\succ$   &  $3.268 \cdot 10^{-4}$  &  $\succ$   &                         &            &  \\
&                     &  ACS   &  $6.646 \cdot 10^{-8}$  &  $\succ$   &  $0.006$                &  $\succ$   &  $> 5\%$                &  $\equiv$  &  \\
\end{tabular}
\label{tab:hybrid_eps}
\end{table}

\begin{table}[htbp]
\centering
\caption{Average number of times that IBMOLS is launched during the search process of PCS and ACS.}
\hspace{500pt} \small \setlength{\tabcolsep}{3pt}
\begin{tabular}{cr||r|rc}
&  Instance           &  \centering PCS  &  \centering ACS  &  \\
\hline
& {\itshape eil51}    &    9.35  &    8.70  &  \\
& {\itshape st70}     &   19.80  &   20.00  &  \\
& {\itshape kroA100}  &    4.60  &    6.25  &  \\
& {\itshape bier127}  &    4.70  &    3.50  &  \\
& {\itshape kroA150}  &   16.35  &   17.45  &  \\
& {\itshape kroA200}  &   23.15  &   22.15  &  \\
& {\itshape pr264}    &   43.00  &   43.15  &  \\
& {\itshape pr299}    &   38.90  &   26.90  &  \\
& {\itshape pr439}    &   18.05  &   12.80  &  \\
& {\itshape pr1002}   &    3.50  &    2.05  &  \\
\end{tabular}
\label{tab:nb_ibmols}
\end{table}

\paragraph{Comparison with Exact Single-objective Results.}
As a last step, we provide a comparison between the results found for the bi-objective version of the RSP proposed in this paper 
and the ones of the mono-objective RSP investigated in~\cite{LL+:04}, where both costs are summed up.
Note that it was not possible to compare our results to the ones of the other formulation of a single-objective RSP investigated in~\cite{LL+:05,MMR:03,RBL:04}, 
where the assignment cost is subject to a constraint, due to the way the bound has been fixed.
In~\cite{LL+:04}, the authors propose an exact method to solve a single-objective RSP aiming at minimizing the sum of the ring cost and of the assignment cost.
But, in order to provide optimal solutions visiting approximately $25$, $50$, $75$ and $100\%$ of the total number of nodes,
they set the ring cost $c_{ij}$ and the assignment cost $d_{ij}$ between two nodes $v_i$ and $v_j$ in the following way:
$c_{ij} = \lceil\alpha l_{ij}\rceil$ and $d_{ij} = \lceil(10-\alpha) l_{ij}\rceil$ with $\alpha \in \{3,5,7,9\}$, 
where $l_{ij}$ denotes the distance between $v_i$ and $v_j$ given in the TSPLIB files.
In order to make our results comparable to their, 
we transform the objective vectors of the best found Pareto front $Z_N^*$ we obtained for every instance into a scalar value as detailed above,
and we compare the minimal one with the optimum%
\footnote{In fact, the authors imposed a time limit for their experiments.
They report the best solution found so far for the instances exceeding this time limit, 
what is the case for the {\itshape kroA200} instance with an $\alpha = 3$, $5$ and $9$ in Table~\ref{tab:compare}.}
found in~\cite{LL+:04}.
Additionally, we compare the best solution visiting every nodes we have found with the optimal TSP solution available in the TSPLIB website%
\footnote{\url{http://www.iwr.uni-heidelberg.de/groups/comopt/software/TSPLIB95/}.}.
Table~\ref{tab:compare} gives the error ratio between the best known value and the best one we have found for every identified (single) objective 
and every benchmark test instance.
In comparison to~\cite{LL+:04}, this ratio is always under $1.5\%$, and is mostly below $1\%$ for every instance.
The optimum is even found for the {\itshape kroA100} and the {\itshape kroA150} instances with an $\alpha = 9$,
and a better solution is found for the {\itshape kroA200} instances with an $\alpha = 3$ and an $\alpha = 9$.
As regards to the optimal TSP solutions (corresponding to solutions visiting every node in the case of the RSP), 
our results are similar for instances with $200$ nodes or less.
For larger instances, they are slightly less good, especially for the {\itshape pr1002} instance, where the error ratio is close to $10\%$.
But, proportionally to the number of nodes, the running time available to solve this test instance is lower than for the other instances.
To sum up, compared with single-objective optimal or near optimal results, 
the search methods we proposed in this paper to solve the B-RSP are quite promising 
with regard to the relatively small computational time available and to the size of the problem instances to be solved.
Nevertheless, it is true that the solutions used to compare the results are the best ones we obtained during the set of all experiments we performed, 
and may not be as good for a single simulation run of a given search method.

\begin{table}[htbp]
\centering	
\caption{Error ratio between the ring cost value of the best known TSP solution 
and the best single objective value found in~\cite{LL+:04} in comparison to the best value found during our experiments.}
\hspace{500pt} \small \setlength{\tabcolsep}{3pt}
\begin{tabular}{cr||r||r|r|r|rc}
&  Instance           &  \multicolumn{1}{c||}{Optimal TSP}  &  \multicolumn{4}{c}{Optimal single-objective}    &  \\
&                     &  \multicolumn{1}{c||}{solution}     &  \multicolumn{4}{c}{RSP solution~\cite{LL+:04}}  &  \\
&  &  &  \multicolumn{1}{c|}{$\alpha=3$}  &  \multicolumn{1}{c|}{$\alpha=5$}  &  \multicolumn{1}{c|}{$\alpha=7$}  &  \multicolumn{1}{c}{$\alpha=9$}  &  \\
\hline
& {\itshape eil51}    &  0.67\%  &     0.67\%  &  0.75\%  &  0.37\%  &     0.69\%  &  \\
& {\itshape st70}     &  0.31\%  &     0.31\%  &  0.49\%  &  0.56\%  &     0.42\%  &  \\
& {\itshape kroA100}  &  0.08\%  &     0.08\%  &  0.12\%  &  0.56\%  &     0.00\%  &  \\
& {\itshape bier127}  &  0.64\%  &     0.64\%  &  0.34\%  &  0.52\%  &     0.01\%  &  \\
& {\itshape kroA150}  &  1.40\%  &     1.40\%  &  1.11\%  &  0.36\%  &     0.00\%  &  \\
& {\itshape kroA200}  &  1.22\%  &  $-$4.82\%  &  1.05\%  &  0.81\%  &  $-$1.44\%  &  \\
& {\itshape pr264}    &  2.26\%  &          -  &       -  &       -  &          -  &  \\
& {\itshape pr299}    &  1.72\%  &          -  &       -  &       -  &          -  &  \\
& {\itshape pr439}    &  2.74\%  &          -  &       -  &       -  &          -  &  \\
& {\itshape pr1002}   &  9.51\%  &          -  &       -  &       -  &          -  &  \\
\end{tabular}
\label{tab:compare}
\end{table}
\section{Conclusion and Perspectives}

A new multi-objective routing problem, the bi-objective ring star problem, has been investigated in this paper for the first time.
It aims to locate a cycle through a subset of nodes of a graph while minimizing a ring cost, related to the length of the cycle, 
and an assignment cost, from non-visited nodes to visited ones.
Despite it is clearly bi-objective, this problem has always been tackled in a single-objective formulation, 
either where both costs are combined~\cite{LL+:04}, or where one cost is regarded as a constraint~\cite{LL+:05,RBL:04}.
As a first step, four population-based metaheuristics have been proposed to approximate the minimal set of efficient solutions for the problem under consideration.
A population-based local search method, IBMOLS, recently proposed in~\cite{BB:07}, has first been designed, in an iterative way, with a variable neighborhood.
Then, two state-of-the-art multi-objective evolutionary algorithms, namely IBEA~\cite{ZK:04} and NSGA-II~\cite{DA+:02}, 
as well as a third one, the simple elitist evolutionary algorithm (SEEA), proposed here for the first time,
have all been fitted for the particular case of the bi-objective ring star problem to be solved.
Experiments were conducted using a set of multiple benchmark test instances.
After having studied the influence of some parameters on the efficiency of the different search methods, 
we compared them to each other and we concluded that the iterative version of IBMOLS was, in general,
significantly better than all the evolutionary algorithms on the problem instances we tackled.
Nevertheless, 
the behavior of a simple search method like SEEA in comparison to the ones of IBEA and NSGA-II is very encouraging with regard to combinatorial problem solving.
As a second step, we designed two cooperative schemes, able to solve any kind of multi-objective optimization problems, 
between SEEA and the non-iterative version of IBMOLS.
Both consist of using SEEA as the main search process and of launching IBMOLS regularly to intensify the search at a given time.
The first hybrid metaheuristic, the periodic cooperative search (PCS), operates a systematic launching of IBMOLS at each of its step.
The second one, the auto-adaptive cooperative search (ACS), launches IBMOLS at a given step of the algorithm only if a certain condition is verified.
Thus, ACS evolves adaptively according to the search scenario and decides by itself, and online, when the cooperation must occur.
In comparison to stand-alone metaheuristics, these two hybrid search methods statistically improves the results on a large number of test instances, 
and particularly on large-size ones.
However, the efficiency difference between PCS and ACS is almost negligible.
This can be explained by the fact that ACS spends time to compute whether the cooperation should occur at a given step of the algorithm. 
This is not the case for PCS, so that the latter devotes all of its computational time to the search process,
even if a cooperation may start at an inopportune time, and then be less gainful.

Although the approaches proposed in this paper are already promising, even compared to the mono-objective approaches of the literature, 
a few future research directions are open.
First is the possibility to improve the population initialization strategy used within all the proposed search methods.
Indeed, the simple one that is used in the paper is a bit rudimentary, as each initial solution has approximately half of its nodes that belong to the cycle.
Surely, an initial population having a well-spread number of visited nodes among its members 
should accelerate the convergence of the search methods toward a well-spread set of non-dominated solutions.
Moreover, the use of a traveling salesman problem heuristic would be helpful in order to improve the ring cost of these initial solutions.
Second, 
we pointed out that the crossover operator designed for the problem under consideration has a tendency to break the parent ring structures in the offspring individuals.
Then, as proposed in~\cite{RBL:04}, we could employ the same heuristic to improve the ring cost of the newly generated solutions.
Third, given that several mutation operators are used within the evolutionary algorithms, 
it is always difficult to determine the appropriate probability selection of each one of them.
Then, a possibility would be to set these rates adaptively, during the search process, according to the efficiency of each operator.
Such a strategy has already been proposed in~\cite{BST:02} in the frame of multi-objective optimization.
Note that the last two points could greatly improve the efficiency of the evolutionary algorithms compared to the local search method, 
what would also be beneficial for the hybrid methods.
From a purely algorithmic point of view, 
it would be interesting to try other quality indicators than the $\epsilon$-indicator within the fitness assignment scheme of IBMOLS.
We know that this binary indicator, although experiencing some difficulties concerning the extreme points of the trade-off surface, 
performs in general much better than others~\cite{BB:07}.
But the indicator-based fitness assignment strategy using the unary hypervolume indicator recently proposed in~\cite{ZBT:07} could be an interesting alternative to explore.
Next, we saw that the stand-alone metaheuristics taking part in the cooperative approach proposed in this paper, SEEA and IBMOLS, can both be launched independently, 
exchanging information via the archive containing the set of non-dominated solutions found so far.
Then, analogously to what has been proposed in COSEARCH~\cite{TB:06} for the mono-objective case, 
these two search agents can be launched in parallel while an adaptive memory stores different kind of information, including the archive.
But, a similar generic model adapted to the case of multi-objective optimization still has to be established.
Finally, we found out that SEEA was a good alternative to state-of-the-art evolutionary algorithms for multi-objective combinatorial optimization 
when a given, relatively small, amount of computational time is allowed.
In the future, we plan to tackle other kinds of combinatorial problems within SEEA to verify if our observations are still valid, 
especially in the case where more than two objectives are involved.
The same remark can be done with regard to the cooperative approaches we proposed.

\section*{Acknowledgment}
A subpart of this work has already been published in~\cite{LJ+:08}.
The authors would like to gratefully acknowledge Matthieu Basseur and Edmund K. Burke for their helpful contribution.

\bibliographystyle{plain}
\bibliography{meta,moo,rsp}

\newpage
\tableofcontents

\end{document}